\documentclass[mathpazo]{cicp}
\usepackage{graphicx}

\usepackage[margin=2.5cm]{geometry}
\usepackage{setspace}
\usepackage[percent]{overpic}
\usepackage{lmodern}
\usepackage{amsmath,amssymb}
\usepackage{amsfonts,amsthm}
\usepackage{lineno}
\usepackage{graphicx}
\usepackage{multirow}
\usepackage{bm, bbm}
\usepackage{colortbl}
\usepackage{booktabs}
\usepackage{fancyvrb}
\usepackage{algorithmicx,algorithm}
\usepackage{algpseudocode}
\usepackage{physics}
\usepackage{xcolor}
\usepackage{pxfonts}
\usepackage[footnotesize,bf]{caption}
\usepackage{subcaption}
\usepackage{mathtools}
\usepackage{hhline}
\usepackage[T1]{fontenc}
\usepackage[hidelinks]{hyperref}
\usepackage{placeins}

\newcommand{\mc}{\mathcal}

\DeclareMathOperator*{\argmin}{arg\,min}

\definecolor{lightergreen}{RGB}{34, 139, 34} %

\graphicspath{{figures/}}

\begin{document}
\title{A hybrid FEM-PINN method for time-dependent partial differential equations}

\author[Author]{Xiaodong Feng\affil{1}\comma\corrauth, Haojiong Shangguan\affil{1}, Tao Tang\affil{2,1}, Xiaoliang Wan\affil{3}}
\address{\affil{1}\ Faculty of Science and Technology, Beijing Normal-Hong Kong Baptist University, Zhuhai 519087, China. \\
	\affil{2}\ School of Mathematics and Statistics, Guangzhou Nanfang College, Guangzhou 510970, China. \\
	\affil{3}\ Department of Mathematics and Center for Computation and Technology, Louisiana State University, Baton Rouge 70803, USA.}
\email{{\tt xiaodongfeng@bnbu.edu.cn} (X. Feng), {\tt shangguanhaojiong@bnbu.edu.cn} (H. Shangguan), {\tt ttang@nfu.edu.cn} (T. Tang), {\tt xlwan@lsu.edu} (X. Wan)}

\begin{abstract}
	In this work, we present a hybrid numerical method for solving evolution partial differential equations (PDEs) by merging the time finite element method with deep neural networks. In contrast to the conventional deep learning-based formulation where the neural network is defined on a spatiotemporal domain, our methodology utilizes finite element basis functions in the time direction where the space-dependent coefficients are defined as the output of a neural network.  We then apply the Galerkin or collocation projection in the time direction to obtain a system of PDEs for the space-dependent coefficients which is approximated in the framework of PINN. The advantages of such a hybrid formulation are twofold: statistical errors are avoided for the integral in the time direction, and the neural network's output can be regarded as a set of reduced spatial basis functions.
	To further alleviate the difficulties from high dimensionality and low regularity, we have developed an adaptive sampling strategy that refines the training set.
	More specifically, we use an explicit density model to approximate the distribution induced by the PDE residual and then augment the training set with new time-dependent random samples given by the learned density model.
	The effectiveness and efficiency of our proposed method have been demonstrated through a series of numerical experiments.
\end{abstract}


\maketitle

\section{Introduction}
Evolution equations, including both time-dependent ordinary and partial differential equations (ODEs / PDEs), are used
in various science and engineering applications, such as fluid mechanics, electromagnetics, quantum mechanics, and elasticity \cite{courant2008methods}.
Many numerical approaches have been developed for such problems, e.g. the finite difference method, the spectral method, and the finite element method.
Recently solving PDEs with deep learning methods has been receiving increasing attention \cite{weinan2021dawning,han2018solving,karniadakis2021physics}. Typical techniques include physics-informed neural networks (PINNs) \cite{raissi2019physics}, the deep Ritz methods \cite{weinan2018deep}, the weak adversarial networks \cite{zang2020weak}, etc.
Although deep learning-based approaches have shown a lot of potential in solving high-dimensional PDEs, there still exist many numerical issues in adapting the neural network approximation to the problem studied.

\subsection{Related work}
In this work, we pay particular attention to the error of PINNs for evolution equations which may grow too fast and limit the application of PINNs in long-term integration.
Many efforts have been made to address this issue. We now briefly review the relevant works.

\paragraph{Improved PINNs:}
To improve the performance of PINNs on long-time integration, many approaches have been developed which mainly focus on seeking a more effective training strategy. A natural idea is to break down the long-time simulation into sub-problems over short time periods. In \cite{meng2020ppinn}, the authors propose a parareal physics-informed neural network (PPINN) that decomposes long-time integration of PDEs into parallelizable short-time subproblems guided by a fast coarse solver, significantly accelerating computation while maintaining accuracy. In \cite{krishnapriyan2021characterizing, wight2020solving}, a time-marching strategy is proposed by splitting the time domain into many small segments, where the training is done segment by segment and the approximate solution at the end of one segment is used as the initial condition for the next segment. In \cite{mattey2022novel}, backward-compatibility PINNs (bc-PINNs) are proposed, where the obtained solution in previous time segments is used as a constraint for the model training in the current time segment. Although these strategies have significantly reduced the difficulty of optimizing long-term simulations, a concrete justification of why such strategies appear to be effective is still missing. In \cite{wang2024respecting}, the authors mention the concept of  the principle of causality and Causal PINNs are developed to incorporate causality into the training process by introducing causality weights. Following this research work, some research efforts \cite{mojgani2023kolmogorov, penwarden2023unified, jung2024ceens, thomas2025causality, wang2026causality} have taken into account the causality. In \cite{mojgani2023kolmogorov}, the authors link PINN training complexity to the Kolmogorov n-width and proposes Lagrangian PINNs (LPINNs) that conform to the causality of convection, reducing training difficulty for convection-dominated PDEs. In \cite{penwarden2023unified}, a unified scalable framework for causal sweeping strategies is developed. In \cite{jung2024ceens}, by using sequential subdomain networks and integral form of PDE, the causality is incorporated into the training process to improve accuracy. In \cite{thomas2025causality}, an effective parametric heat transfer solver is developed by proposing a novel loss function that respects causality via adaptive collocation. In \cite{wang2026causality}, this work enhances PINNs for spatiotemporal PDEs by synergistically combining causality-informed training with residual-based adaptive refinement, achieving superior accuracy in capturing complex interface evolution.

\paragraph{Evolution deep neural networks (EDNNs):}
EDNNs \cite{du2021evolutional, bruna2024neural} are formulated with the Dirac-Frenkel variational principle to train networks by minimizing the residual sequentially in time, where the model parameters are time-dependent rather than global in  the whole space-time domain. Traditional time-marching methods can be used to update the model parameters. Compared to PINNs, the training of EDNNs is more expensive while it is more flexible to adapt EDNNs to constraints such as Hamiltonian conservation \cite{schwerdtner2024nonlinear}. The efficiency of EDNNs is improved in \cite{kao2024petnns} by making a part of model parameters time-dependent and in \cite{berman2024randomized} by updating randomized sparse subsets of model parameters at each time step. Energy-dissipative evolutionary deep operator neural networks \cite{zhang2024energy} adopts the EDNNs paradigm to construct an energy-dissipative operator learning framework that approximates solution operators while strictly preserving energy stability. In \cite{kast2024positional}, EDNNs are extended to parametric PDEs on complex geometries using Laplace-Beltrami embeddings and architectural boundary enforcement, while enhancing performance via Krylov solvers and stiff integrators. In \cite{kim2025multi}, Multi-EDNN methods comprising coupled EDNN (CEDNN) for decoupled multi-variable systems and distributed EDNN (DEDNN) for partitioned domains is provided to efficiently solve complex PDEs with enhanced scalability.

\paragraph{Operator learning:}
Physics-informed DeepONet \cite{wang2023long, wang2021learning} can be used to learn a solution operator over a short time interval $t\in[0,\Delta t]$. Starting with $n=2$, the model's prediction at $n\Delta t$ can be obtained from the trained model using the approximate solution at $(n-1)\Delta t$ as the input. TI-DeepONet \cite{nayak2025ti} embeds numerical time-integrators into operator learning to enforce causality and enable stable. Other operator learning examples include auto-regressive networks \cite{geneva2020modeling}, transformer \cite{geneva2022transformers, wu2025coast}, diffusion model \cite{lippe2023pde} etc. In \cite{gaby2024neural}, the authors propose a novel operator learning framework that approximates solution operators by learning a control vector field in the parameter space of a reduced-order model, enabling efficient evolution of arbitrary initial conditions.

\paragraph{Hybrid strategies:} These approaches try to hybridize classical numerical methods with deep learning techniques by either adapting neural networks to augment classical PDE solvers \cite{bruno2022fc,dresdner2022learning} or adapting classical numerical approaches to improve the performance of PINNs \cite{chiu2022can}. For example, in \cite{chiu2022can}, a coupled automatic and numerical differentiation approach is proposed to take advantage of the regularization induced by numerical discretization.
In \cite{gu2022deep}, a deep adaptive basis Galerkin approach is proposed where the orthogonal polynomial expansion is employed in time direction and the expansion coefficients are modeled as the output of a deep neural network.

\subsection{Our contribution}
The main contributions of this work are summarized as follows:
\begin{itemize}
	\item We have developed a hybrid numerical method by merging the time finite element method with deep neural networks. The approximate solution is a linear combination of the time finite element basis functions, where the coefficients are given by the output of a neural network. We subsequently apply Galerkin or collocation projection to eliminate the time and use the neural network to approximate the coefficients of the basis functions. This strategy has some advantages: First, the numerical difficulties induced by random sampling and causality are avoided in the time direction since the projection can be done accurately.
	      Second, all the coefficients define a set of reduced basis functions on the computation domain, which are learned through the neural network. The approximate solution can also be regarded as a time-dependent linear combination of these reduced basis functions, which shares similarities with the low-rank representation in the study of high-dimensional problems.

	\item We have proposed a deep adaptive sampling strategy to enhance the numerical efficiency. Mesh refinement in the time direction is straightforward. Particular attention needs to be paid to the random sampling in the physical space, especially for high-dimensional and low-regularity problems. Using a spatially conditional bounded KRnet and a discrete distribution in the time direction, we have constructed a joint density model to learn the distribution induced by the PDE residual, based on which new time-dependent samples are generated to refine the training set.
\end{itemize}
The remainder of the paper is organized as follows. In Section \ref{section:preliminaries}, we provide an overview of PINNs. Then we briefly introduce the existing work in solving evolution PDEs and review classic finite element methods for the first-order ordinary differential equations, including Galerkin and collocation frameworks. In Section \ref{section:methodology}, we introduce a hybrid FEM-PINN method
and also develop some adaptive sampling strategies.
Several numerical examples are demonstrated in Section \ref{section:numerical_experiments} to test the performance of the proposed method. The paper is concluded in Section \ref{section:conclusion}.

\section{Preliminaries}\label{section:preliminaries}

\subsection{Physics-informed neural networks (PINNs)}
We begin with a brief overview of physics-informed neural networks (PINNs).
We consider a general time-dependent PDE
\begin{equation}
	u_t(\bm{x},t) + \mathcal{N}[u](\bm{x},t)=f(\bm{x},t),\quad  \bm{x}\in\Omega\subset \mathbb{R}^d, t\in[0,T],
	\label{general_PDE}
\end{equation}
subject to the initial and boundary conditions
\begin{equation}
	\begin{aligned}
		 & u(\bm{x},0) =g(\bm{x}),\quad \bm{x}\in \Omega,                                       \\
		 & \mathcal{B}[u](\bm{x},t) = b(\bm{x}, t),\quad \bm{x}\in\partial \Omega,\; t\in[0,T],
	\end{aligned}
	\label{general_pde_condition}
\end{equation}
where $\mathcal{N}[\cdot]$ is a linear or nonlinear differential operator, and $\mathcal{B}[\cdot]$ is a boundary operator corresponding to Dirichlet, Neumann, Robin or periodic boundary conditions.

Following the original work of Raissi et al. \cite{raissi2019physics}, we represent the unknown solution $u(\bm{x}, t)$ with a deep neural network $u_{\theta}(x, t)$, where $\theta$ denotes all tunable parameters (e.g. weights and biases). Then, a physics-informed model can be trained by minimizing the following composite loss function
\begin{equation}
	\mathcal{L}(\theta) = \lambda_{ic} \mathcal{L}_{ic}(\theta) + \lambda_{bc}\mathcal{L}_{bc}(\theta) + \lambda_{r}\mathcal{L}_{r}(\theta),
\end{equation}
where
\begin{equation}
	\begin{aligned}
		 & \mathcal{L}_{ic}(\theta) = \frac{1}{N_{ic}}\sum_{i=1}^{N_{ic}}\left|u_{\theta}(\bm{x}_{ic}^{i}, 0) - g(\bm{x}_{ic}^{i})\right|^2,                                                                          \\
		 & \mathcal{L}_{bc}(\theta) = \frac{1}{N_{bc}}\sum_{i=1}^{N_{bc}}\left|\mathcal{B}[u_{\theta}]( \bm{x}_{bc}^{i}, t_{bc}^{i}) - b(\bm{x}_{bc}^{i}, t_{bc}^{i})\right|^2,                                       \\
		 & \mathcal{L}_{r}(\theta)=\frac{1}{N_r}\sum_{i=1}^{N_r}\left|\frac{\partial u_{\theta}}{\partial t}(\bm{x}_r^i, t_r^i) + \mathcal{N}[u_{\theta}](\bm{x}_{r}^{i}, t_{r}^{i}) - f(\bm{x}_r^i, t_r^i)\right|^2.
	\end{aligned}
\end{equation}
Here $\{\bm{x}_{ic}^{i}\}_{i=1}^{N_{ic}}$, $\{t_{bc}^{i},\bm{x}_{bc}^{i}\}_{i=1}^{N_{bc}}$ and $\{t_{r}^{i}, \bm{x}_{r}^{i}\}_{i=1}^{N_r}$ can be the vertices of a fixed mesh or points that are randomly sampled
at each iteration of a gradient descent algorithm. The gradients with respect to both the input variables $(t,\bm{x})$ and the model parameters $\theta$ can be efficiently computed via automatic differentiation \cite{griewank2008evaluating}. Moreover, the hyper-parameters $\{\lambda_{ic}, \lambda_{bc}, \lambda_{r}\}$ allow the flexibility of assigning a varying learning rate to each loss term to balance their interplay during the training process, which may be user-specified or tuned automatically.

\subsection{Continuous time finite element method}
Evolution PDEs are often approximated by spatial finite elements together with a sequential time-marching scheme. Another choice is to construct a finite element approximation space on the space-time domain. We briefly review the time finite element method for first-order ordinary differential equations. Consider the following model problem:
\begin{equation}
	\begin{aligned}
		 & u'(t)  + \mathcal{N}[u](t) = f(t), \quad t\in [0,T], \\
		 & u(0) = 0,
	\end{aligned}
	\label{model_equation}
\end{equation}
where $\mathcal{N}$ is a linear or nonlinear operator and $f(t)\in L^2(I)$ with $I=[0,T]$.
\subsubsection{Galerkin projection}
We let $X \coloneqq \{u\in H^{1}(I):u(0)=0\}$ be the trial space
and $Y \coloneqq L^2(I)$ the test space.
The primal variational formulation of \eqref{model_equation} is as follows.
\begin{equation}
	\left\{
	\begin{aligned}
		 & \mbox{Find}\; u \in X \; \mbox{such that}                     \\
		 & (u', v) + (\mathcal{N}[u], v) = (f, v),\quad \forall v \in Y,
	\end{aligned}
	\label{variational_form}
	\right.
\end{equation}
where $(\cdot,\cdot)$ indicates the inner product of two functions. For approximation, we consider the Galerkin projection with $X_{N}\subset X$ and $Y_{N}\subset Y$, i.e.,
\begin{equation}
	\left\{
	\begin{aligned}
		 & \mbox{Find } u_N \in Y_{N} \mbox{ such that }                                     \\
		 & (u_N' , v_N)  + (\mathcal{N}[u_N], v_N) = (f, v_N) ,\quad \forall v_N\,\in Y_{N}.
	\end{aligned}
	\right.
	\label{galerkin_form}
\end{equation}
Define
\begin{equation}
	X_N = \mathrm{span}\{\phi_j(t)|0\leq j\leq N\}, \quad Y_N =  \mathrm{span}\{\psi_j(t)|0\leq j\leq N\},
\end{equation}
where $\phi_j(t)$ and $\psi_j(t)$ are finite element basis functions. Let
\begin{equation}
	u_N(t) = \sum_{i=0}^{N} \tilde{u}_{i}\phi_i(t) \in X_N
	\label{model_galerkin_expression}
\end{equation}
be the approximate solution with undetermined coefficients $\tilde{u}_i$. Taking $v_N=\psi_j$ for $j=0,\cdots,N$ in (\ref{galerkin_form}) leads to the following system
\begin{equation}
	\sum_{i=0}^{N}\big(\partial_t \phi_i(t),\psi_j(t)\big) \tilde{u}_i+ \bigg( \mc{N}\Big[\sum_{i=0}^{N}\phi_i(t)\tilde{u}_{i}\Big],\psi_j(t) \bigg)=(f(t),\psi_j(t)), \quad \forall j=0,1, \cdots, N.
	\label{galerkin_system}
\end{equation}
The inner products in \eqref{galerkin_system} can be accurately evaluated with the Gaussian quadrature formulas, where the degree of exactness is determined by the nonlinearity of $\mathcal{N}$.
\subsubsection{Collocation projection}\label{finite_element_collocation_method}
Collocation projection provides a flexible strategy especially when $\mathcal{N}$ is nonlinear \cite{douglas1973finite,percell1980c}.
Let $\{s_k\}_{k=1}^K$ be Gaussian-type quadrature points on the reference interval $[0,1]$
\begin{equation}
	0\leq s_1<s_2<\cdots<s_K\leq 1.
\end{equation}
Consider a partition of $[0,T]$ with
\begin{equation}
	0=t_0<t_1<\cdots <t_{\hat{M}}=T, \quad h_i = t_i - t_{i-1},\;  i=1,\cdots,\hat{M}.
\end{equation}
Define
\begin{equation}\label{eqn:grid_s}
	s_{m,k} = t_{m-1} + h_ms_k,\quad 1\leq k \leq K,\quad 1\leq m\leq \hat{M}.
\end{equation}
We seek the approximate solution by enforcing the equation on the collocation points, i.e.,
\begin{equation}
	\left\{
	\begin{aligned}
		 & \mbox{Find } u\in X_N\cap C^1(I)  \mbox { such that}                                            \\
		 & \partial _t u(s) + \mathcal{N}[u](s) = f(s), \quad s\in\cup_{m=1}^{\hat{M}}\{s_{m,k}\}_{k=1}^K,
	\end{aligned}
	\right.
\end{equation}
where $X_N\cap C^1(I)$ defines a finite element approximation space with $C^1$ elements. It is shown in \cite{Boor1973} by selecting the collocation points carefully collocation projection yields the same order of accuracy as Galerkin projection.  Typical piecewise polynomials with at least $C^1$ regularity include piecewise cubic Hermite polynomials and cubic spline functions.

\section{Methodology} \label{section:methodology}
Now we are ready to present our approach for evolution equations (\ref{general_PDE}).
We aim to seek an approximate solution of the following form
\begin{equation}
	u_N(\bm{x},t;\theta) = \sum_{i=0}^{N}\omega_i(\bm{x};\theta)\phi_i(t),
	\label{neural_galerkin_formulation}
\end{equation}
where $\{\phi_i(t)\}_i$ is a pre-specified set of time finite element basis functions, $\omega_i: \mathbb{R}^d\to \mathbb{R}$ are modeled by the output of a neural network $\bm{\omega}(\bm{x},\theta):\mathbb{R}^d\rightarrow\mathbb{R}^{N+1}$, and $\theta$ includes all tunable model parameters.
More precisely, $\omega(\bm{x},\theta)$ is a fully-connected neural network defined as
\begin{equation}
	\bm{\omega}(\bm{x};\theta) \coloneqq a^T h_{L-1} \circ h_{L-2}\circ \cdots \circ h_1(\bm{x}),\quad \forall\bm{x}\in\mathbb{R}^d,
\end{equation}
where $L\in\mathbb{N}^+$, $a\in \mathbb{R}^{M_{L-1}\times (N+1)}, \; h_{\ell}(\bm{x}_{\ell})\coloneqq \sigma(W_{\ell}\bm{x}_{\ell}+ b_{\ell})$ with $W_{\ell}\in \mathbb{R}^{M_{\ell}\times M_{\ell-1}}\; (M_0\coloneqq d)$ and $b_{\ell}\in \mathbb{R}^{M_\ell}$ for $\ell=1,2,\cdots,L-1$. Then $\theta \coloneqq \{a, W_{\ell}, b_{\ell}:1\leq \ell \leq L-1 \}$.
$\sigma(\bm{x})$ is an activation function which acts on $\bm{x}$ componentwisely to return a vector of the same size as $\bm{x}$. We let $M_{\ell}=M$ be a fixed number for all $\ell$ and $\mathcal{F}_{L,M}$ the set consisting of all $\bm{\omega}$ with depth $L$ and width $M$.
\subsection{A hybrid FEM-PINN method}
We consider the following hypothesis space
\begin{equation}
	\mathcal{U}_{N}\coloneqq \left\{
	u_{N}(\bm{x},t;\theta) = \sum_{i=0}^{N}\omega_i(\bm{x};\theta)\phi_i(t), \,\bm{\omega}=(\omega_0,\cdots,\omega_N)\in \mathcal{F}_{L,M}
	\right\}.
\end{equation}
The Galerkin projection along the time direction yields that
\begin{equation}
	\left\{
	\begin{array}{l}
		\mathrm{Find} \; u_{N}\in \mathcal{U}_{N} \mbox{ such that} \\
		(\partial_t u_{N}(\bm{x},\cdot), v_N) + (\mathcal{N}[u_{N}](\bm{x},\cdot), v_N) = (f(\bm{x},\cdot), v_N), \quad \forall v_N\in \mathrm{span}\{\psi_j(t)|0\leq j\leq N\}, \quad \forall\,\bm{x}\in \Omega,
	\end{array}
	\right.
\end{equation}
where $(\cdot,\cdot)$ indicates the inner product of two functions with respect to time. More specifically, we have
\begin{equation}
	\sum_{i=0}^{N}\big(\partial_t \phi_i(t),\psi_j(t)\big) \omega_i(\bm{x};\theta) + \bigg( \mc{N}\Big[\sum_{i=0}^{N}\phi_i(t)\omega_i(\bm{x};\theta)\Big],\psi_j(t) \bigg)=(f(\bm{x},t),\psi_j(t)), \quad \forall j=0,1, \cdots, N,\quad \forall \bm{x}\in\Omega.
	\label{galerkin_system_neural_linear}
\end{equation}
If $\mathcal{N}$ is linear with respect to $\omega$ and time-independent, the above system can be further simplified as follows:
\begin{equation}
	\sum_{i=0}^{N}\big(\partial_t \phi_i(t),\psi_j(t)\big) {\omega}_i(\bm{x};\theta)+
	\sum_{i=0}^{N}\big(\phi_i(t),\psi_j(t)\big)\mathcal{N}({\omega}_i(\bm{x};\theta))
	=(f(\bm{x}, t),\psi_j(t)), \quad \forall j=0,1, \cdots, N.
	\label{NN_linear_finite_element_martix_form}
\end{equation}

Now let us turn to the collocation projection along the time direction.
Let $S_t=\cup_{m=1}^{\hat{M}}\{s_{m,k}\}_{k=1}^K$, where $s_{m,k}$ is defined in equation \eqref{eqn:grid_s}. The collocation formulation can be written as
\begin{equation}\label{eqn:sys_collocation}
	\left\{
	\begin{aligned}
		 & \mbox{Find } u_N \in \mathcal{U}_N \mbox{ such that }                                                                                    \\
		 & \partial_t u_N(\bm{x},s;\theta) + \mathcal{N}[u_N](\bm{x},s;\theta) = f(\bm{x}, s),\quad\forall s\in S_t, \quad\forall \bm{x}\in \Omega.
	\end{aligned}
	\right.
\end{equation}
More specifically,
\begin{equation}
	\sum_{i=0}^{N}\partial_t \phi_i(s) \omega_i(\bm{x};\theta) + \mathcal{N}\left[\sum_{i=0}^{N}\phi_i(s)\omega_i(\bm{x};\theta)\right] = f(\bm{x},s), \quad \forall s\in S_t,\quad\forall \bm{x}\in\Omega.
	\label{collocation_system_neural_linear}
\end{equation}

\begin{remark}
	The Galerkin and collocation projections yield respectively two systems of PDEs for $\omega_i(\bm{x};\theta)$. Due to the hybrid form of $u_N(\bm{x},t;\theta)$, all integrals for the Galerkin projection in the time direction can be done accurately by Gaussian quadrature formulas. Since the temporal basis functions are polynomials, collocation projection is also effective \cite{Boor1973}. We then mainly focus on the integration in the physical space.
\end{remark}

We subsequently approximate the PDE systems \eqref{galerkin_system_neural_linear} and \eqref{collocation_system_neural_linear} in the framework of PINNs.
More specifically, we consider the following minimization problem
\begin{equation}
	\min\limits_{\theta} \mathcal{L}(\theta) = \mathcal{L}_{r}(\theta) + \gamma_1 \mathcal{L}_{ic}(\theta) + \gamma_2 \mathcal{L}_{bc}(\theta),
	\label{NN_continuos_loss}
\end{equation}
where
\begin{equation}
	\begin{aligned}
		 & \mathcal{L}_{ic}(\theta) = \Vert u(\bm{x},0;\theta) - g(\bm{x})\Vert _{L^2(\Omega)}^2,\quad \mathcal{L}_{bc}(\theta) = \Vert\mathcal{B}[u](\bm{x},t;\theta) - b(\bm{x},t)\Vert _{L^2(\partial \Omega\times [0,T])}^2,
	\end{aligned}
	\label{individual_loss}
\end{equation}
with $0<\gamma_1,\gamma_2<\infty$ being penalty parameters. For system \eqref{galerkin_system_neural_linear}, we define
\begin{equation}
	\begin{aligned}
		 & \mathcal{L}_r(\theta) = \sum_{j=0}^{N} \mathcal{L}^g_{r,j}(\theta),                                                                                                                                                                                                                                                        \\
		 & \mathcal{L}^g_{r,j}(\theta) = \left\Vert \sum_{i=0}^{N}\big(\partial_t \phi_i(t),\psi_j(t)\big) {\omega}_i(\bm{x};\theta)+ \bigg( \mc{N}\Big[\sum_{i=0}^{N}\phi_i(t)\omega_i(\bm{x};\theta)\Big],\psi_j(t) \bigg)-(f(\bm{x}, t),\psi_j(t))\right\Vert^2_{L^2(\Omega)}=\left\|r_j^g(\bm{x};\theta)\right\|^2_{L^2(\Omega)}.
	\end{aligned}
	\label{galerkin_continuous_residual}
\end{equation}
For system \eqref{collocation_system_neural_linear}, we define
\begin{equation}
	\begin{aligned}
		 & \mathcal{L}_r(\theta) = \sum_{j=1}^{|S_t|}\mathcal{L}^c_{r,j}(\theta),                                                                                                                                                                                                                \\
		 & \mathcal{L}^c_{r,j}(\theta) = \left\Vert \sum_{i=0}^{N} \partial_t \phi_i(s_j){\omega}_i(\bm{x};\theta) + \mathcal{N}\left[\sum_{i=0}^{N}{\omega_i}(\bm{x};\theta)\phi_i(s_j)\right] - f(\bm{x}, s_j)\right\Vert^2_{L^2(\Omega)}=\left\|r_j^c(\bm{x};\theta)\right\|^2_{L^2(\Omega)},
	\end{aligned}
\end{equation}
where we order all collocation points in $S_t$ as $s_j$ with $j=1,\ldots,|S_t|$.

\begin{figure}[H]
	\centering
	\includegraphics[width=0.7\linewidth]{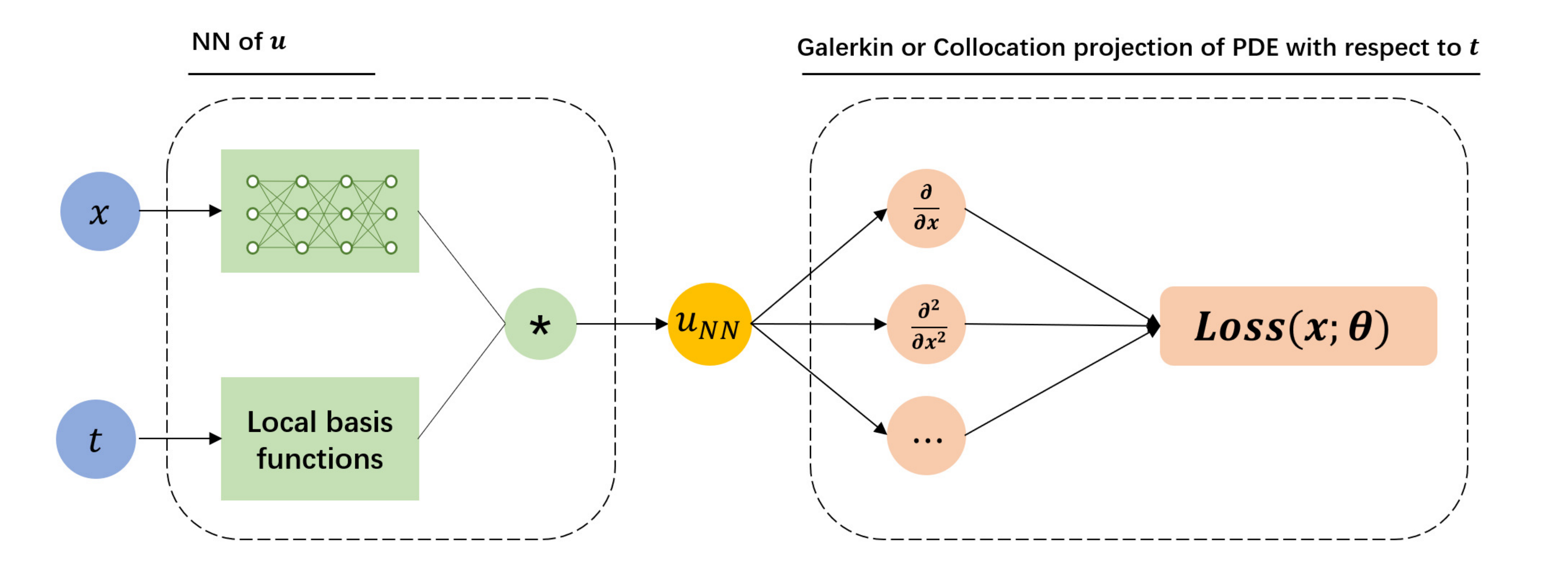}
	\caption{Schematic of the proposed approach.}
	\label{fig:schematic}
\end{figure}

We note that $|S_t|\geq (N+1)$ in general since there are $N+1$ time finite element basis functions. For simplicity, we let $|S_t|=N+1$ and consider the following form
\begin{equation}
	\mathcal{L}_r(\theta) = \sum_{j=0}^{N}\mathcal{L}_{r,j}(\theta),\quad \mathcal{L}_{r,j}(\theta) = \left\|r_j(\bm{x};\theta)\right\|^2_{L^2(\Omega)},
	\label{def_r_i}
\end{equation}
which are shared by both the Galerkin and collocation projections, i.e., $r_j=r_j^g\textrm{ or }r_j^c$.
The loss functional (\ref{NN_continuos_loss}) is usually discretized numerically before the optimization with respect to $\theta$ is addressed. In practice, one often chooses uniformly distributed collocation points $S_r = \{S_{r,j}\}_{j=0}^{N}=\left\{\{\bm{x}_{r,j}^{(i)}\}_{i=1}^{N_{r,j}}\right\}_{j=0}^{N}$, $S_{ic}=\{\bm{x}_{ic}^{(i)}\}_{i=1}^{N_{ic}}$ on $\Omega$ and $S_{bc} = \left\{(\bm{x}_{bc}^{(i)}, t_{bc}^{(i)})\right\}_{i=1}^{N_{bc}}$ on $\partial \Omega\times[0,T]$ for the discretization of the three terms in the objective functional (\ref{NN_continuos_loss}), leading to the following empirical loss
\begin{equation}
	\widehat{\mathcal{L}}(\theta) = \sum_{j=0}^{N}\Vert r_j(\bm{x};\theta)\Vert^2_{N_{r,j}} + \hat{\gamma}_1 \Vert \mathcal{B}[u](\bm{x},t;\theta) - b(\bm{x},t)\Vert _{N_{ic}} ^2 + \hat{\gamma}_2 \Vert u(\bm{x},0;\theta) - g(\bm{x})\Vert^2_{N_{bc}},
	\label{eqn:empricial_loss}
\end{equation}
where $0<\hat{\gamma}_1, \hat{\gamma}_2<\infty$, and
\begin{equation}
	\Vert u(\bm{x}) \Vert_{N_{r,j}} = \left(\frac{1}{N_{r,j}}\sum_{i=1}^{N_{r,j}}u^2(\bm{x}_{r,j}^{(i)})\right)^{\frac{1}{2}},\,\Vert u(\bm{x},0)\Vert _{N_{ic}} = \left(\frac{1}{N_{ic}}\sum_{i=1}^{N_{ic}}u^2(\bm{x}_{ic}^{(i)},0)\right)^{\frac{1}{2}},\,\Vert u(\bm{x},t)\Vert_{N_{bc}} = \left(\frac{1}{N_{bc}}\sum_{i=1}^{N_{bc}}u^2(\bm{x}_{bc}^{(i)}, t_{bc}^{(i)})\right)^{\frac{1}{2}}.
	\label{empirical_loss}
\end{equation}
We then seek an estimator $\hat{\theta}$ by minimizing the empirical loss (\ref{eqn:empricial_loss}) via stochastic gradient descent methods, i.e.,
\begin{equation}\label{eqn:opt_discrete}
	\hat{\theta}=\argmin_{\theta} \widehat{\mathcal{L}}(\theta).
\end{equation}
As suggested by \cite{wight2020solving,krishnapriyan2021characterizing,mattey2022novel}, we can also employ a time-marching strategy to reduce optimization difficulties. Specifically, we partition the temporal domain $[0,T]$ into sub-domains $[0,\Delta t], [\Delta t, 2\Delta t],\cdots, [T-\Delta t, T]$.
Neural networks are then trained on each sub-domain successively with the initial conditions given by the same model trained on previous sub-domains.
The schematic of the proposed approach is shown in Figure \ref{fig:schematic}, and the corresponding algorithm is summarized as follows.
\begin{algorithm}
	\caption{Solving time-dependent PDEs by hybrid FEM-PINN method}
	\label{alg:1}
	\begin{algorithmic}
		\State \textbf{Input:} terminal time $T$, number of segment $N_{\mathrm{segment}}$, number of basis function $N + 1$, number of epoch $N_e$, training data $S_r = \cup_{j=0}^NS_{r,j}=\cup_{j=0}^N\{\bm{x}_{r,j}^{(i)}\}_{i=1}^{N_{r,j}}$, $S_{ic}=\{\bm{x}_{ic}^{(i)}\}_{i=1}^{N_{ic}}$ and $S_{bc} = \{(\bm{x}_{bc}^{(i)}, t_{bc}^{(i)})\}_{i=1}^{N_{bc}}$, initial learning rate $l_r^{(0)}$, decay rate $\eta$, step size $n_s$
		\For{$k=1,\cdots,N_{\mathrm{segment}}$}
		\State Construct $N+1$ local finite element basis functions in the interval $[(k-1)T/N_{\mathrm{segment}}, kT/N_{\mathrm{segment}}]$.
		\State Set $l_r$ to $l_r^{(0)}$.
		\For{$j=1,\cdots,N_e$}

		\If{$j \% n_s==0$}
		\State $l_r=\eta*l_r$.
		\EndIf
		\State Divide $\{S_{r,j}\}_{j=0}^{N}$, $S_{ic}$ and $S_{bc}$ into $n_r, n_{ic}, n_{bc}$ mini-batches  $\{S_{r,j}^{ib}\}_{j=0}^{N}, S_{ic}^{ib}, S_{bc}^{ib}$ randomly, respectively.
		\For{$ib = 1,\cdots,n$}
		\State Compute the loss function ${\widehat{\mathcal{L}}}(\theta)$ \eqref{eqn:empricial_loss} for mini-batch data  $\{S_{r,j}^{ib}\}_{j=0}^{N}, S_{ic}^{ib}$ and  $S_{bc}^{ib}$.
		\State Update $\theta$ by stochastic gradient descent method.
		\EndFor
		\EndFor
		\State Compute the prediction for new initial training data set $\{\bm{x}_{ic}^{(i)}, t_{ic}^{(i)}\}_{i=1}^{N_{ic}}$, where $t_{ic}^{(i)} = kT/N_{\mathrm{segment}}$.
		\State Update the above prediction values as $g(\bm{x}_{ic}^{(i)})$.
		\State Update the corresponding initial and boundary training data $S_{ic}, S_{bc}$.
		\State Save the parameter $\theta$ as $\theta_k$.
		\EndFor
		\State \textbf{Output:} predicted solutions $\{u(\bm{x},t_k;\theta_k)\}_{k=1}^{N_{\mathrm{segment}}}$, where $t_k\in [(k-1)T/N_{\mathrm{segment}}, kT/N_{\mathrm{segment}}], k=1,\cdots,N_{\mathrm{segment}}$.
	\end{algorithmic}
\end{algorithm}

\subsection{Error analysis}
We conduct error analysis for the proposed hybrid FEM-PINN method in the case of linear time-independent operator $\mathcal{N}$ and Galerkin projection. The error analysis for the collocation projection and nonlinear operator $\mathcal{N}$ will be left for future work.

Let \( u_N(\bm{x}, t) \) denote the temporal finite element approximation of the solution with $N$ finite element basis functions, \( u_{N,M}(\bm{x}, t;\theta^*) \) the best approximation of \( u_N(\bm{x}, t) \) by a two-layer fully-connected neural network of width \( M \), \( u_{N,M}(\bm{x}, t;\theta_P^*) \) the best approximation obtained by minimizing the empirical loss over \( P \) training samples and \( u_{N,M}(\bm{x}, t;\theta_P) \) the actual neural network obtained by minimizing the empirical loss. Then, the total error can be decomposed as follows:
\begin{equation}
	\begin{aligned}
		     & \| u(\bm{x}, t) - u_{\mathrm{N,M}}(\bm{x}, t;\theta_P)\|                                                                                                                                                              \\
		\leq & \underbrace{\| u(\bm{x}, t) - u_{N}(\bm{x}, t)\| }_{\mbox{Finite element approximation error}} + \underbrace{ \|u_N(\bm{x}, t) - u_{N,M}(\bm{x}, t;\theta^*)\|}_{\mbox{Neural network approximation error}}           \\
		     & + \underbrace{\|u_{N,M}(\bm{x}, t;\theta^*) - u_{N,M}(\bm{x}, t;\theta_P^*)\|}_{\mbox{Statistical error}} + \underbrace{\|u_{N,M}(\bm{x}, t;\theta_P^*) - u_{N,M}(\bm{x}, t;\theta_P)\|}_{\mbox{Optimization error}},
	\end{aligned}
	\label{eqn:error_decomposition}
\end{equation}
The first term $\|u(\bm{x}, t) -u_N(\bm{x}, t)\|$ quantifies the discretization error due to the use of a finite-dimensional function space (e.g., piecewise polynomial basis) in time for approximating the true solution $u(\bm{x}, t)$.  The second term $\|u_N(\bm{x}, t) - u_{N,M}(\bm{x}, t)\|$ representation error of the neural network in approximating the projection $u_N(\bm{x}, t)$, i.e., how well the neural network with optimized parameters $\theta^*$ can represent the finite element solution. The third term $\|u_{N,M}(\bm{x}, t;\theta^*) - u_{N,M}(\bm{x}, t;\theta)\|$ captures the statistical error, which quantifies the generalization gap due to finite training data, assuming the neural network has been perfectly optimized with no optimization error. The last term $\|u_{N,M}(\bm{x}, t;\theta_P^*) - u_{N,M}(\bm{x}, t;\theta_P)\|$ measures the optimization error, which quantifies the gap between the best possible neural network approximation $u_{N,M}(\bm{x}, t;\theta_P^*)$ and the actual neural network obtained by minimizing the empirical loss $u_{N,M}(\bm{x}, t;\theta_P)$. In the following, we will focus on the analysis of the first three terms, since the optimization error is generally difficult to analyze due to the non-convexity of the high-dimensional loss landscape and the use of stochastic optimization algorithms.

\subsubsection{Barron-type functions and approximation property}
In the error analysis, we consider the approximation property of two layer fully-connected neural networks for Barron-type functions, since there are many existing results showing that dimension-independent approximation rates with respect to the width can be achieved for such functions. Given the following class of two-layer fully-connected neural networks with width $M$:
\begin{equation}
	\mathcal{F}_{\sigma}(M) \coloneqq \left\{\mathrm{NN}(\bm{x}; \theta) = \sum_{m=1}^Mc_m\sigma(\bm{w}_m^Tx+b_m) : c_i, b_i \in \mathbb{R}, \bm{w}_i \in \mathbb{R}^d, m=1,\cdots,M\right\},
\end{equation}
where $\theta = \{c_m, \bm{w}_m, b_m\}_{m=1}^M$ denotes the trainable parameters and $\sigma$ is a non-linear activation function. Then, we say a function $f:\mathbb{R}^d\to\mathbb{R}$ belongs to the Barron space $\mathcal{B}^s$ if it has bounded s-th order Barron norm:
\begin{equation}
	\Vert f \Vert_{\mathcal{B}^s} = \int_{\mathbb{R}^d}(1+\Vert\bm{w}\Vert_1)^s\vert\mathcal{F}[f](\bm{w})\vert d\bm{w},
\end{equation}
where $\mathcal{F}[f](\bm{w})$ is the Fourier transform of $f$.
The following lemma shows differentiable functions with compact support and sufficiently high order derivatives are in $\mathcal{B}^s$.
\begin{lemma}[Lemma 3.4, \cite{gu2022deep}]\label{lemma:Barron_space}
	Let $s \geq 1$ and $f\in C_c^\gamma(\mathbb{R}^d)$ with $\gamma \geq d/2 + s$, then $f \in \mathcal{B}^s$ in the sense that
	\begin{equation}
		\Vert f \Vert_{\mathcal{B}^s} \le C(\gamma, s)\int_{\mathbb{R}^d}\vert f\vert^2 + \vert D^{\gamma}f\vert^2 \mathrm{d}\bm{x} < \infty.
	\end{equation}
\end{lemma}
Then we introduce that any function in the Barron space $\mathcal{B}^s$ can be approximated by a two-layer fully-connected neural network in Sobolev norm.
\begin{lemma}[Corollary 1, \cite{siegel2020approximation}]\label{lemma:NN_approximation}
	Let $s\in\mathbb{N}$. Suppose the activation function of $\mathcal{F}_{\sigma}(M)$ is set as the sigmoidal function $\sigma(x) = (1+e^{-x})^{-1}$, arctan $\sigma(x) = \arctan(x)$, hyperbolic tangent $\sigma(x) = \tanh(x)$ or SoftPlus $\sigma(x) = \log(1 + e^x)$. Then for any $f \in \mathcal{B}^{s+1}$, we have
	\begin{equation}
		\inf_{w\in\mathcal{F}_{\sigma}(M)}\Vert w-f \Vert_{H^s(\Omega)} \le \vert \Omega\vert^{\frac{1}{2}}C(s, diam(\Omega), \sigma)\Vert f\Vert_{\mathcal{B}^{s+1}}M^{-\frac{1}{2}}.
	\end{equation}
\end{lemma}

\subsubsection{Path norm and Rademacher complexity}
In the context of statistical learning theory, Rademacher complexity provides a powerful tool for deriving generalization bounds and it characterizes the worst-case estimation error by controlling how much the empirical risk deviates from the true risk. Consequently, bounding the Rademacher complexity of a hypothesis space allows us to establish data-dependent upper bounds on the generalization gap, thereby rigorously quantifying the trade-off between model complexity and statistical accuracy.

Given a hypothesis space $\mathcal{H}$ and a sample set $S = \{\bm{x}_p\}_{p=1}^P$ drawn i.i.d. from a distribution $\mathcal{D}\subset\mathbb{R}^d$, the empirical Rademacher complexity of $\mathcal{H}$ with respect to $S$ is defined as
\begin{equation}
	\widehat{\mathcal{R}}_S(\mathcal{H}) \coloneqq \mathbb{E}_{\bm{\zeta}}\left[\sup_{h\in\mathcal{H}}\frac{1}{p}\sum_{p=1}^P \zeta_p h(\bm{x}_p)\right],
\end{equation}
where $\bm{\zeta} = (\zeta_1, \cdots, \zeta_P)$ are i.i.d. Rademacher random variables taking values in $\{-1, 1\}$ with equal probability. The following lemma provides a generalization bound in terms of the empirical Rademacher complexity.
\begin{lemma}[Theorem 3.3, \cite{mohri2018foundations}]\label{lemma:generalization_bound}
	Let $\mathcal{H}$ be a hypothesis space of functions mapping from $\mathcal{X}$ to $[0,B^2]$. For any $0 < \delta < 1$, with probability at least $1-\delta$ over the choice of an i.i.d. sample set $S = \{\bm{x}_p\}_{p=1}^P$ of size $n$ drawn from $\mathcal{D}$, the following inequality holds for all $h\in\mathcal{H}$:
	\begin{equation}
		\sup_{h\in\mathcal{H}}\vert\mathbb{E}_{\bm{x}\sim\mathcal{D}}[h(\bm{x})] - \frac{1}{P}\sum_{p=1}^P h(\bm{x}_p)\vert \leq 2\widehat{\mathcal{R}}_S(\mathcal{H}) + 3B\sqrt{\frac{\log(2/\delta)}{2P}}.
	\end{equation}
\end{lemma}

In \cite{li2020complexity}, the authors introduced the path norm for two layer fully-connected neural networks with general activation functions and established upper bounds on the empirical Rademacher complexity of the hypothesis space defined by such networks in terms of the path norm. Consider the following neural network space
\begin{equation}
	\mathcal{W}_Q(\mathcal{F}_{\sigma}(M)) \coloneqq \{x\mapsto \omega(x;\theta):\omega(x;\theta) \in \mathcal{F}_{\sigma}(M), \Vert \theta\Vert_{\mathcal{P}} \le Q\},
\end{equation}
where the path norm are defined as follows:
\begin{equation}
	\Vert\theta\Vert_{\mathcal{P}} \coloneqq \sum_{m=1}^M\vert c_m\vert (\Vert\bm{w}_m\Vert_1 + \vert b_m\vert + 1), \quad \theta\in\Theta \coloneqq \{\theta:\Vert \theta\Vert_{\mathcal{P}} \le Q\}.
\end{equation}
The following lemma provides an upper bound on the empirical Rademacher complexity of $\mathcal{W}_Q(\mathcal{F}_{\sigma}(M))$.
\begin{lemma}[Theorem 3, \cite{li2020complexity}]\label{lemma:Rademacher_complexity}
	Consider the general activation function $\sigma: \mathbb{R} \to \mathbb{R}$. Assume that $\sigma(x)$ is continuous and twice weakly differentiable on $\mathbb{R}$, and $\sigma''(x)$ is Riemann integrable on any compact subset of $\mathbb{R}$. Define
	\begin{align}
		\gamma_0(\sigma) & = \int_{\mathbb{R}}\vert \sigma''(x)\vert(\vert x \vert + 1)\mathrm{d}x, \\
		g(x)             & = \vert \sigma(x)\vert + (\vert x\vert + 2)\vert \sigma'(x)\vert,        \\
		\gamma(\sigma)   & = \gamma_0(\sigma) + \inf_{x\in\mathbb{R}}g(x).
	\end{align}
	If $\gamma(\sigma) < \infty$, then we have
	\begin{equation}
		\widehat{\mathcal{R}}_S(\mathcal{W}_Q(\mathcal{F}_{\sigma}(M))) \leq 2\gamma(\sigma)Q\sqrt{\frac{2\ln(2d+2)}{P}}.
	\end{equation}
\end{lemma}

\subsubsection{Error estimates for heat equation}
For simplicity, we consider the error analysis for the following heat equation with homogeneous Dirichlet boundary condition and initial condition, which is a special case of the general PDE system \eqref{general_PDE} and \eqref{general_pde_condition} with $\mathcal{N} = -\Delta$ and $g = u_0, b = 0$:
\begin{equation}
	\left\{
	\begin{aligned}
		 & \partial_t u(\bm{x},t) -\Delta u(\bm{x},t) = f(\bm{x},t), & (\bm{x},t)\in \Omega\times (0,T],            \\
		 & u(\bm{x},0) = u_0(\bm{x}),                                & \bm{x}\in \Omega,                            \\
		 & u(\bm{x}, t) = 0,                                         & (\bm{x}, t)\in \partial \Omega \times (0,T], \\
	\end{aligned}
	\right.
	\label{heat_equation}
\end{equation}
where $\Omega\subset \mathbb{R}^d$ be a bounded domain with Lipschitz boundary $\partial \Omega$, and let $T>0$ be a fixed final time.

For convenience, we incorporate the initial and boundary conditions into the hypothesis space $\mathcal{U}_N$ by representing the solution as
\begin{equation}
	u_{N, M}(\bm{x}, t;\theta) = u_0(\bm{x})\phi_0(t) + \sum_{j=1}^N \nu(\bm{x})\omega_j(\bm{x}; \theta)\phi_j(t),
	\label{eqn:solution_embedding}
\end{equation}
where $\omega_j(\bm{x}; \theta) \in \mathcal{F}_{\sigma}(M), j=1,\cdots,N$ and $\nu(\bm{x})$ is a prescribed smooth function governing the boundary condition, namely, $\nu|_{\partial \Omega} = 0$. We suppose $u/\nu < \infty$ on $\partial\Omega$ for each $t\in[0, T]$ and normalize $\nu$ by multiplying some positive constant so that
\begin{equation}
	\max\{\Vert \nu\Vert_{\infty}, \Vert D\nu \Vert_{\infty}, \Vert D^2\nu\Vert_{\infty} \} = 1.
\end{equation}

\begin{theorem}\label{theorem:error_estimate_heat_equation}
	For the heat equation \eqref{heat_equation}, we assume that $f\in L^2(0,T;L^2(\Omega))$ and $u_0\in H_0^1(\Omega)\cap H^2(\Omega)$. Let $u(\bm{x}, t)\in  H^2(0,T;L^2(\Omega)) \cap H^1(0,T;H_0^1(\Omega)\cap H^2(\Omega))$ be the unique strong solution and $u_{N,M}(\bm{x}, t;\theta_P^*)$ be the best approximation obtained by minimizing the empirical loss over $P$ training samples $S = \{\bm{x}_p\}_{p=1}^P$ drawn i.i.d. from the uniform distribution on $\Omega$. If the time finite element space $\mathcal{V}_N$ consists of piecewise linear functions and the activation function of $\mathcal{F}_{\sigma}(M)$ is set as hyperbolic tangent $\sigma(x) = \tanh(x)$, then for any $0 < \delta < 1$, with probability at least $1-\delta$, we have
	\begin{equation}
		\begin{aligned}
			     & \| u(\bm{x}, t) - u_{N,M}(\bm{x}, t;\theta_P^*)\|_{L^2(0,T;L^2(\Omega))}                                                                        \\
			\leq & C \left(N^{-1} + N d^{\frac{1}{2}} M^{-\frac{1}{2}} + \sqrt{N^2 d M^{-1} + N^2d \sqrt{\frac{\ln d}{P}} + \sqrt{\frac{\ln(2/\delta)}{P}}}\right)
		\end{aligned}
	\end{equation}
	where $C>0$ is a constant depending on $\Omega, T, \sigma, u, \nu$ and independent of $N, M, P, d$.
\end{theorem}
\begin{proof}
	The proof of Theorem \ref{theorem:error_estimate_heat_equation} is given in Appendix A.
\end{proof}
Note that the error bound in Theorem \ref{theorem:error_estimate_heat_equation} is not tight and the convergence rate can be improved by a more careful analysis. For example, when the piecewise linear basis functions are employed, the expected convergence rate of the finite element approximation error is O$(N^{-2})$, whereas the current analytical techniques only yield O$(N^{-1})$. However, the main purpose of Theorem \ref{theorem:error_estimate_heat_equation} is to show that the error $\Vert u_{N, M}(\bm{x}, t;\theta_P^*) - u(\bm{x}, t)\Vert_{L^2(0,T;L^2(\Omega))}$ can be bounded by a sum of three terms corresponding to the finite element approximation error, neural network approximation error and statistical error, and the convergence rates of these three errors with respect to $N$, $M$ and $P$ are all polynomial. If we suppose the optimization error is negligible, the error can be made arbitrarily small by selecting sufficiently large $N$, $M$ and $P$.

\subsection{Deep adaptive sampling method}
From the error analysis, the approximation error is determined by the hypothesis space $\mathcal{U}_N$, while the statistical error depends on the training samples. Here, random samples are used for the integration in the physical space. To reduce the statistical errors, we consider the adaptive sampling method \cite{tang2023pinns,gao2023failure}. For simplicity, we only consider the interior residual $\mathcal{L}_r(\theta)$, and the time interval is $[0,1]$. As suggested in \cite{tang2023pinns}, we relax the objective function $\mathcal{L}_r(\theta)$ as:
\begin{equation}\label{eqn:Lr_p}
	\widetilde{\mathcal{L}}_{r}(\theta) = \sum_{i=0}^{N} \lambda_i\widetilde{\mathcal{L}}_{r,i}(\theta) = \sum_{i=0}^{N} \lambda_i\int_{\Omega} r_i^2(\bm{x};\theta)p_i(\bm{x})d\bm{x} \approx \frac{1}{N_{r}}\sum_{i=0}^{N} \sum_{j=1}^{N_r} \lambda_ir_i^2(\bm{x}^{(i)}_j;\theta),
\end{equation}
where $\lambda_i>0$, $\sum_{i=0}^N\lambda_i=1$, the set $\{\bm{x}^{(i)}_j\}_{j=1}^{N_r}$ is generated with respect to the probability density function $p_i(\bm{x})>0$ instead of a uniform distribution. We associate $\widetilde{\mathcal{L}}_{r,i}(\theta)$ with a weight $\lambda_i$. The minimizer of $\widetilde{\mathcal{L}}_r(\theta)$ is also the solution to the problem if the exact solution is included in the hypothesis space. To reduce the error induced by the Monte Carlo approximation, we may adjust $p_i(\bm{x})$ to make %
the residuals $r_i^2(\bm{x};\theta)$ nearly uniform.
To do this, we refine the training set gradually by adding new samples according to the distribution induced by $r_i^2(\bm{x};\theta^{(k)}$), where $k$ indicates the adaptivity iteration and $\theta^{(k)}$ is the optimal model parameter given by the previous training set. Once the training set is updated, the model will be retrained, based on which the training set will be refined again. In a nutshell, the model and the training set are updated alternately.

The main problem is that we need to obtain new samples from $N+1$ distributions induced by $r_i^2(\bm{x};\theta^{(k)})$.
To handle this issue, we adopt a time-dependent density estimation strategy. Specifically, we augment the spatial variable in different terms ($\widetilde{\mathcal{L}}_{r,i}$) with an extra dimension $s$,  and consider a set $\{s_i\}_{i=0}^N$ of grid points on the time interval.
For the Galerkin approach, $s_i$ can be interpreted as pre-defined nodes of finite element mesh; for the collocation approach, $s_i$ can be viewed as a reordering of pre-defined Gaussian nodes $s_{m,k}$ (see Section \ref{finite_element_collocation_method}). We define
a weighted empirical measure
\begin{equation*}
	\delta_\lambda(A)=\sum_{i=0}^N\lambda_i\delta_{s_{i}}(A)
\end{equation*}
for any $A\subset[0,1]$ with $\delta_{s_i}$ being the Dirac measure
and let $r(\bm{x},s;\theta)$ be an interpolation function satisfying
\begin{equation}
	r(\bm{x},s;\theta) = r_i(\bm{x};\theta) \quad \mathrm{if} \; s = s_i, \quad i=0,\cdots,N.
\end{equation}
Let $p_{\bm{X},S}(\bm{x},s)=p_{\bm{X}|S}(\bm{x}|s)p_S(s)$ be a joint PDF. Choosing $p_S(s)ds=\delta_{\lambda}(ds)$.
We have
\begin{equation}\label{eqn:Lr_pj}
	\int\int _{\Omega} r^2(\bm{x},s;\theta)p_{\bm{X},S}(\bm{x},s)\mathrm{d}\bm{x} \mathrm{d}s=\sum_{i=0}^N\int_{\Omega}r_i^2(\bm{x};\theta)p_{\bm{X}|S}(\bm{x}|s_i)\lambda_i\mathrm{d}\bm{x},
\end{equation}
which is consistent with equation \eqref{eqn:Lr_p} if $p_{\bm{X}|S}(\bm{x}|s_i)=p_i(\bm{x})$. Using $p_{\bm{X},S}(\bm{x},s)$, the objective functional is discretized as
\begin{equation}
	\widetilde{\mathcal{L}}_r(\theta)\approx\frac{1}{N_r}\sum_{i=1}^{N_r}r^2(\bm{x}^{(i)},s^{(i)};\theta),
	\label{eqn:adaptive_loss_function}
\end{equation}
where $\{(\bm{x}^{(i)},s^{(i)})\}_{i=1}^{N_r}$ are sampled from $p_{\bm{X},S}(\bm{x},s)$. We will use a density model with the form $p_{\bm{X}|S}(\bm{x}|s)p_S(s)$ to approximate the distribution induced by $r^2(\bm{x},s;\theta)$. New samples from the trained density model will be added to the training set for refinement.

\subsubsection{Model $p_S(s)$}
Without loss of generality, we assume that $s\in[0,1]$. We aim to find a invertible transformation $z=f(s)$ such that
\begin{equation}
	\quad p_S(s) = p_Z(f(s))|\det \nabla_s f|, \quad Z\sim \mathcal{U}[0, 1],
\end{equation}
where $\mathcal{U}$ denotes the uniform distribution. We use the bounded polynomial spline layer $f_{\mathrm{poly}}$ \cite{muller2019neural} to parameterize $f$. Specifically, let $0=l_0<l_1<\cdots<l_{m-1}<l_m=1$ be a given partition of the unit interval and $\{k_j\}_{j=0}^m$ be the corresponding weights satisfying $\sum_j k_j=1$. A piecewise linear polynomial can be defined as follows:
\begin{equation}
	p(s) =\frac{k_{j+1}-k_j}{l_{j+1}-l_j}(s-l_j) + k_j,\quad \forall s\in[l_{j},l_{j+1}].
	\label{non_pdf_eq}
\end{equation}
Then the corresponding cumulative probability function $f_{\mathrm{poly}}$ admits the following formulation:
\begin{equation}
	f_{\mathrm{poly}}(s) = \frac{k_{j+1}-k_j}{2(l_{j+1}-l_j)}(s-l_j)^2+k_j(s-l_j)+\sum_{i=0}^{j-1}\frac{k_{i}+k_{i+1}}{2}(l_{i+1}-l_i),\quad \forall s\in[l_{j},l_{j+1}].
	\label{non_cdf_eq}
\end{equation}
To satisfy $\int_0^1p(s){\rm d}s=1$, we can model $\{k_j\}_{j=0}^{m}$ as
\begin{equation}
	k_j = \frac{\exp(\tilde{k}_j)}{C}, \ \forall j = 0, \ldots, m,
\end{equation}
where $\theta_{f,1} = \{\tilde{k}_j\}_{j=0}^m$ are trainable parameters and $C$ is a normalization constant:
\begin{equation}
	C = \sum_{i=0}^{m-1}\frac{(\exp(\tilde{k}_i)+\exp(\tilde{k}_{i+1}))(l_{i+1}-l_i)}{2}.
\end{equation}

Notice that the polynomial spline layer $f_{\mathrm{poly}}(\cdot;\theta_{f,1})$ (\ref{non_pdf_eq})-(\ref{non_cdf_eq}) yields explicit monotonous expressions, and its inverse can be readily computed. Then an explicit PDF model $p_{\mathrm{poly}}(s;\theta_{f,1})$ can be obtained by letting $f=f_{\mathrm{poly}}$, i.e.,
\begin{equation}
	p_{\mathrm{poly}}(s;\theta_{f,1}) = p_Z(f_{\mathrm{poly}}(s))\left|\det \nabla _{s}f_{\mathrm{poly}}\right|.
\end{equation}
\subsubsection{Model $p_{X|S}(\textbf{x}|s)$}
For $\bm{x}\in\mathbb{R}^d$, we seek a invertible transformation $\bm{z}=f(\bm{x}, s)\in\mathbb{R}^d$ for any given $s$ such that
\begin{equation}
	p_{\bm{X}|S}(\bm{x}|s) = p_{\bm{Z}|S}(\bm{z}|s)\left| \frac{\partial f(\bm{x},s)}{\partial \bm{x}}\right|, \quad \bm{Z}|S\sim \mathcal{U}[-1, 1]^{d},\quad  \forall s.
	\label{conditonal_change_of_variable}
\end{equation}
Here we employ conditional bounded KR-net $f_{\mathrm{B-KRnet}}(\cdot, s)$ \cite{zeng2023bounded} to parameterize $f(\cdot, s)$. The basic idea of conditional bounded KRnet is to define the structure of $f(\bm{x},s)$ in terms of the Knothe-Rosenblatt rearrangement. The transformation $f(\cdot, s)$ inspired by the Knothe-Rosenblatt (K-R) rearrangement \cite{carlier2010knothe} has a low-triangular structure
\begin{equation}
	\bm{z} = f(\bm{x},s) = \left[\begin{array}{l}
			f_1(x_1,s)     \\
			f_2(x_1,x_2,s) \\
			\vdots         \\
			f_d(x_1,\cdots,x_d,s)
		\end{array}\right].
\end{equation}
The sub-transformations $f_1,\cdots,f_d$ consist of polynomial spline layers and coupling layers \cite{dinh2016density}. More details can be found in \cite{zeng2023bounded,feng2022solving}. Let $f_{\mathrm{B-KRnet}}(\cdot,s;\theta_{f,2})$ indicate the conditional invertible transport map induced by  bounded KR-net, where $\theta_{f,2}$ includes the model parameters. Then an explicit PDF model $p_{\mathrm{B-KRnet}}(\bm{x},s;\theta_{f,2})$ can be obtained by letting $f=f_{\mathrm{B-KRnet}}$ in equation \eqref{conditonal_change_of_variable}
\begin{equation}
	p_{\mathrm{B-KRnet}}(\bm{x}|s;\theta_{f,2}) = p_{\bm{Z}}(f_{\mathrm{B-KRnet}}(\bm{x},s)) \left|\det \nabla_{\bm{x}} f_{\mathrm{B-KRnet}}\right|.
\end{equation}

\subsubsection{Adaptive sampling approach}
Now we model a continuous joint density distribution $p_{\theta_f}(\bm{x},t)$
\begin{equation}
	p_{\theta_f}(\bm{x},t) = p_{\mathrm{poly}}(t;\theta_{f,1}) p_{\mathrm{B-KRnet}}(\bm{x}|t;\theta_{f,2}),
\end{equation}
where $\theta_f = \{\theta_{f,1}, \theta_{f,2}\}$. To seek the "optimal" parameter $\theta_f$, we can minimize the following objective
\begin{equation}
	\begin{aligned}
		D_{\mathrm{KL}}(\hat{r}_{\theta}(\bm{x},t) || p_{\theta_f}(\bm{x},t))
		 & = D_{\mathrm{KL}}\big(\hat{r}_{\theta}(\bm{x},t)|| p_{\mathrm{poly}}
		(t;\theta_{f,1})p_{\mathrm{B-KRnet}}(\bm{x}|t;\theta_{f,2})\big)                                              \\
		 & =\iint \hat{r}_{\theta}(\bm{x},t) \log \left(\hat{r}_{\theta}(\bm{x},t)\right) \mathrm{d}\bm{x}\mathrm{d}t
		- \iint \hat{r}_{\theta}(\bm{x},t) \log \left( p_{\mathrm{poly}}(t;\theta_{f,1})p_{\mathrm{B-KRnet}}(\bm{x}|t;\theta_{f,2})\right) \mathrm{d}\bm{x}\mathrm{d}t,
	\end{aligned}
	\label{density_kl}
\end{equation}
where $D_{\mathrm{KL}}$ indicates the Kullback-Leibler (KL) divergence and $\hat{r}_{\theta}(\bm{x},t)\propto r^2(\bm{x},t;\theta)$ is the induced measure by continuous residual squared $r^2(\bm{x},t;\theta)$
\begin{equation*}
	r(\bm{x},t;\theta) = \partial_t u_N(\bm{x},t;\theta) + \mathcal{N}[u_N](\bm{x},t;\theta) - f(\bm{x},t).
\end{equation*}

The first term on the right-hand side in \eqref{density_kl} corresponds to the negative differential entropy of
$\hat{r}_{\theta}(\bm{x},t)$, which does not affect the optimization with respect to
$\theta_f$. So minimizing the KL divergence is equivalent to minimizing the cross
entropy between $\hat{r}_{\theta}(\bm{x},t)$ and $p_{\theta_f}(\bm{x},t)$
\begin{equation}
	\begin{aligned}
		H(\hat {r}_{\theta}(\bm{x},t), p_{\theta_f}(\bm{x},t)) &
		= - \iint \hat{r}_{\theta}(\bm{x},t) \log \left( p_{\mathrm{poly}}(t;\theta_{f,1})
		p_{\mathrm{B-KRnet}}(\bm{x}|t;\theta_{f,2})\right) \mathrm{d}\bm{x}\mathrm{d}t,
	\end{aligned}
\end{equation}
Since the samples from $\hat{r}_{\theta}(\bm{x},t)$ are not available, we approximate the cross entropy using the importance sampling technique:
\begin{equation}
	H(\hat{r}_{\theta}(\bm{x},t),p_{\theta_f}(\bm{x},t))\approx -\frac{1}{N_r}\sum_{i=1}^{N_r}
	\frac{\hat{r}_{\theta}(\bm{x}_i, t_i)}{p_{\mathrm{poly}}(t_i;\hat{\theta}_{f,1})
		p_{\mathrm{B-KRnet}}(\bm{x}_i|t_i;\hat{\theta}_{f,2})}
	\left(\log p_{\mathrm{poly}}(t_i;\theta_{f,1}) + \log p_{\mathrm{B-KRnet}}
	(\bm{x}_i|t_i;{\theta}_{f,2})\right),
	\label{cross_entropy}
\end{equation}
where
\begin{equation}
	t_i\sim p_{\mathrm{poly}}(\cdot;\hat{\theta}_{f,1}),\quad  \bm{x}_i\sim p_{\mathrm{B-KRnet}}(\cdot|t_i;\hat{\theta}_{f,2}).
	\label{sample_from_joint_distribution}
\end{equation}
The choice of $\hat{\theta}_{f}=\{\hat{\theta}_{f,1},\hat{\theta}_{f,2}\}$ will be specified as follows.

Once obtaining the well-trained parameters $\theta_f^* = \{\theta_{f,1}^*, \theta_{f,2}^*\}$, one can refine initial training set
via adding new samples and update the weights in \eqref{eqn:Lr_p}. Specifically, note that the
residuals $r_i(\bm{x};\theta)$ in \eqref{def_r_i} are only needed at time $s_i$,
we need to construct the corresponding discrete distributions.
For system \eqref{galerkin_system_neural_linear} in Galerkin projection, we use the following discrete distribution
\begin{equation}
	p_{\mathrm{dis}}(s) =\left\{ \begin{array}{rl}
		\frac{\int_{\Omega}r_i^2(\bm{x};\theta)\mathrm{d}\bm{x}}{\sum_{i=0}^{N}\int_{\Omega}r_i^2(\bm{x};\theta)\mathrm{d}\bm{x}}, & s=s_i,0\leq i\leq N \\
		0,                                                                                                                         & \mathrm{otherwise}.
	\end{array}\right.
	\label{galerkin_discrete_density}
\end{equation}

For system \eqref{collocation_system_neural_linear} in collocation projection, we simply use the following discrete distribution
\begin{equation}
	p_{\mathrm{dis}}(s;\theta_{f,1}^*) = \left\{
	\begin{array}{rl}
		\int_{0}^{(s_0+s_1)/2}p_{\mathrm{poly}}(s;\theta_{f,1}^*)\mathrm{d}s,                   & s=s_0,              \\
		\int_{(s_{i-1}+s_i)/2}^{(s_i+s_{i+1})/2}p_{\mathrm{poly}}(s;\theta_{f,1}^*)\mathrm{d}s, & s = s_i,0<i<N       \\
		\int_{(s_{N-1}+s_N)/2}^1p_{\mathrm{poly}}(s;\theta_{f,1}^*)\mathrm{d}s,                 & s=s_N,              \\
		0,                                                                                      & \mathrm{otherwise}.
	\end{array}
	\right.
\end{equation}Then the weights in \eqref{eqn:Lr_p} can be determined via the discrete distribution $p_{\mathrm{dis}}(\cdot;\theta_{f,1}^*)$.
For collocation projection, it is straightforward to refine the
training set. We first generate $N_{\mathrm{new}}$ samples $\{t_j\}_{j=1}^{N_{\mathrm{new}}}$ via the
well-trained model
$p_{\mathrm{dis}}(\cdot;\theta^*_{f,1})$, then we generate the corresponding spatial samples
$\bm{x}_j\sim p_{\mathrm{B-KRnet}} (\cdot|t_j;\theta^*_{f,2})$ for each $t_j$.
After that we reorder newly generated data into the original spatial space
\begin{equation}
	\bm{x}_j \in S^{(i)}_{r,\mathrm{new}}, \quad \mathrm{if}\; t_j = s_i, \forall i=0,\cdots,N,\; j=1,\cdots,N_{\mathrm{new}}.
\end{equation}
For Galerkin projection, recall that we need to use the same training set for the basis function
with same support (see \eqref{galerkin_continuous_residual}). Hence for each $s_i$ in \eqref{galerkin_discrete_density}, one
firstly identies its corresponding basis function and support, here it is denoted as $\tilde{I}_i$. Then we generate
spatial samples $\bm{x}_{i,k}\sim p_{\mathrm{B-KRnet}}(\cdot|t_{i,k};\theta_{f,2}^*)$, where $t_{i,k}$ are Gaussian quadrature points in the interval
$\tilde{I}_i$. The total number of generated spatial points for each $s_i$ should be proportional to the discrete weights $p_{\mathrm{dis}}(s_i)$. Again, we recoder newly
generated data into the original spatial space
\begin{equation}
	\bm{x}_{i,k} \in S^{(i)}_{r,\mathrm{new}}, \quad \forall i=0,\cdots,N.
\end{equation}

We are now ready to present our algorithms.
Let $\{S_{r,0}^{(i)}\}_{i=0}^{N}, S_{ic}$ and $S_{bc}$ be three sets of collocation points that are
uniformly sampled from $\Omega^{N+1}$, $\Omega\times \{0\}$ and $\partial \Omega\times [0,T]$, respectively.
Using $\{S_{r,0}^{(i)}\}_{i=0}^{N}, S_{ic}$ and $S_{bc}$, we minimize the empirical
loss \eqref{eqn:empricial_loss} to obtain $u(\bm{x},t;\theta^{*, (1)})$. With
$u(\bm{x},t;\theta^{*,(1)})$, we minimize the cross entropy \eqref{cross_entropy}
to get $p_{\theta^{*,(1)}_f}(\bm{x},t)$, where we use uniform samples for importance
sampling. To refine the training set, a new set
$S_{r,\mathrm{new}}=\{S^{(i)}_{r,\mathrm{new}}\}_{i=0}^{N}$ is generated according to the model $p_{\theta_f^{*,(1)}}(\bm{x},t)$, and new training set can be updated as
$S_{r,1} = S_{r,0}\cup S_{r,\mathrm{new}}$. Then we continue to update the approximate solution $u(\bm{x},t;\theta^{*,(1)})$ using $S_{r,1}$ as the training set, which yields a refined model $u(\bm{x},t;\theta^{*,(2)})$. Starting from $k=2$, we seek $p_{\theta^{*,k}}(\bm{x},t)$ using the previous approach. We repeat the procedure to obtain an adaptive algorithm (see Algorithm \ref{alg:2}).
\begin{algorithm}
	\caption{Deep adaptive sampling method for hybrid FEM-PINNs method}
	\label{alg:2}
	\begin{algorithmic}
		\State \textbf{Input:} number of basis functions $N+1$, number of adaptive iteration $N_{\mathrm{adaptive}}$, number of adaptive training epochs $N_a$, number of newly added points $N_{\mathrm{new}}$, initial weights $\{\lambda_i^{(0)}\}_{i=0}^N$, training data $S_{r}^{(0)}$, $S_{ic}$ and $S_{bc}$
		\State Solve time-dependent PDE via Algorithm \ref{alg:1} to obtain $u(\bm{x}, t;\theta^{*,(0)})$.
		\For{$k=0,\cdots,N_{\mathrm{adaptive}}-1$}
		\State // Train probability density model.
		\For {$j=1,\cdots,N_a$}
		\State Sample $m$ samples from $S_{r}^{(k)}$.
		\State Update $p(\bm{x}, t;\theta_{f}^{*, (k+1)})$ by descending the stochastic gradient of \eqref{cross_entropy}
		\EndFor
		\State // Update weights and refine training set.
		\State Update weights $\{\lambda_i^{(k+1)}\}_{i=0}^N$ and generate new training samples $S_{r}^{(\mathrm{new})}$ through $p(\bm{x}, t;\theta_{f}^{*, (k+1)})$.
		\State Update training set $S_{r}^{(k+1)} = S_{r}^{(k)}\cup S_{r}^{(\mathrm{new})}$.
		\State Solve time-dependent PDE via Algorithm \ref{alg:1} to obtain $u(\bm{x}, t;\theta^{*,(k+1)})$.
		\EndFor
		\State \textbf{Output:} predicted solution $u(\bm{x},t;\theta^*)$.
	\end{algorithmic}
\end{algorithm}

\section{Numerical experiments}\label{section:numerical_experiments}
In this section, we conduct some numerical experiments to demonstrate the effectiveness of the proposed method, including one convection equation, one Allen-Cahn equation, one two-dimensional and complex geometry test problem, one two-dimensional and low regularity test problem, and two high-dimensional linear or nonlinear problems. The low-dimensional benchmarks are included mainly to verify the accuracy of the framework and to facilitate comparison with standard examples widely used in the machine-learning literature, even though classical numerical solvers are typically sufficient for such problems. The high-dimensional problems, on the other hand, are intended to illustrate the applicability of the method in more challenging scenarios where conventional approaches face significant computational difficulties. Throughout all benchmarks, we will employ the fully-connected neural network equipped with hyperbolic tangent activation functions (Tanh) and initialized using the Glorot normal scheme \cite{glorot2010understanding}. All neural networks are trained via stochastic gradient descent using the Adam optimizer with default settings \cite{kingma2014adam} and an exponential learning rate decay with a decay-rate of 0.9 every 1,000 training iterations. All experiments are implemented by JAX \cite{bradbury2018jax}.

In order to test the validity of the method, we use the following relative $L_2$ error:
\begin{equation}
	err_{L_2} = \frac{\sqrt{\sum_{i=1}^{\mathrm{num}}|u(\bm{x}_i,t_i;\theta) - u(\bm{x}_i,t_i)|^2}}{\sqrt{\sum_{i=1}^{\mathrm{num}}|u(\bm{x}_i,t_i)|^2}},
\end{equation}
where $\mathrm{num}$ represents the total number of test points chosen randomly in the domain, and $u(\bm{x}_i,t_i;\theta)$ and $u(\bm{x}_i, t_i)$ represent the predicted and the exact solution values, respectively. Five independent numerical experiments are conducted with varying initializations, and the mean and standard deviation of the relative $L_2$ error are reported.
\subsection{Convection equation}
\begin{figure}[H]
	\centering
	\includegraphics[width=0.3\linewidth]{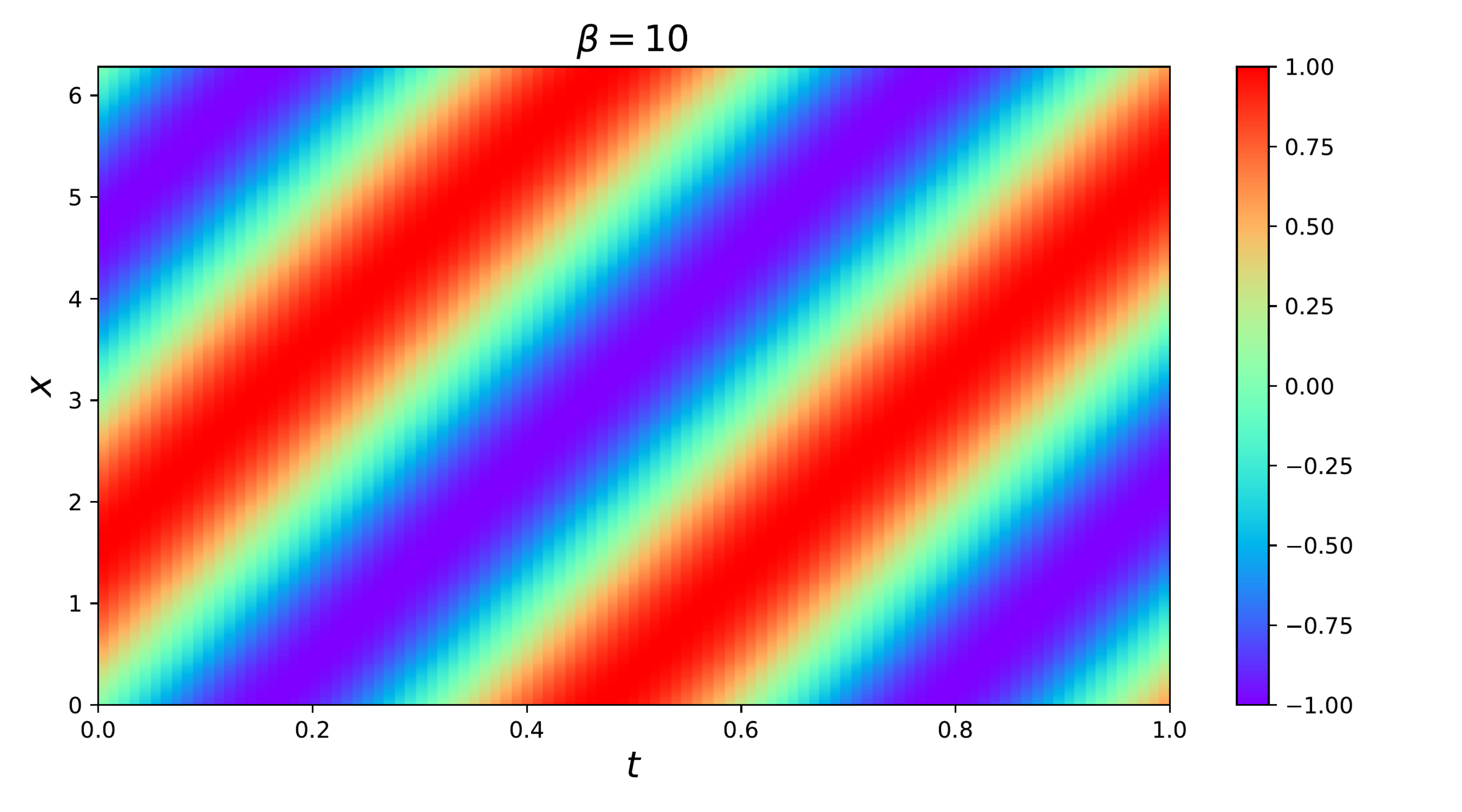}
	\includegraphics[width=0.3\linewidth]{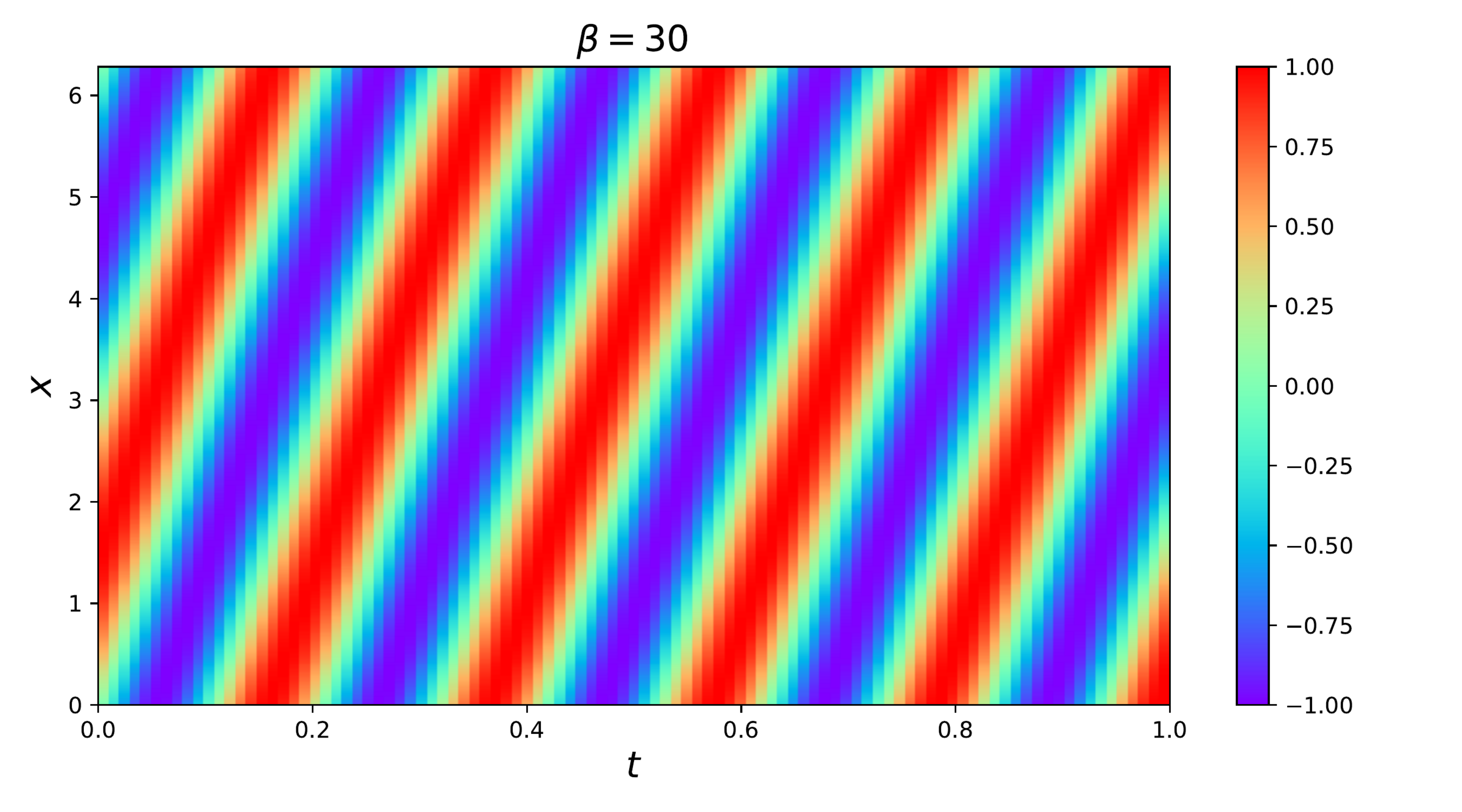}
	\includegraphics[width=0.3\linewidth]{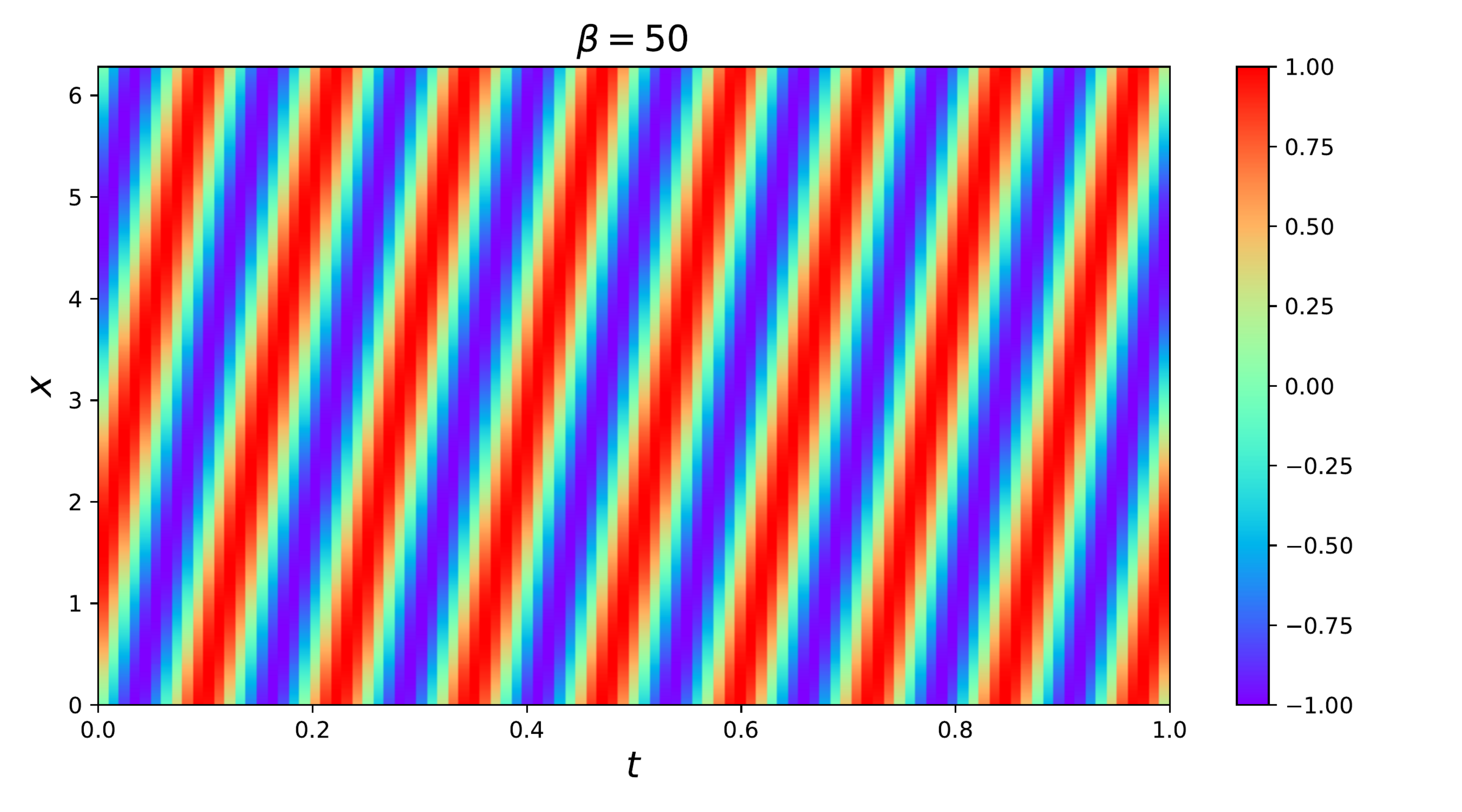}
	\caption{Convection equation. Exact solutions. From left to right: $\beta=10,30$ and $50$.}
	\label{fig:exactbeta}
\end{figure}

We start with a one-dimensional linear convection equation of the form
\begin{equation}
	\frac{\partial u}{\partial t} + \beta \frac{\partial u}{\partial x} = 0, \quad (x,t)\in [0, 2\pi] \times [0,1],
\end{equation}
subject to periodic boundary conditions and an initial condition $u(0,x)=\sin(x)$.
The precise solutions for varying values of $\beta$ are depicted in Figure \ref{fig:exactbeta}. One can observe that as $\beta$ increases, the solution exhibits increasingly pronounced temporal variations.

\begin{figure}[H]
	\centering
	\includegraphics[width=0.45\linewidth]{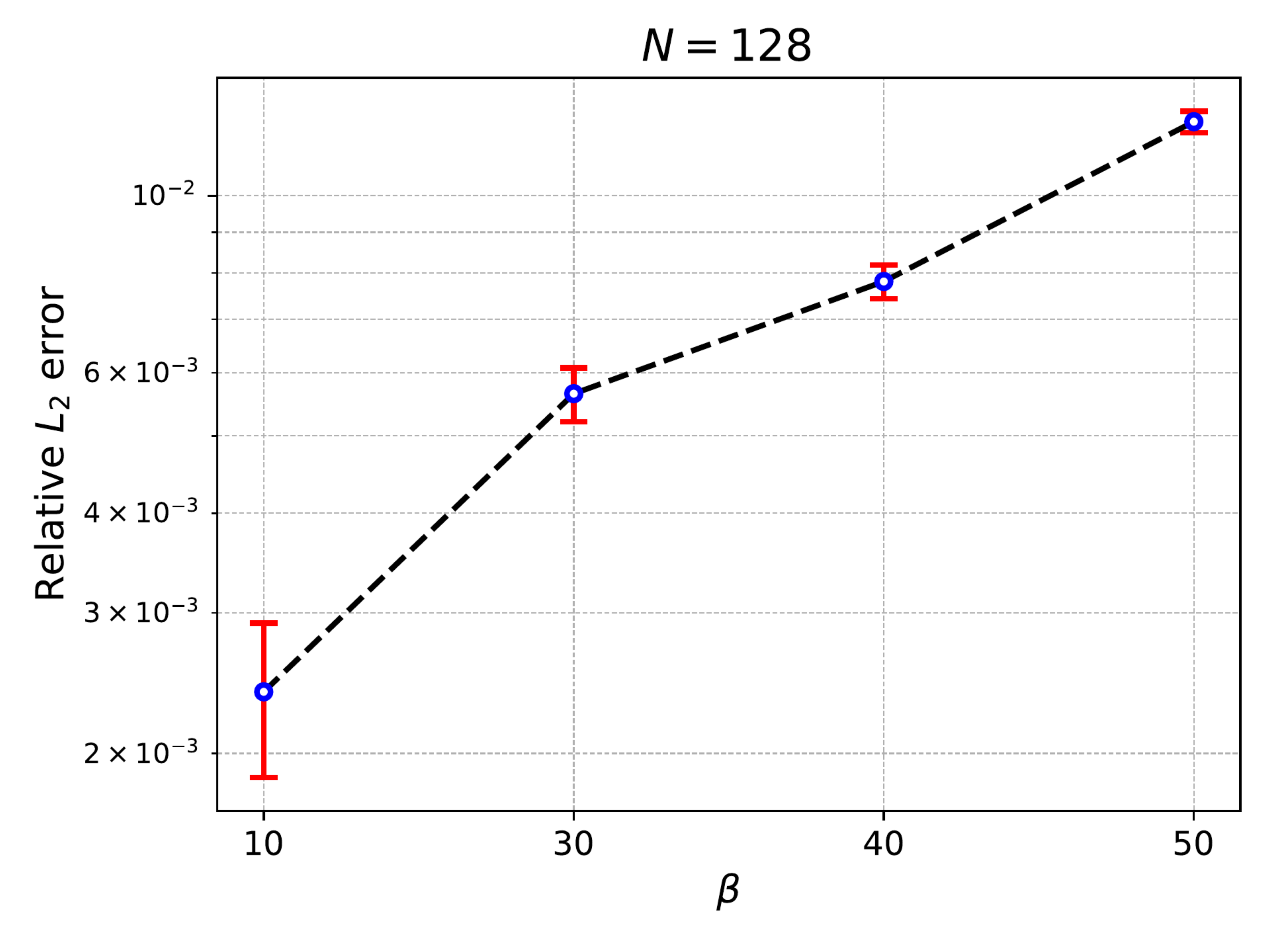}
	\includegraphics[width=0.45\linewidth]{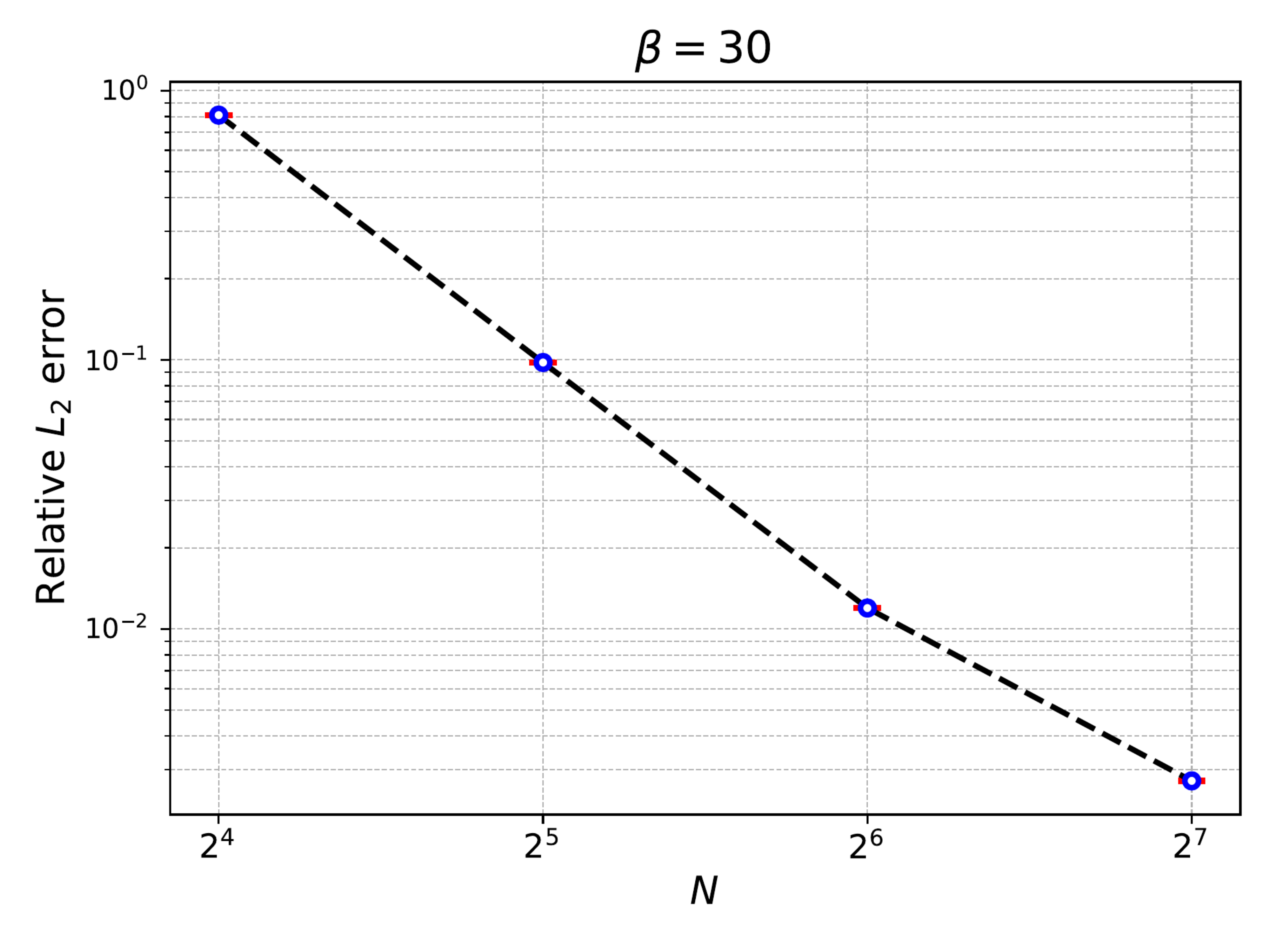}
	\caption{Convection equation. Numerical results with linear finite element discretization.}
	\label{fig:convection}
\end{figure}

We initially evaluate the performance with  linear finite element basis functions in Galerkin projection. The latent coefficients, denoted as $\{\omega_i(\bm{x};\theta)\}$, are represented using a fully-connected neural
network with tanh activation function, 4 hidden layers and 128 neurons per hidden layer. To simplify the training
object, we strictly impose the periodic boundary conditions by embedding the input coordinates into Fourier expansion
(see \ref{app_periodic_boundary_condition}). Notice that the above equation is linear, we use the linear
form of loss function \eqref{NN_linear_finite_element_martix_form} and set $N_{\mathrm{segment}}=1$. We create a uniform
mesh of size 400 in the spatial computational domain $[0,2\pi]$, yielding 400 initial points and 400 collocation points for enforcing the PDE residual. We proceed by training the
resulting model via full-batch gradient descent using the Adam optimizer for 40,000 iterations. As shown in Figure \ref{fig:convection}, for fixed $N$,
the relative $L_2$ error gradually increases as $\beta$ increases; for fixed $\beta=30$, the relative $L_2$ error exhibits linear convergence concerning the number of basis functions ($N$). Particularly, when $\beta=30$, the proposed method attains a remarkable relative $L_2$ error of $2.85e-3$.

\begin{figure}[H]
	\centering
	\includegraphics[width=0.45\linewidth]{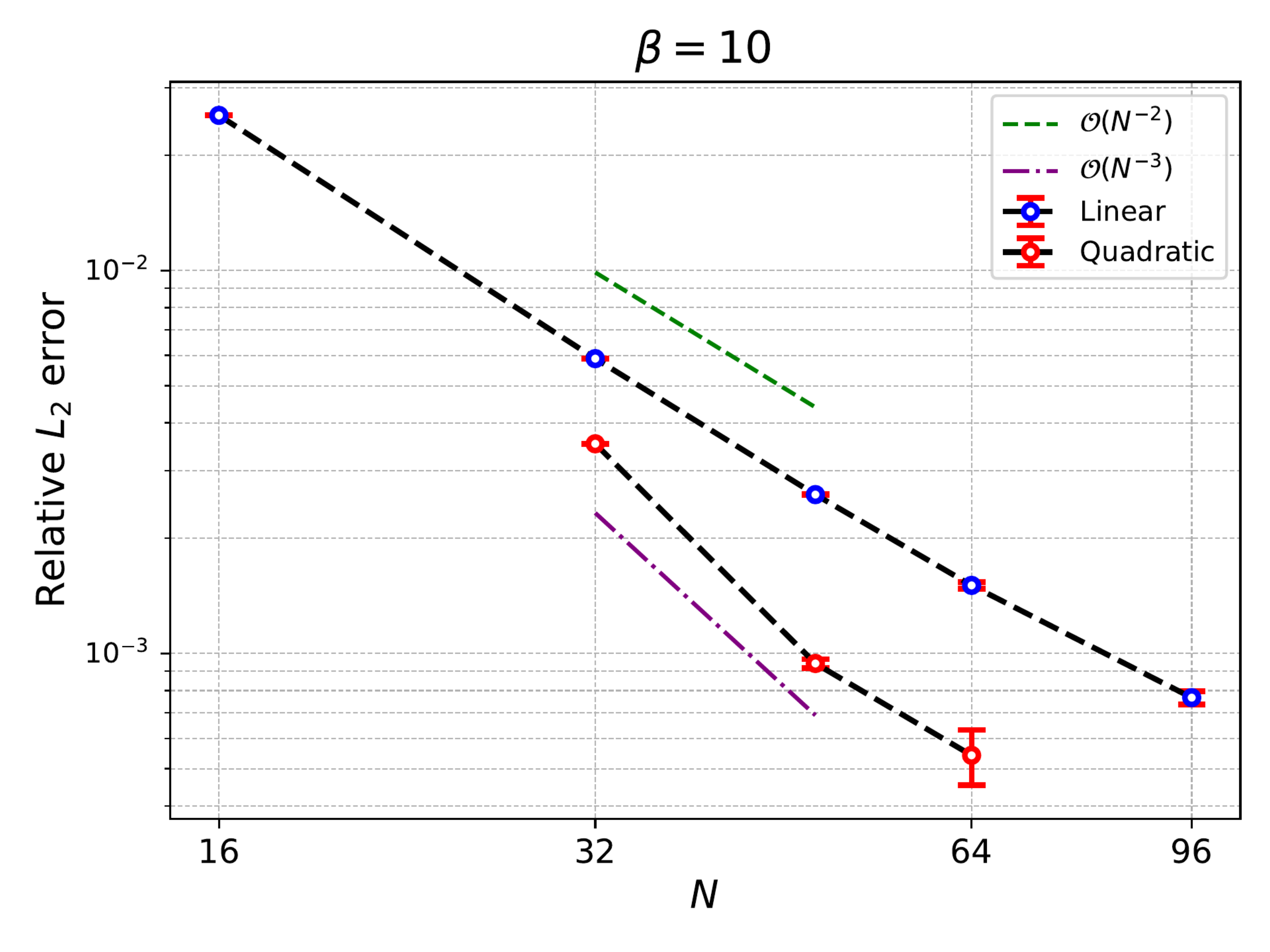}
	\includegraphics[width=0.45\linewidth]{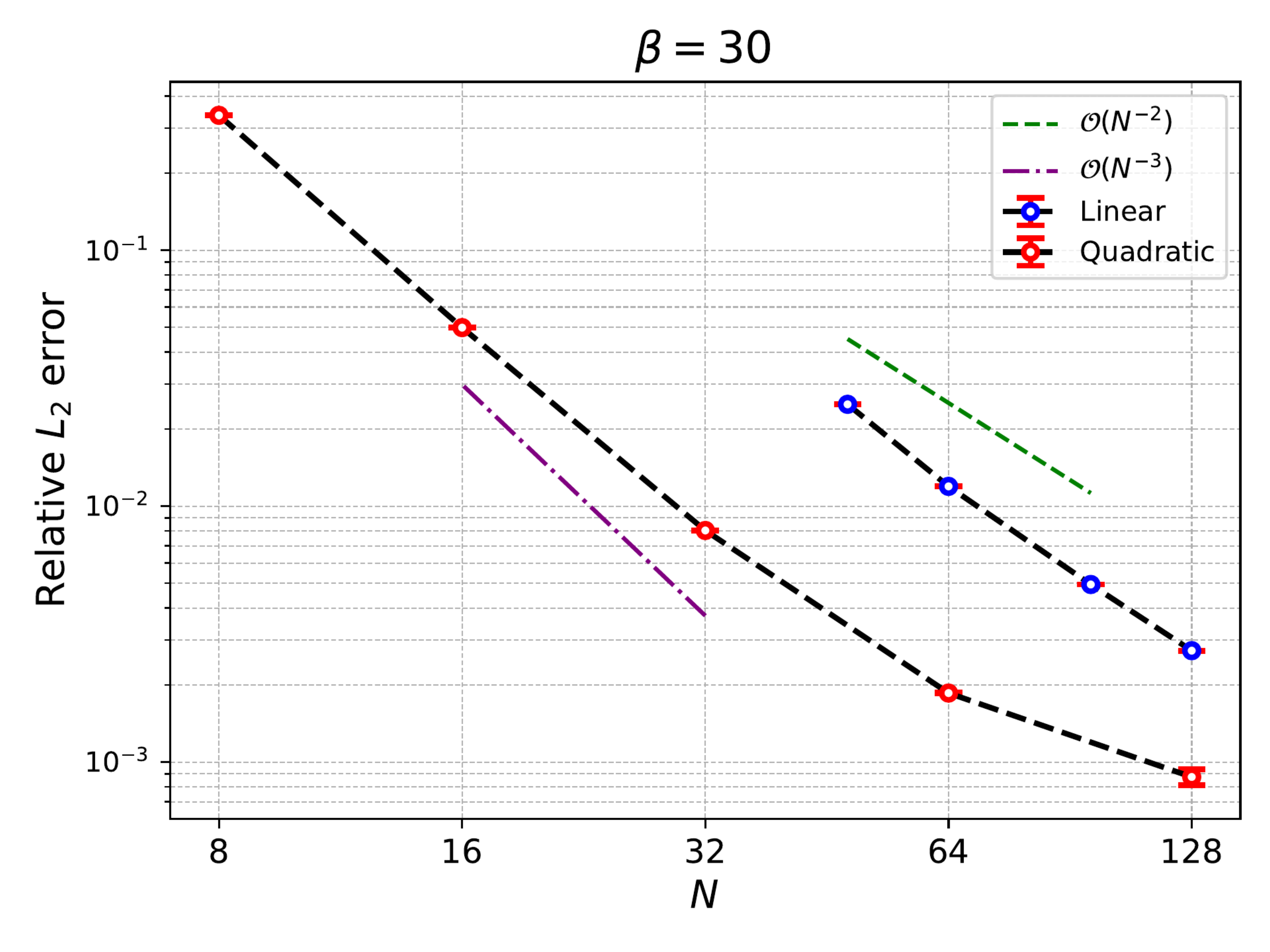}
	\caption{Convection equation. Convergence rate for linear finite element basis functions versus quadratic finite element basis functions.}
	\label{fig:convectionlinearquadratic}
\end{figure}

Furthermore, we investigate the impact of linear and quadratic finite element basis functions on the performance of the proposed model. Specifically, we consider two cases: $\beta$ = 10 with the number of mesh elements $N$ ranging from 16 to 96, and $\beta$ = 30 with $N$ ranging from 8 to 128. We then train the proposed model under same hyperparameter configurations. Figure \ref{fig:convectionlinearquadratic} presents a summary of the relative $L_2$ errors observed in the trained models. It is not surprising that the error diminishes when we replace linear basis functions with quadratic ones,  in accordance with the classical results of finite element theory. More specifically, when $N$ is relatively small, the relative $L_2$ errors exceed the expected theoretical convergence rates ($\mathcal{O}(h^2)$ for linear basis functions and $\mathcal{O}(h^3)$ for quadratic basis functions). Additionally, when $N$ is relatively large, the corresponding error decreases at a slower rate than the theoretical convergence rate. This suggests that for larger $N$, the spatial error begins to dominate due to the increased difficulty of optimization, while the temporal error becomes less significant.

To investigate the effects of increasing the number of basis functions, we compare the relative $L_2$ errors given by various settings
when $\beta=50$. It is seen in Table \ref{convection_relative_error_galerkin_collocations} that relative $L_2$ error decreases as $N$ increases for all cases.

\begin{table}[H]
	\centering
	\resizebox{1.0\linewidth}{!}{
		\begin{tabular}{ccccc}
			\toprule
			$N$                              & 32                      & 64                      & 96                      & 128                     \\
			\midrule
			Galerkin (linear)                & 1.69e+00 $\pm$ 2.38e-05 & 7.73e-02 $\pm$ 1.23e-05 & 1.91e-02 $\pm$ 1.06e-05 & 8.69e-03 $\pm$ 1.18e-05 \\
			Galerkin (quadratic)             & 3.83e-02 $\pm$ 8.53e-06 & 5.58e-03 $\pm$ 1.15e-05 & 2.26e-03 $\pm$ 2.67e-05 & 1.13e-03 $\pm$ 3.74e-05 \\
			Collocation (Hermite polynomial) & 1.97e-01 $\pm$ 1.78e-06 & 1.18e-02 $\pm$ 6.76e-07 & 2.34e-03 $\pm$ 2.17e-04 & 7.50e-04 $\pm$ 2.92e-04 \\
			Collocation (spline)             & 6.44e-01 $\pm$ 2.29e-05 & 1.83e-02 $\pm$ 1.13e-06 & 3.16e-03 $\pm$ 2.31e-05 & 9.83e-04 $\pm$ 7.18e-07 \\
			\bottomrule
		\end{tabular}}
	\caption{Convection equation ($\beta=50$). Relative $L_2$ error (mean $\pm$ std) at different $N$ for different methods.}
	\label{convection_relative_error_galerkin_collocations}
\end{table}
In Table \ref{convection_relative_error_DABG} we compare the performance of the proposed method with the DABG \cite{gu2022deep}
at different $\beta$. One can observe that for these two methods, when $\beta$ is relatively small, as $N$ increases, the error first decreases and then increases; when $\beta$ is relatively large, the error consistently decreases as $N$ increases. We conjecture that reason for this phenomenon is that when $\beta$ is relatively small (i.e., the solution is smoother), small $N$ can achieve good accuracy, while larger $N$ may lead to greater optimization challenges. Another message from Table \ref{convection_relative_error_DABG} is that the proposed approach exhibits reduced sensitivity to temporal frequency variations in the solution compared to DABG.
\begin{table}[!h]
	\centering
	\resizebox{1.0\linewidth}{!}{
		\begin{tabular}{cccccccc}
			\toprule
			\multicolumn{2}{c}{$N$}     & 8                   & 16       & 32       & 64       & 96       & 128                 \\
			\midrule
			\multirow{2}{*}{$\beta=10$} & Galerkin(quadratic) & 1.55e-02 & 3.52e-03 & 1.03e-03 & 4.67e-04 & 6.91e-04 & 3.91e-01 \\
			                            & DABG                & 1.36e-02 & 1.25e-04 & 8.73e-04 & 1.04e-05 & 3.33e+00 & 1.32e+00 \\
			\midrule
			\multirow{2}{*}{$\beta=30$} & Galerkin(quadratic) & 3.37e-01 & 4.98e-02 & 8.03e-03 & 1.90e-03 & 9.03e-04 & 5.04e-04 \\
			                            & DABG                & 1.12e+00 & 2.87e-01 & 2.93e-04 & 4.34e-04 & 5.75e-05 & 1.06e-04 \\
			\midrule
			\multirow{2}{*}{$\beta=50$} & Galerkin(quadratic) & 4.67e+00 & 2.15e-01 & 3.83e-02 & 5.17e-03 & 2.08e-03 & 1.32e-03 \\
			                            & DABG                & 1.05e+00 & 1.03e+00 & 3.30e-03 & 1.44e-03 & 6.46e-04 & 6.39e-04 \\
			\bottomrule
		\end{tabular}}
	\caption{Convection equation. Different relative $L_2$ errors between the proposed method and DABG.}
	\label{convection_relative_error_DABG}
\end{table}

We also compare the performance and computational time of the proposed method with causal PINNs. We set the same hyper parameters to compare these two methods, including the architecture of neural network, the initial learning rate, decay rate, the number of training points in space, and iterations. The $\epsilon$ in causal PINNs is set to $\{1e-2, 1e-1, 1e0, 1e1, 1e2\}$. All runtime statistics were computed on the same hardware, a Nvidia Tesla V100 w/ 32 GB memory. It is shown in Table \ref{convection_relative_error_causal}  that as $\beta$ increases, the proposed method demonstrates superior accuracy compared to causal PINNs. The proposed method (Galerkin projection) runs approximately 100 times faster than causal PINNs. Such a significant difference in efficiency is due to the distinct computational graphs adopted by these two methods. For the proposed method, because the residual originates from the linear system defined in \eqref{NN_linear_finite_element_martix_form}, we can pre-calculate the finite element sparse matrix along the time direction. We consider the computational cost of automatic differentiation in calculating the residual. For any given spatial point $\bm{x}$, causal PINNs need $n$ back-propagation computations for $n$ temporal points, while our method only requires a single back-propagation computation and one matrix-vector multiplication of the above-mentioned finite element sparse matrix and the neural network output.
\begin{table}[!h]
	\centering
	\begin{tabular}{ccccc}
		\toprule
		\multicolumn{2}{c}{$\beta$}   & 10                  & 30       & 50                  \\
		\midrule
		\multirow{2}{*}{Error}        & Galerkin(quadratic) & 2.95e-04 & 1.32e-03 & 3.57e-03 \\
		                              & Causal PINNs        & 7.58e-04 & 1.12e-02 & 1.65e+00 \\
		\midrule
		\multirow{2}{*}{Running time} & Galerkin(quadratic) & 57       & 62       & 63       \\
		                              & Causal PINNs        & 6639     & 8087     & 8057     \\
		\bottomrule
	\end{tabular}
	\caption{Convection equation. Relative $L_2$ error and running time between the proposed method and causal PINNs.}
	\label{convection_relative_error_causal}
\end{table}

\subsection{Allen-Cahn equation}
The next example aims to illustrate the effectiveness of time marching technology in our proposed method. Consider the following Allen-Cahn equation
\begin{equation}
	\begin{array}{rcl}
		u_t -c_1^2 \nabla ^2 u +f(u) & = & 0, \quad x\in[-1, 1],t\in [0,1], \\
		f(u)                         & = & c_2 (u^3 - u),                   \\
		u(x,0)                       & = & x^2 \cos(\pi x),                 \\
		u(1,t)                       & = & u(-1, t),                        \\
		u_x(1, t)                    & = & u_x(-1, t),
	\end{array}
\end{equation}
where $c_1^2=0.0001$ and $c_2=5$.
We take the number of linear mesh elements as 100 and represent the latent coefficients $\{\omega_i(\bm{x};\theta)\}$ by a fully-connected neural network with tanh activation function, 4 hidden layers and 128 neurons per hidden layer. Similarly, we strictly impose the periodic boundary conditions by embedding the input coordinates into Fourier expansion (see \ref{app_periodic_boundary_condition}). We create a uniform mesh of size 1,000 in the spatial computational domain $[-1, 1]$, yielding 1,000 initial points and 1,000 collocation points for enforcing the PDE residual. We proceed by training the resulting model via full-batch gradient descent using the Adam optimizer for
the total 100,000 iterations. The resulting $L_2$ errors for different numbers of time segments are shown in Table \ref{table_allen_cahn_linear}. We observe the relative $L_2$ error significantly decreases as the number of time segments increases. One can see that the predicted solution achieves an excellent agreement with the ground truth, yielding an approximation error of $7.67e-3$ measured in the relative $L_2$ norm.

\begin{table}[H]
	\centering
	\begin{tabular}{cccc}
		\toprule
		$N_{\mathrm{segment}}$ & 1  & 2                       & 4                       \\
		\midrule
		Error                  & -- & 1.17e-02 $\pm$ 4.36e-03 & 7.67e-03 $\pm$ 1.27e-03 \\
		\bottomrule
	\end{tabular}
	\caption{\centering Allen-Cahn equation. Relative $L_2$ error (mean $\pm$ std) for different segments. (Linear finite element basis function with $N=100$)}
	\label{table_allen_cahn_linear}
\end{table}

\begin{figure}[H]
	\centering
	\includegraphics[width=0.3\linewidth]{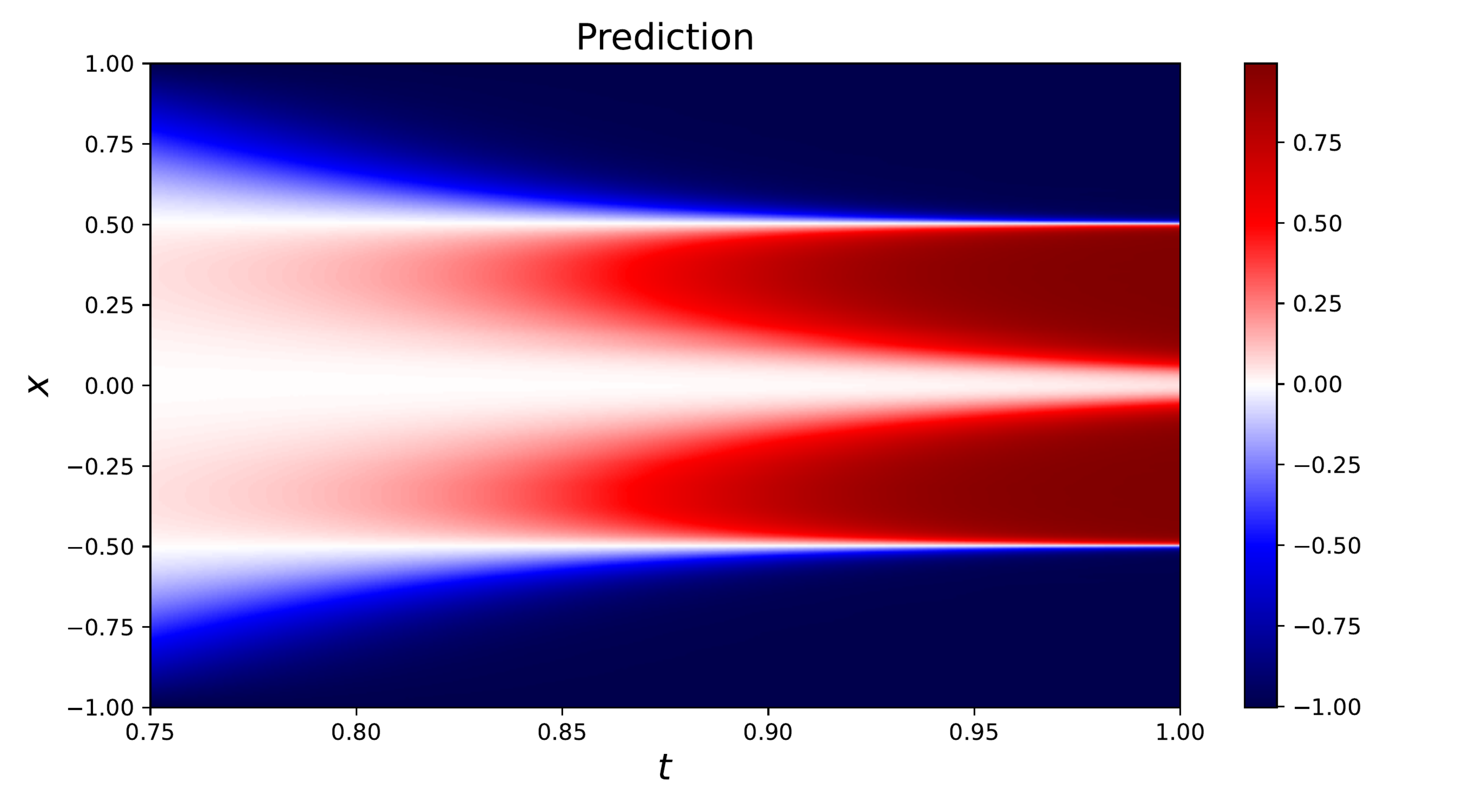}
	\includegraphics[width=0.3\linewidth]{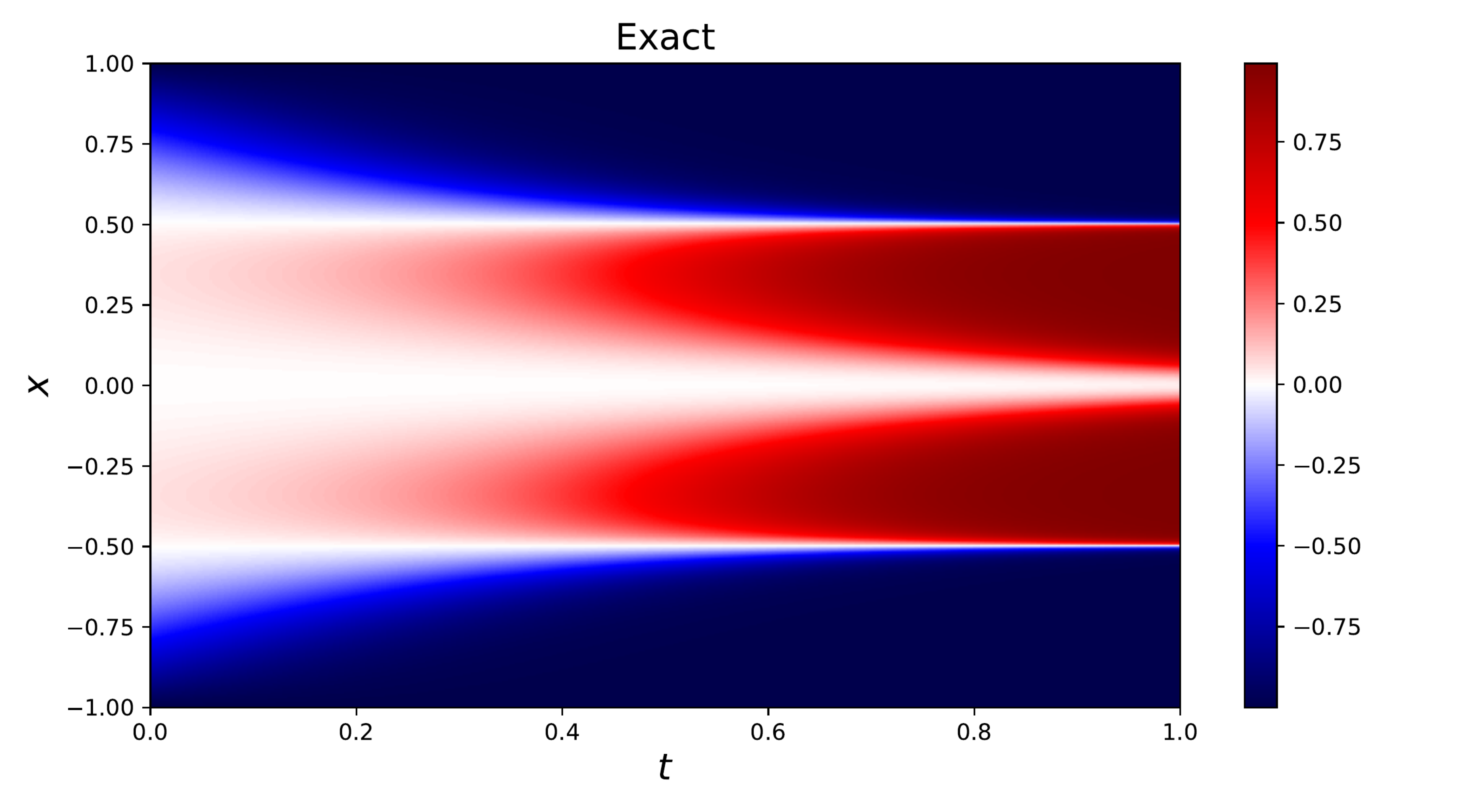}
	\includegraphics[width=0.3\linewidth]{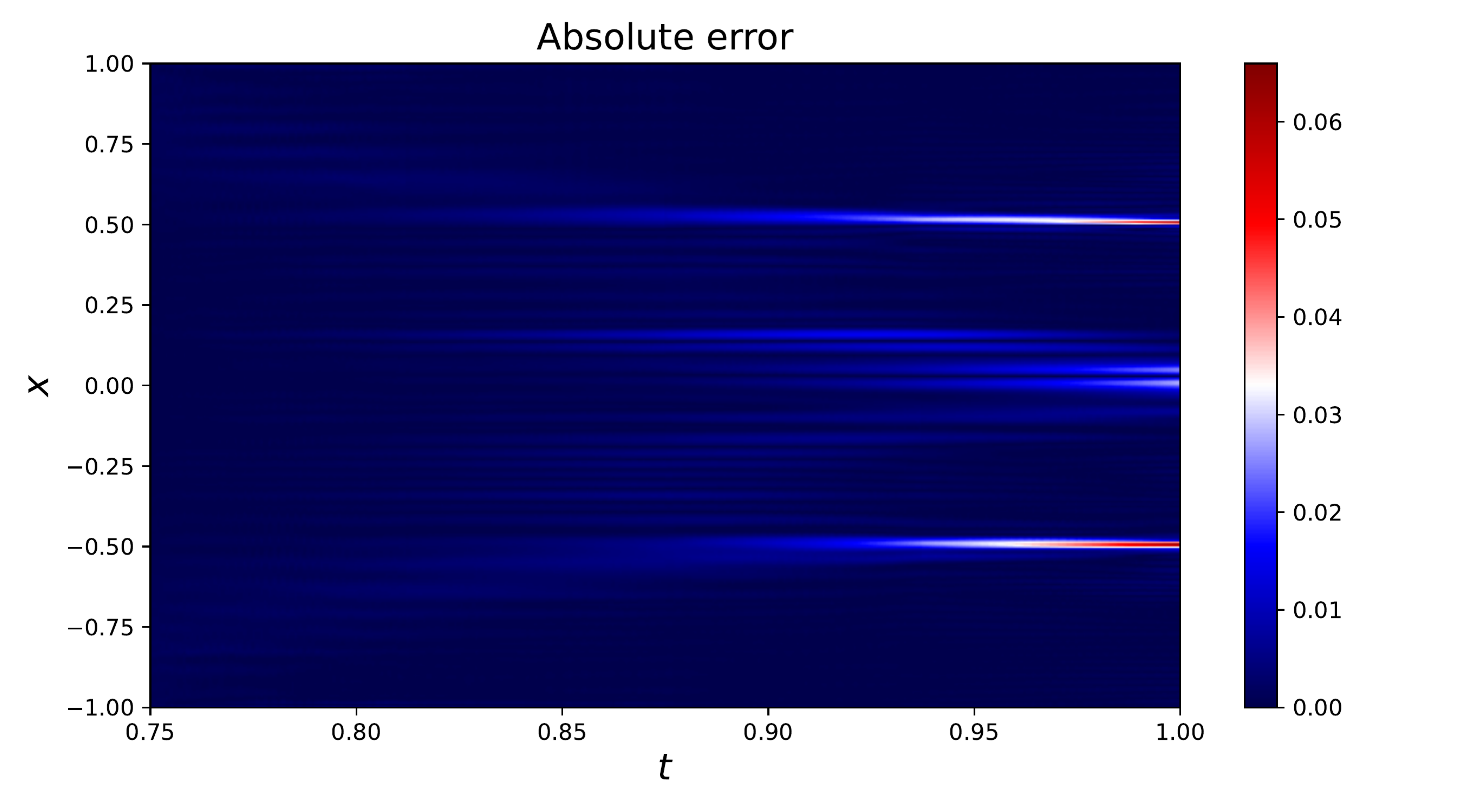}
	\caption{Allen Cahn equation. From left to right: the prediction and absolute error. Relative $L_2$ error is 5.08e-3 (Quadratic basis function, $N_{\mathrm{segment}}$=4, $N$=30).}
	\label{fig:allencahnabsnquadratic}
\end{figure}

In addition, we employ quadratic basis functions and set $N$ to 30. All other parameters remain consistent with the previous configuration. In Table \ref{table_allen_cahn_quadratic} similar results are observed compared to the previous scenario.  Moreover, we also present a representative predicted solution in Figure \ref{fig:allencahnabsnquadratic}, one can see that the predicted solution is in good agreement with the reference solution.

\begin{table}[H]
	\centering
	\begin{tabular}{cccc}
		\toprule
		$N_{\mathrm{segment}}$ & 1                      & 2                       & 4                       \\
		\midrule
		Error                  & 2.25e-01$\pm$ 8.51e-03 & 7.57e-03 $\pm$ 1.68e-03 & 5.08e-03 $\pm$ 1.87e-03 \\
		\bottomrule
	\end{tabular}
	\caption{\centering Allen-Cahn equation. Relative $L_2$ error (mean $\pm$ std) for different segments. (Quadratic finite element basis function with $N=30$)}
	\label{table_allen_cahn_quadratic}
\end{table}

\begin{table}[H]
	\centering
	\resizebox{1.0\linewidth}{!}{
		\begin{tabular}{cccc}
			\toprule
			$N$                              & 10                      & 20                      & 30                      \\
			\midrule
			Galerkin (quadratic)             & 4.07e-03 $\pm$ 3.65e-04 & 7.87e-03 $\pm$ 1.43e-03 & 7.57e-03 $\pm$ 1.68e-03 \\
			Collocation (Hermite polynomial) & 4.12e-03 $\pm$ 2.31e-04 & 5.51e-03 $\pm$ 3.64e-04 & 7.65e-03 $\pm$ 1.04e-03 \\
			Collocation (spline)             & 3.81e-03 $\pm$ 4.79e-04 & 3.89e-03 $\pm$ 5.39e-04 & 3.56e-03 $\pm$ 6.73e-04 \\
			\bottomrule
		\end{tabular}}
	\caption{\centering Allen-Cahn equation. Relative $L_2$ error (mean $\pm$ std) for different approaches ($N_{\mathrm{segment}}=2$)}
	\label{fig:allencahn_galerkin_collocation}
\end{table}

We summarize some relative $L_2$ errors for $N_{\mathrm{segment}}=2$ in
Table \ref{fig:allencahn_galerkin_collocation}. We observe that the collocation framework with piecewise spline functions has achieved the best accuracy.   Moreover, we compare %
the Galerkin projection with causal PINNs in Table \ref{fig:allen_cahn_cpu_time}. Again, we set the same hyper parameters to compare these two methods. The $\epsilon$ in causal PINNs is set to $\{1e-2, 1e-1, 1e0, 1e1, 1e2\}$. It is shown that the accuracy of these two methods is comparable. Note that due to the presence of non-linear term $f(u)$, we must employ a Gaussian quadrature formula of a sufficient degree of exactness for the Galerkin projection in the time direction. Even so, the proposed method runs approximately 40 times faster than causal PINNs.
\begin{table}[!h]
	\centering
	\begin{tabular}{ccc}
		\toprule
		                                 & Galerkin(quadratic) & Causal PINNs \\
		\midrule
		Running time                     & 71                  & 2524         \\
		Error ($N_{\mathrm{segment}}$=1) & 1.95e-02            & 1.33e-02     \\
		Error ($N_{\mathrm{segment}}$=2) & 3.58e-03            & 3.73e-03     \\
		\bottomrule
	\end{tabular}
	\caption{Allen-Cahn equation. Relative $L_2$ error and running time for the proposed method and causal PINNs.}
	\label{fig:allen_cahn_cpu_time}
\end{table}

\subsection{Heat equation on complex geometry}
\begin{figure}[H]
	\centering
	\includegraphics[width=0.9\linewidth]{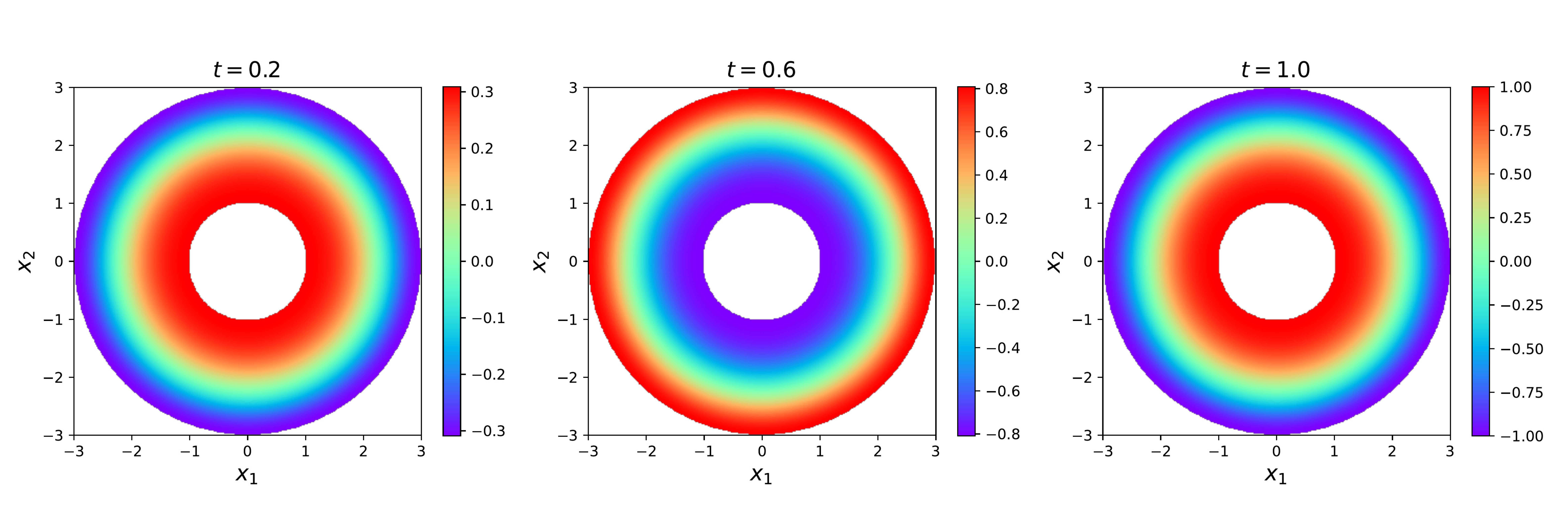}
	\caption{Heat equation on complex geometry. Exact solutions for different time. From left to right: $t = 0.2, 0.6$ and $1.0$.}
	\label{fig:exactannulus}
\end{figure}

In this example, we aim to evaluate the performance of the proposed method on a multi-connected domain, specifically an annulus domain, to assess its capability in handling complex topological structures. We consider the two-dimensional heat equation
\begin{equation}
	\begin{aligned}
		 & u_t - \Delta u = f, \quad  \mbox{in } \Omega \times (0, 1],                          \\
		 & u(\bm{x}, 0) = \sin\left(\frac{\pi}{2}(x_1^2+x_2^2)\right), \quad \mbox{in } \Omega, \\
		 & \frac{\partial u}{\partial \bm{n}} = 0, \quad \mbox{on } \partial\Omega\times[0, 1],
	\end{aligned}
\end{equation}
in an annulus domain $\Omega = \{\bm{x}\in \mathbb{R}^d: 1 < |\bm{x}| < 3\}$. Here the heat source is set to be
\begin{equation}
	\small
	f(x_1, x_2, t) = -8\pi \sin(8\pi t) \sin\left(\frac{\pi}{8} (x_1^2 + x_2^2 + 3) \right) - \frac{\pi}{2} \cos(8\pi t) \cos\left(\frac{\pi}{8} (x_1^2 + x_2^2 + 3) \right) + (\frac{\pi}{4})^2 (x_1^2 + x_2^2) \cos(8\pi t) \sin\left(\frac{\pi}{8} (x_1^2 + x_2^2 + 3) \right).
\end{equation}
The true solution is given by
\begin{equation}
	u(x_1, x_2, t) = \cos(8\pi t)\sin\left(\frac{\pi}{8}(x_1^2+x_2^2+3)\right),
\end{equation}
and solutions at $t = 0.2, 0.6$ and $1.0$ are plotted in Figure \ref{fig:exactannulus}.
We set $N_{\mathrm{segment}}$ to be 1. The latent coefficients, denoted as $\{\omega_i(\bm{x}; \theta)\}$, are represented using a fully-connected neural network with tanh activation function, $4$ hidden layers and $128$ neurons per hidden layer. We randomly generate $2,000$ collocation points in the domain $\Omega$ for enforcing the PDE residual, $400$ points on the inner boundary of the annulus, $1,200$ points on the outer boundary of the annulus and $2,000$ initial points. We proceed by training the resulting model via full-batch gradient descent using the Adam optimizer for $40,000$ iterations. The corresponding results from $5$ independent runs with different network initializations for different numbers of mesh elements are summarized in Table \ref{table:annulus_error}. Furthermore, the resulting relative $L_2$ error (mean $\pm$ std) is $2.10e+00 \pm 1.91e-01$, which exhibited significantly lower accuracy compared to the proposed method.
\begin{table}[H]
	\centering
	\resizebox{1.0\linewidth}{!}{
		\begin{tabular}{ccccc}
			\toprule
			N                                & 10                      & 20                      & 40                      & 80                      \\
			\midrule
			Galerkin (linear)                & 1.27e+00 $\pm$ 1.27e-03 & 1.44e-01 $\pm$ 5.39e-04 & 3.53e-02 $\pm$ 9.29e-04 & 2.04e-02 $\pm$ 9.68e-04 \\
			\midrule
			Galerkin (quadratic)             & 9.82e-02 $\pm$ 1.26e-04 & 2.96e-02 $\pm$ 9.75e-04 & 2.01e-02 $\pm$ 8.27e-04 & 2.50e-02 $\pm$ 9.78e-04 \\
			\midrule
			Collocation (Hermite polynomial) & 2.52e-02 $\pm$ 1.96e-03 & 9.53e-03 $\pm$ 6.13e-03 & 9.23e-03 $\pm$ 4.87e-03 & 1.44e-02 $\pm$ 9.68e-03 \\
			\midrule
			Collocation (spline)             & 5.94e-01 $\pm$ 4.71e-03 & 1.19e-02 $\pm$ 4.53e-03 & 6.64e-03 $\pm$ 2.50e-03 & 9.54e-03 $\pm$ 2.14e-03 \\
			\bottomrule
		\end{tabular}
	}
	\caption{\centering Heat equation on complex geometry. Relative $L_2$ errors (mean $\pm$ std) at different N.}
	\label{table:annulus_error}
\end{table}

\subsection{High-dimensional linear problem}
\label{high-dimensional-linear-problem}
To demonstrate the effectiveness of the proposed approaches in tackling high-dimensional PDEs, we consider the following parabolic equation with non-constant coefficients
\begin{equation}
	\begin{aligned}
		 & u_t - \nabla \cdot (a(\bm{x})\nabla u)=f, \quad \mbox{in } \Omega\times (0,1], \\
		 & u(\bm{x},0)=0,\quad \mbox{in }\Omega,                                          \\
		 & u = 0,\quad \mbox{on }\partial \Omega\times [0,1],
	\end{aligned}
\end{equation}
with $a(\bm{x})=1+ |\bm{x}|^2/2$. The domain is set to be the 20-D unit ball $\Omega = \{\bm{x}\in\mathbb{R}^{20}\mid |\bm{x}|<1\}$, and the true solution is set to be
\begin{equation}
	u(\bm{x},t) = \sin \big(\sin (2\pi \omega t)(|\bm{x}|^2 - 1)\big).
\end{equation}
We set $\omega$ to be 3 and $N_{\mathrm{segment}}$ to be 1. We represent the latent coefficients $\{{\omega}_i(\bm{x};\theta)\}$ by a fully-connected neural network with tanh activation function, 5 hidden layers, and 128 neurons per hidden layer. We impose exactly the Dirichlet boundary conditions by transforming the output into the following form:
\begin{equation}
	\tilde{{\omega}}_i(\bm{x};\theta) = {\omega}_{i}(\bm{x};\theta) (\vert \bm{x}\vert -1).
\end{equation}
In addition, we set $\omega_0(\bm{x};\theta)$ to be zero such that the initial conditions are exactly satisfied. To obtain a set of training data for evaluating PDE residual, we randomly sample 100,000 collocation points in $\Omega$. Since the problem is linear, equation \eqref{NN_linear_finite_element_martix_form} is used to define the loss. We set the size of mini-batch to 5,000 and train the model via mini-batch stochastic descent with the Adam optimizer for 40,000 iterations. The corresponding results for different numbers of mesh elements are summarized in Table \ref{high-linear-result}. We observe that the resulting relative $L_2$
error is $9.67e-4$, which is more accurate than the ones in recent work \cite{gu2022deep} ($3.07e-3$ for the deep adaptive basis Galerkin approach and $4.54e-2$ for the PINNs).

\begin{table}[H]
	\centering
	\resizebox{1.0\linewidth}{!}{
		\begin{tabular}{ccccc}
			\toprule
			N                                & 10                      & 20                      & 40                      & 80                      \\
			\midrule
			Galerkin (linear)                & 1.04e+00 $\pm$ 3.11e-05 & 2.48e-01 $\pm$ 7.97e-05 & 5.75e-02 $\pm$ 1.91e-04 & 1.48e-02 $\pm$ 2.51e-04 \\
			Galerkin (quadratic)             & 4.45e-02 $\pm$ 1.25e-04 & 7.61e-03 $\pm$ 2.73e-05 & 1.74e-03 $\pm$ 2.07e-05 & 9.67e-04 $\pm$ 2.03e-05 \\
			Collocation (Hermite polynomial) & 2.28e-02 $\pm$ 3.28e-06 & 1.41e-03 $\pm$ 1.87e-05 & 6.58e-04 $\pm$ 3.48e-05 & 1.04e-03 $\pm$ 8.73e-05 \\
			Collocation (spline)             & 4.89e-02 $\pm$ 6.66e-05 & 1.42e-02 $\pm$ 3.44e-05 & 6.42e-04 $\pm$ 3.51e-05 & 8.89e-04 $\pm$ 1.07e-04 \\
			\bottomrule
		\end{tabular}
	}
	\caption{\centering {High-dimensional linear problem. Relative $L_2$ error (mean $\pm$ std) for different $N$.}}
	\label{high-linear-result}
\end{table}

\subsection{High-dimensional nonlinear problem}
In this case, we solve the Allen-Cahn equation
\begin{equation}
	\begin{aligned}
		 & u_t - \Delta u + u^3 - u = f_{AC}, \quad \mbox{in }\Omega\times [0,1], \\
		 & u(\bm{x},0) = 0, \quad \mbox{in }\Omega,                               \\
		 & u = 0, \quad \mbox{on } \partial \Omega \times [0,1].
	\end{aligned}
	\label{high_AC_example}
\end{equation}
The domain is set as $\Omega = \{\bm{x}\in \mathbb{R}^{20}: |\bm{x}|<1\}$, and the true solution is set to be
\begin{equation*}
	u(\bm{x},t) = \sin(\sin(2\pi \omega t)(|\bm{x}|^2 -1)).
\end{equation*}
Here we set $\omega$ to be 3 and $N_{\mathrm{segment}}$ to be 1. We represent the latent coefficients $\{\omega_i(\bm{x};\theta)\}$ by a fully-connected neural network with tanh activation function, 5 hidden layers and 128 neurons per hidden layer. The initial condition and Dirichlet boundary condition can be exactly embedded into the network structure (as described in Section \ref{high-dimensional-linear-problem}). We set $K$ in equation \eqref{eqn:grid_s} to be 30 and randomly sample 100,000 collocation points for evaluating PDE residual. We set mini-batch to be 5,000 and proceed by training the resulting model via stochastic gradient descent using the Adam optimizer for 40,000 iterations. The obtained errors for different number of finite element $N$ are shown in Table \ref{high-nonlinear-result}. The resulting relative $L_2$ error is $1.16e-3$, which is more accurate than the recent ones in \cite{gu2022deep} ($2.07e-3$ for the deep adaptive basis Galerkin approach and $7.79e-2$ for the PINNs).

\begin{table}[H]
	\centering
	\resizebox{1.0\linewidth}{!}{
		\begin{tabular}{ccccc}
			\toprule
			N                    & 10                      & 20                      & 40                      & 80                      \\
			\midrule
			Galerkin (quadratic) & 4.48e-02 $\pm$ 5.34e-05 & 7.92e-03 $\pm$ 7.63e-05 & 1.89e-03 $\pm$ 3.21e-04 & 1.16e-03 $\pm$ 3.57e-04 \\			Collocation (Hermite polynomial) & 2.48e-02 $\pm$ 6.38e-05 & 2.21e-03 $\pm$ 5.92e-04 & 1.72e-03  $\pm$ 4.46e-04& 1.33e-03 $\pm$ 5.31e-04
			\\
			Collocation (spline) & 5.04e-02 $\pm$ 6.21e-06 & 1.61e-02 $\pm$ 9.75e-05 & 2.42e-03 $\pm$ 3.18e-04 & 2.31e-03 $\pm$ 2.36e-04
			\\
			\bottomrule
		\end{tabular}
	}
	\caption{\centering {High-dimensional nonlinear problem. Relative $L_2$ error (mean $\pm$ std) for different $N$.}}
	\label{high-nonlinear-result}
\end{table}

\subsection{Low regularity test problem}
\begin{figure}[H]
	\centering
	\includegraphics[width=0.9\linewidth]{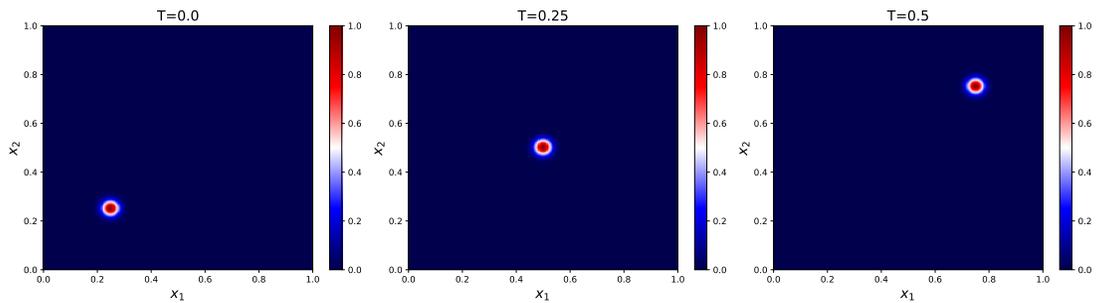}
	\caption{Low regularity test problem. Exact solutions for different time ($\beta$=1000). From left to right: $t=0,0.25$ and $0.5$.}
	\label{fig:exactsolution1000}
\end{figure}

Our last example aims to demonstrate the effectiveness of the proposed adaptivity strategy. We consider the following low-regularity equation
\begin{equation*}
	\begin{aligned}
		 & u_t - \Delta u =f\quad \mbox{in }\Omega\times (0,T]                         \\
		 & u(x_1, x_2, 0) = g(x_1, x_2)\quad \mbox{in } \Omega\times \{0\},            \\
		 & u(x_1, x_2, t) = h(x_1,x_2, t)\quad \mbox{on } \partial \Omega\times (0,T].
	\end{aligned}
\end{equation*}
where $\Omega= [0, 1]^2, T=0.5$ and the true solution is chosen as
\begin{equation*}
	u(x_1,x_2,t)=\exp\big(-\beta\left[(x_1 - t - 1/4)^2 + (x_2 - t - 1/4)^2\right]\big).
\end{equation*}
This solution has one peak at the point $(t+1/4, t+1/4)$ and decreases rapidly away from $(t+1/4,t+1/4)$, as shown in Figure \ref{fig:exactsolution1000}.

We first consider the case with $\beta=200$. We represent the latent coefficients
$\{\tilde{\omega}_i(\bm{x};\theta)\}$ by a fully-connected neural network with
tanh activation function, 5 hidden layers and 128 neurons per hidden layer.
To simplify the training objective, we precisely embed the Dirichlet
boundary conditions into the neural network architecture through the
following transformation
\begin{equation*}
	\tilde{\omega}_i(\bm{x};\theta) =  x_1x_2(1-x_1)(1-x_2)\omega_i(\bm{x};\theta).
\end{equation*}
For B-KRnet, we take 8 CDF coupling layers, and two fully connected layers with 32
neurons for each CDF coupling layer. The corresponding activation function is tanh. To assess the effectiveness of our adaptive sampling
approach, we generate a uniform meshgrid with size $256\times 256\times
	100$ in $[0,1]^2\times [0,0.5]$  and compute relative $L_2$ error on these grid points.

For collocation projection, we set the number of mesh elements to 20, and randomly sample 1,000 collocation points in $\Omega$ as our initial training set. We let $N_{\mathrm{adaptive}}=5$ and
$N_{\mathrm{new}}=500$. The number of epochs for training
$u(\bm{x},t;\theta)$ and $p_{\theta_f}(\bm{x},t)$ is set to 40,000 and
5,000, respectively.

\begin{figure}[H]
	\centering
	\includegraphics[width=0.4\linewidth]{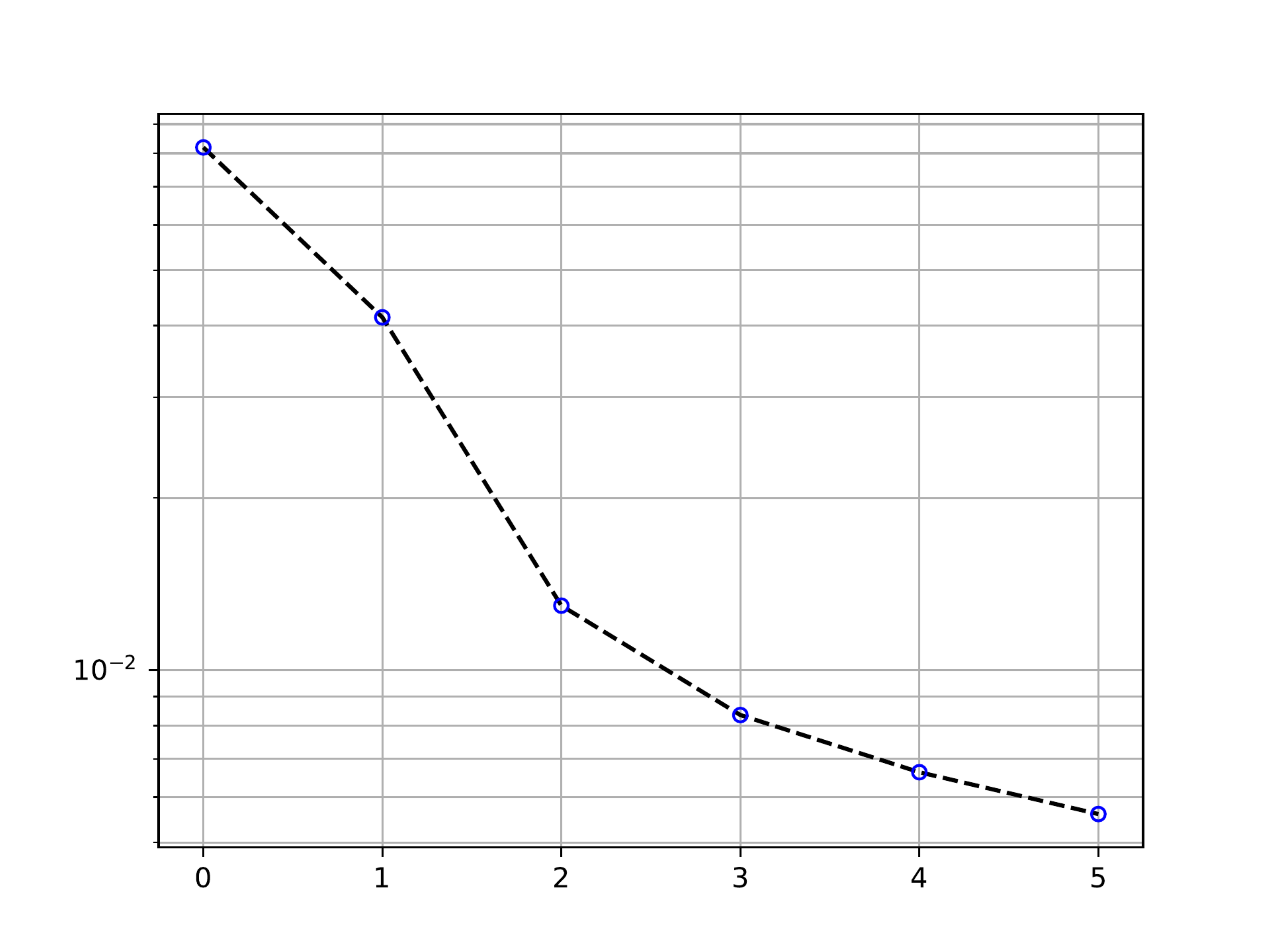}
	\caption{Low regularity test problem. Relative $L_2$ error of different adaptivity iterations. ($\beta = 200$, collocation projection)}
	\label{fig:collocation_error_beta_200}
\end{figure}

\begin{figure}[H]
	\centering
	\includegraphics[width=0.3\linewidth]{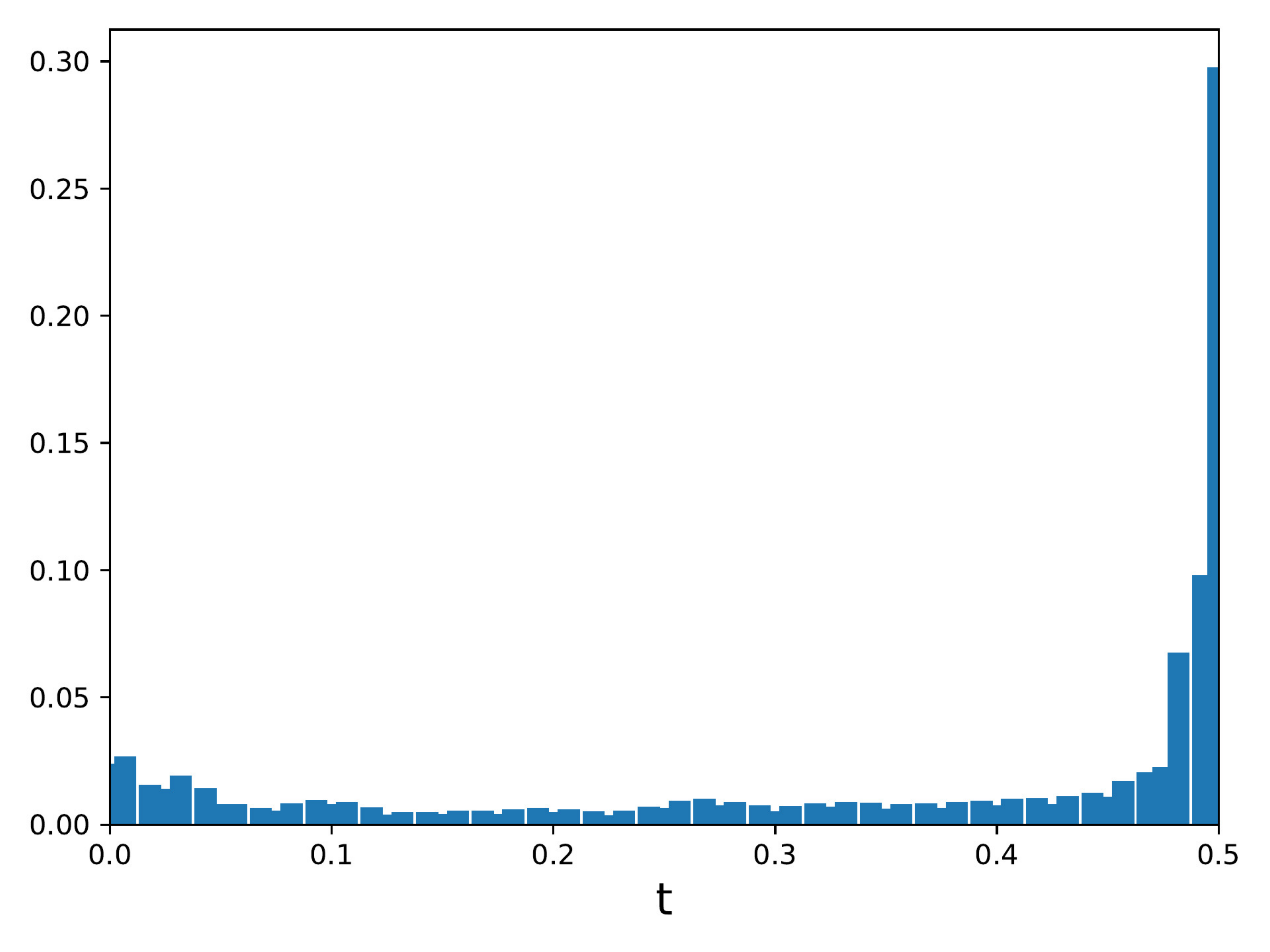}
	\includegraphics[width=0.3\linewidth]{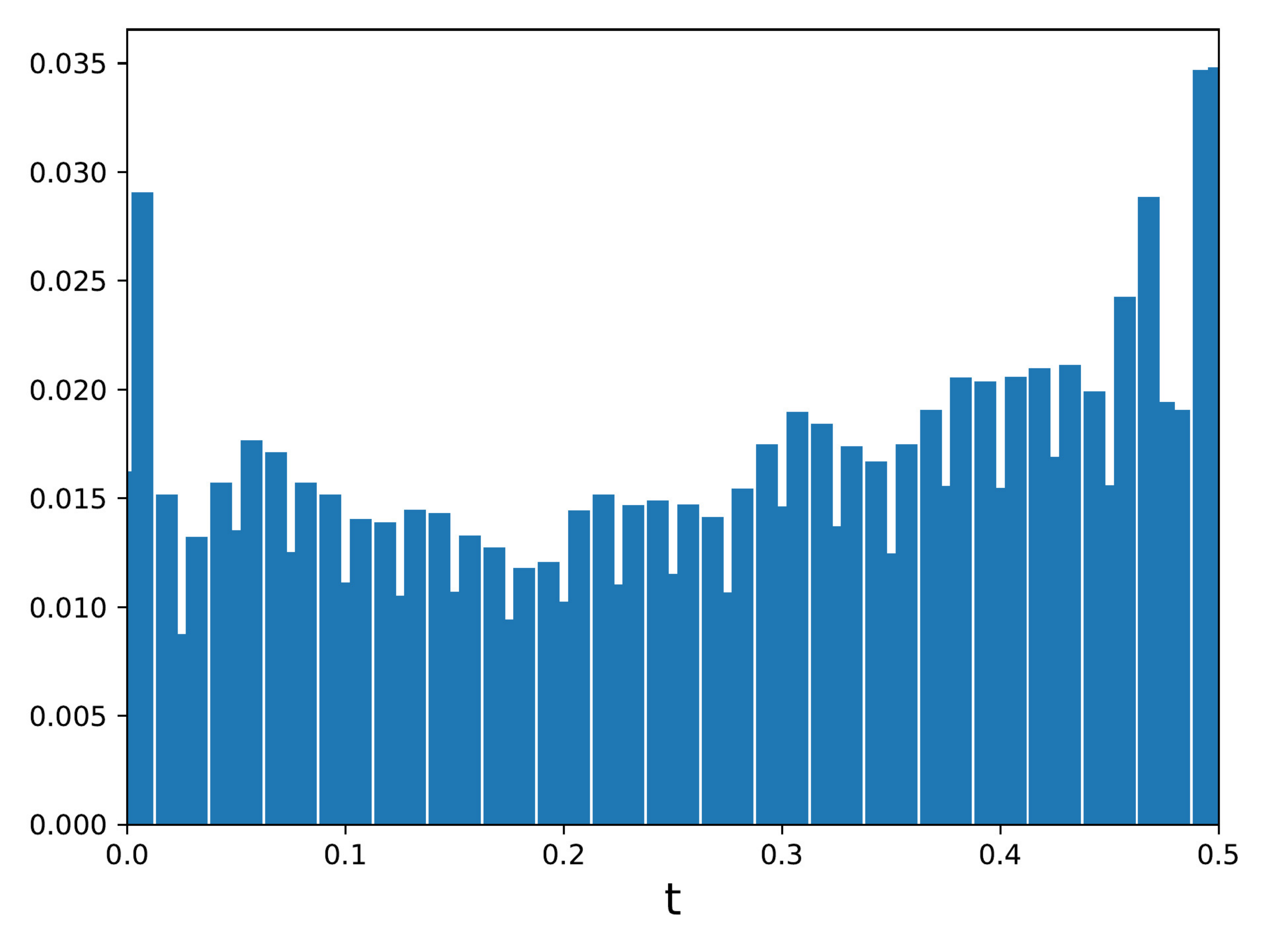}
	\includegraphics[width=0.3\linewidth]{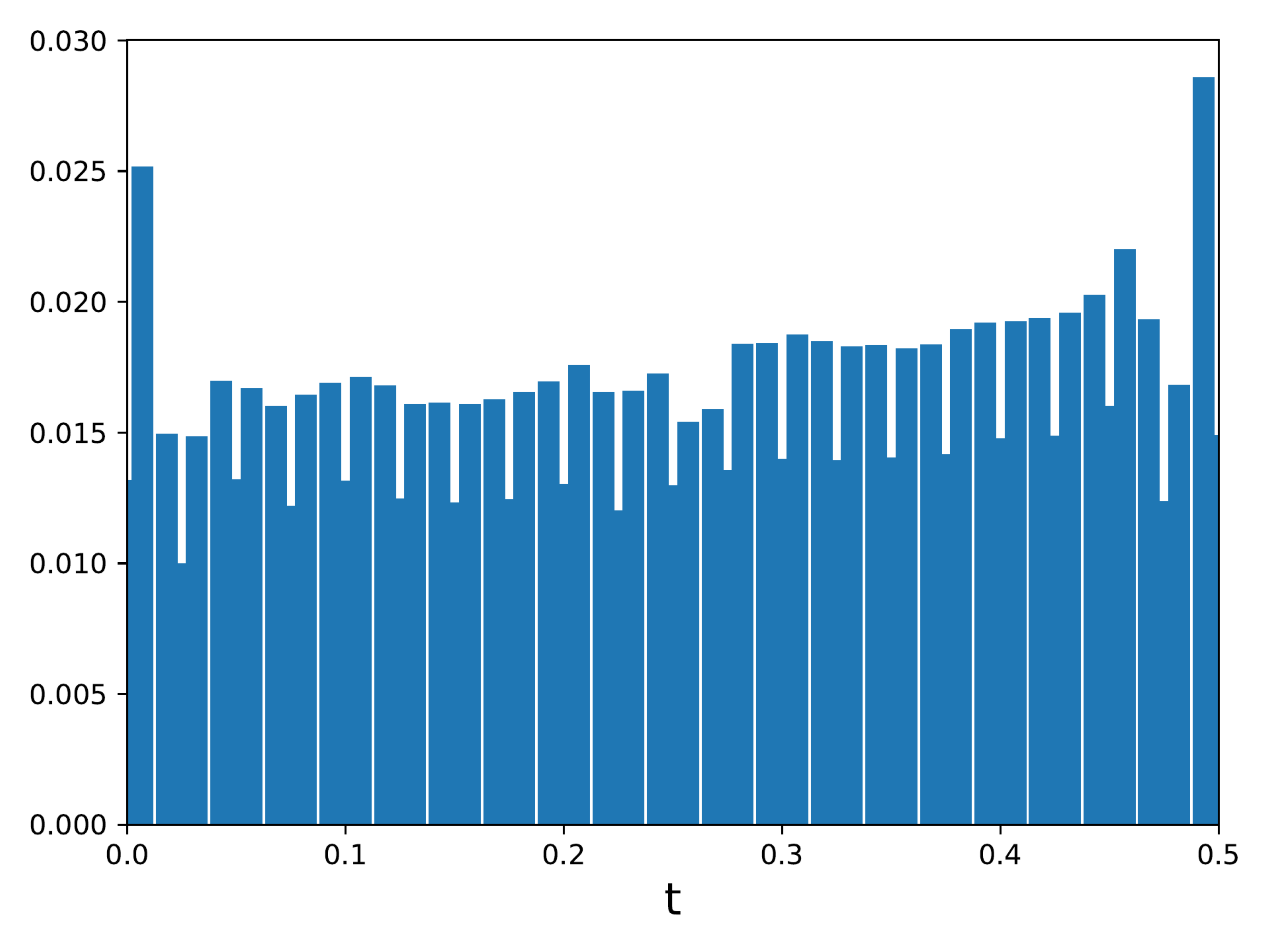}
	\caption{Low regularity test problem. Temporal density of the 1st, 3rd, 5th adaptivity iteration. ($\beta=200$, collocation projection)}
	\label{fig:collocation_temporal_density_beta_200}
\end{figure}

\begin{figure}[H]
	\centering
	\includegraphics[width=0.3\linewidth]{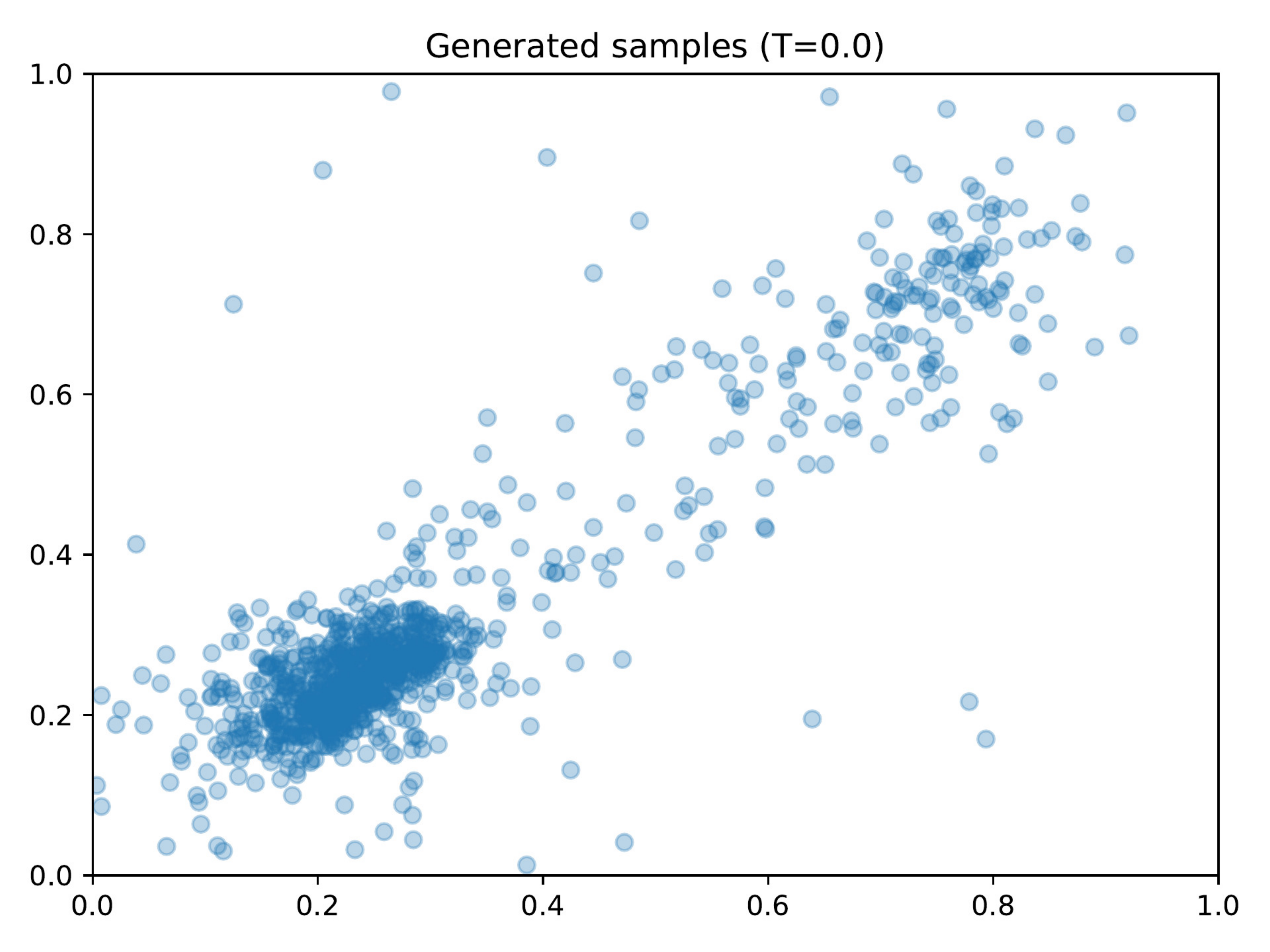}
	\includegraphics[width=0.3\linewidth]{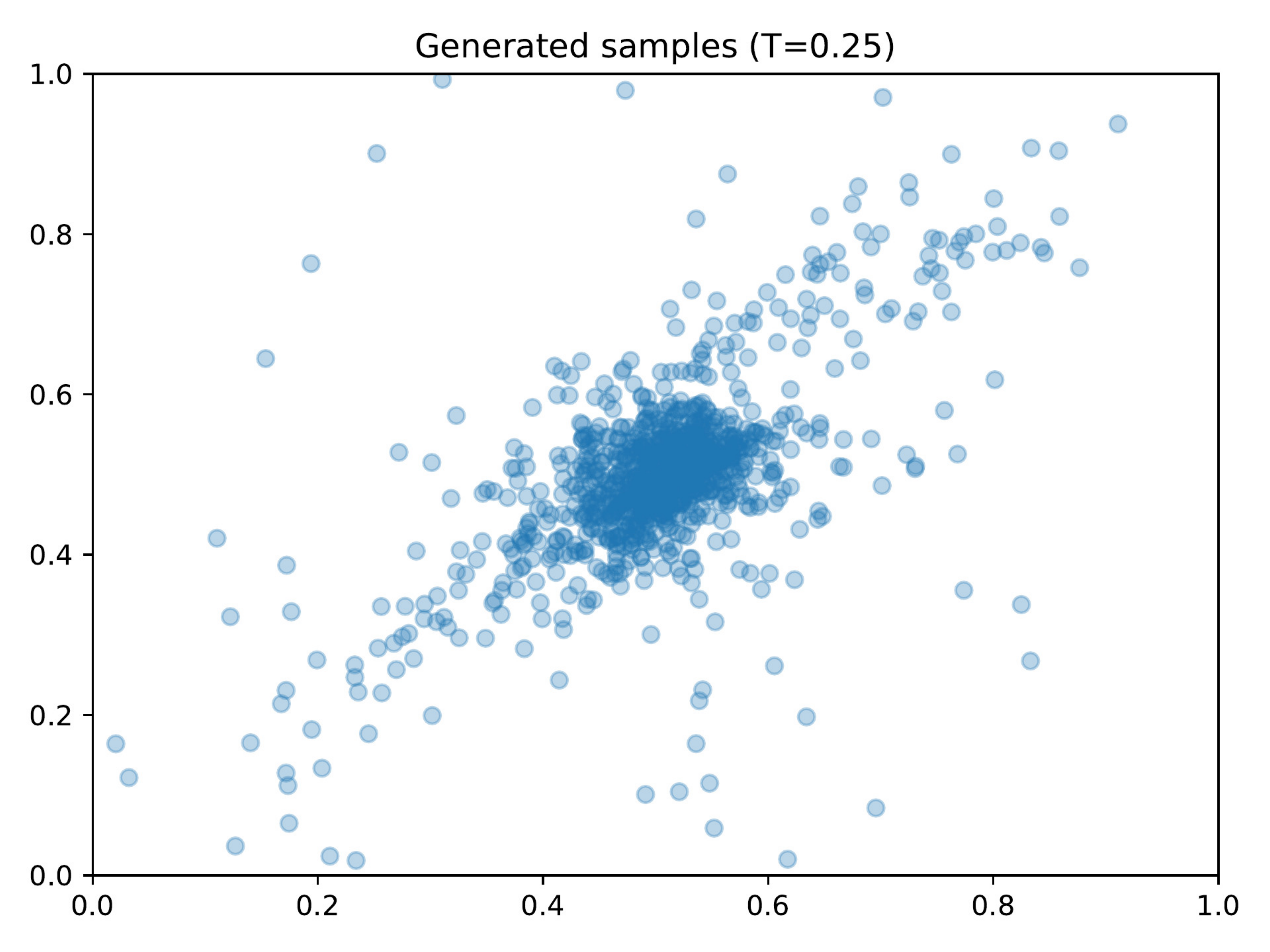}
	\includegraphics[width=0.3\linewidth]{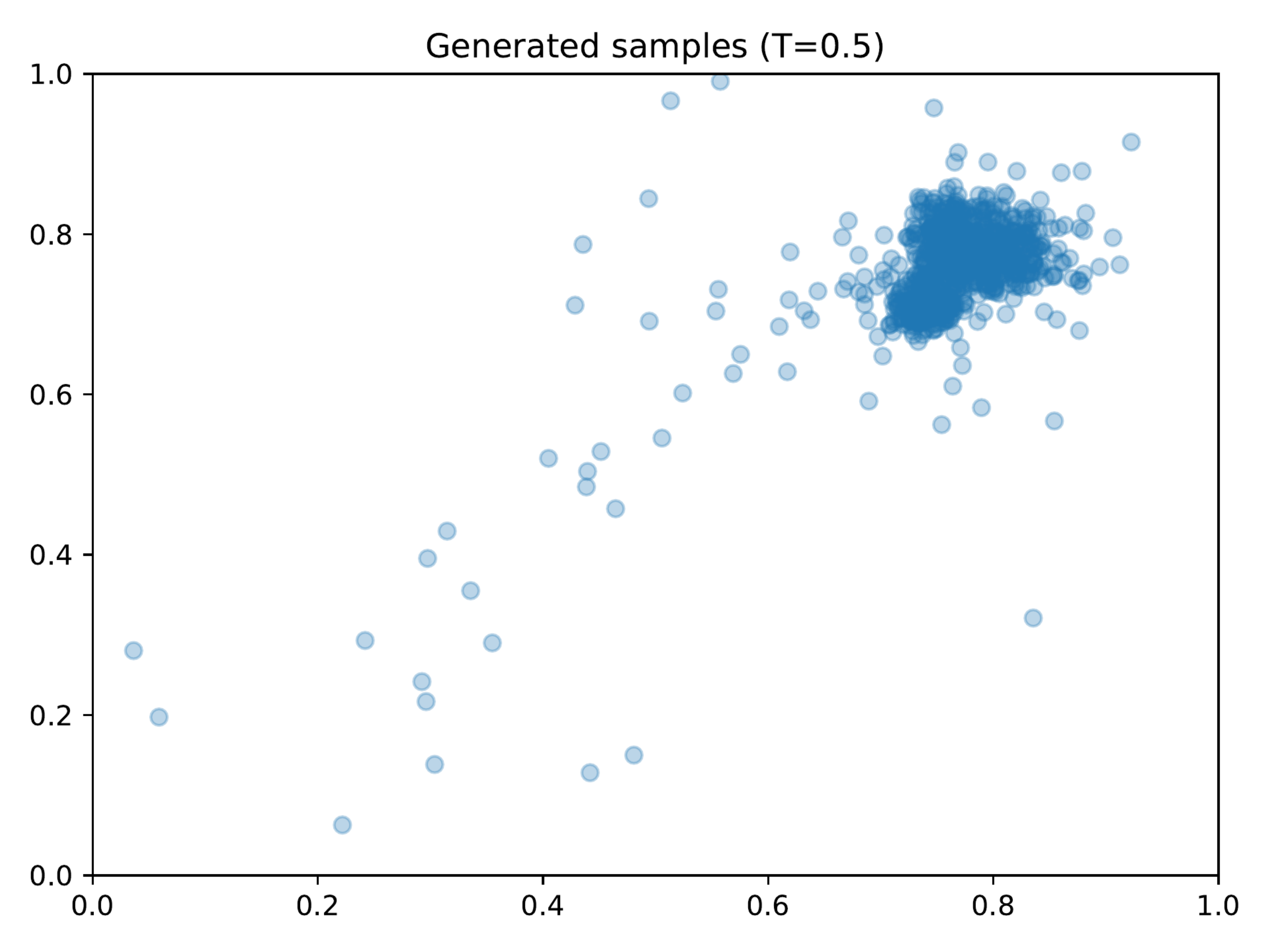}

	\includegraphics[width=0.3\linewidth]{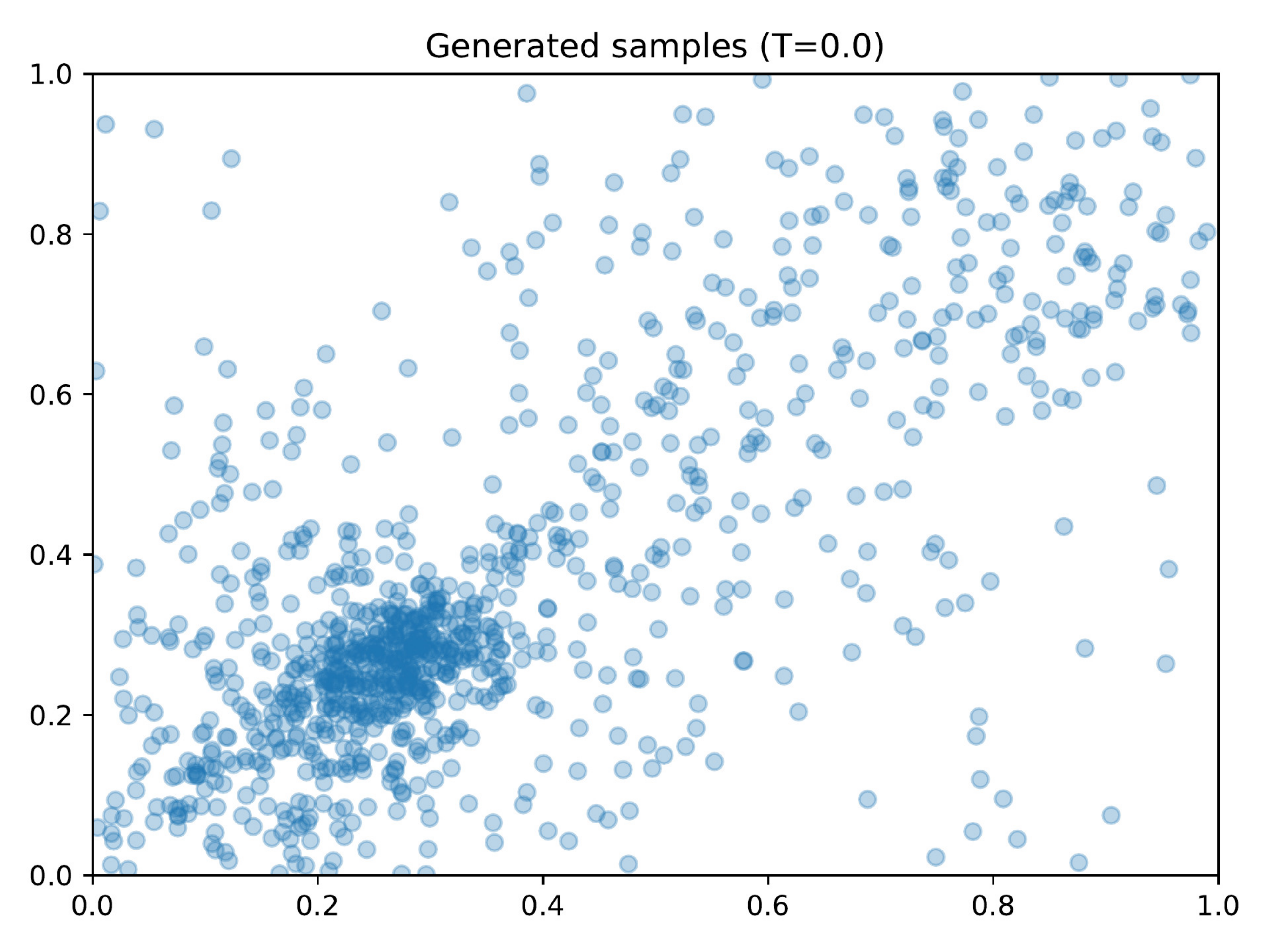}
	\includegraphics[width=0.3\linewidth]{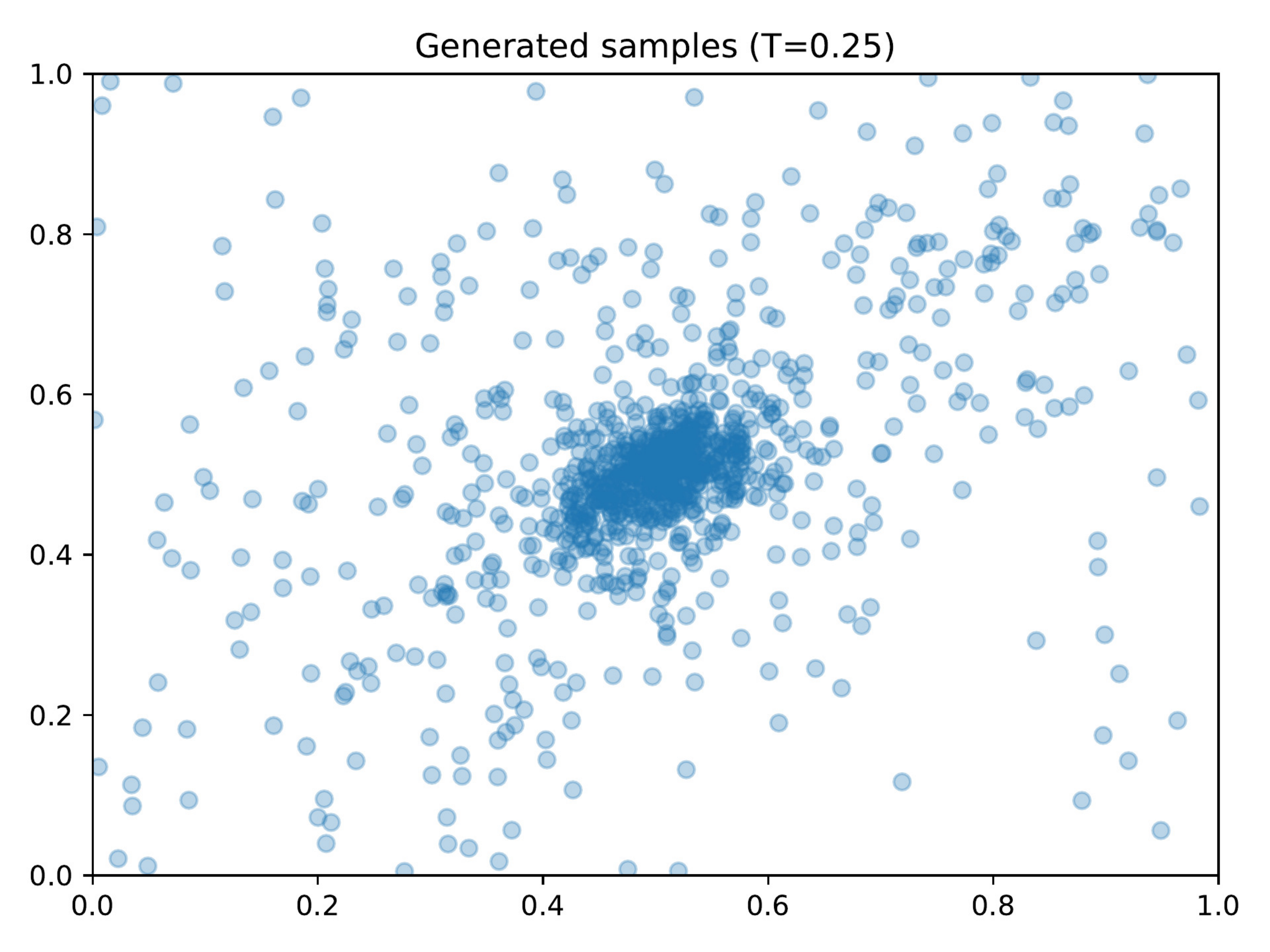}
	\includegraphics[width=0.3\linewidth]{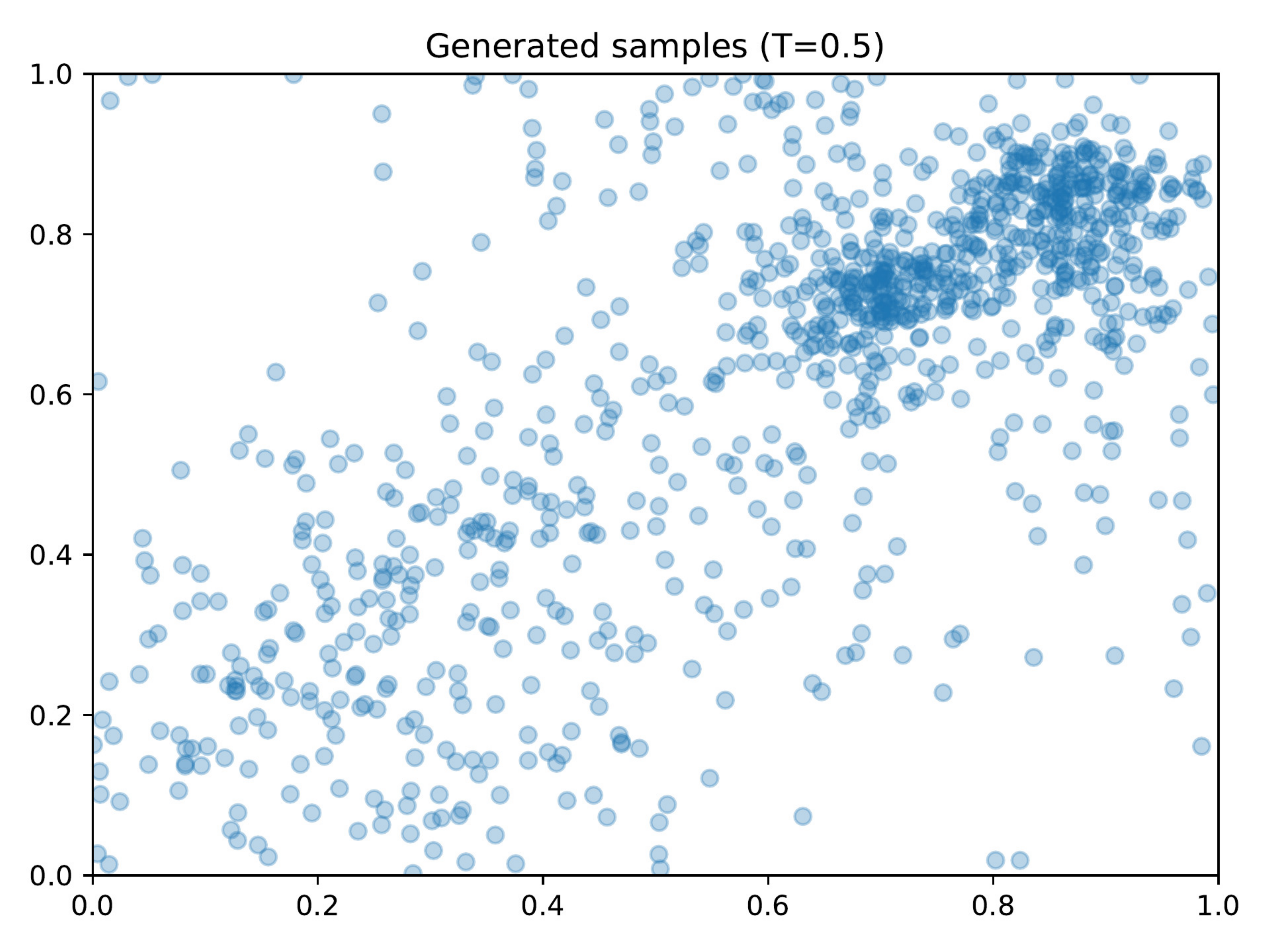}
	\caption{Low regularity test problem. From top to bottom: The generated samples of the 1st and 5th adaptivity iteration. ($\beta=200$, collocation projection)}
	\label{fig:collocation_generated_samples_beta_200}
\end{figure}

In Figure \ref{fig:collocation_error_beta_200}, we plot the approximation errors with respect to the adaptivity iterations. It is clearly seen that the error rapidly decreases as the adaptivity iteration step $k$ increases.
Figure \ref{fig:collocation_temporal_density_beta_200} shows the evolution of temporal distribution $p_{\mathrm{dis}}(t)$
for $k=1,3,5$. It is seen that the largest density at the first adaptivity iteration is around terminal time. After the training set is refined,
the error profile becomes more flat at $k=3$ and $5$. Similar results are observed in Figure \ref{fig:collocation_generated_samples_beta_200}, which indicates that the
distribution of residual becomes more uniform as adaptivity iteration increases.

For Galerkin projection, we set the number of piecewise quadratic basis functions
to 20, the remaining hyper-parameter settings remain the same. The corresponding numerical results, included in \ref{app:adaptivity_figure}, are very similar to the results for collocation projection.

\begin{figure}[H]
	\centering
	\includegraphics[width=0.3\linewidth]{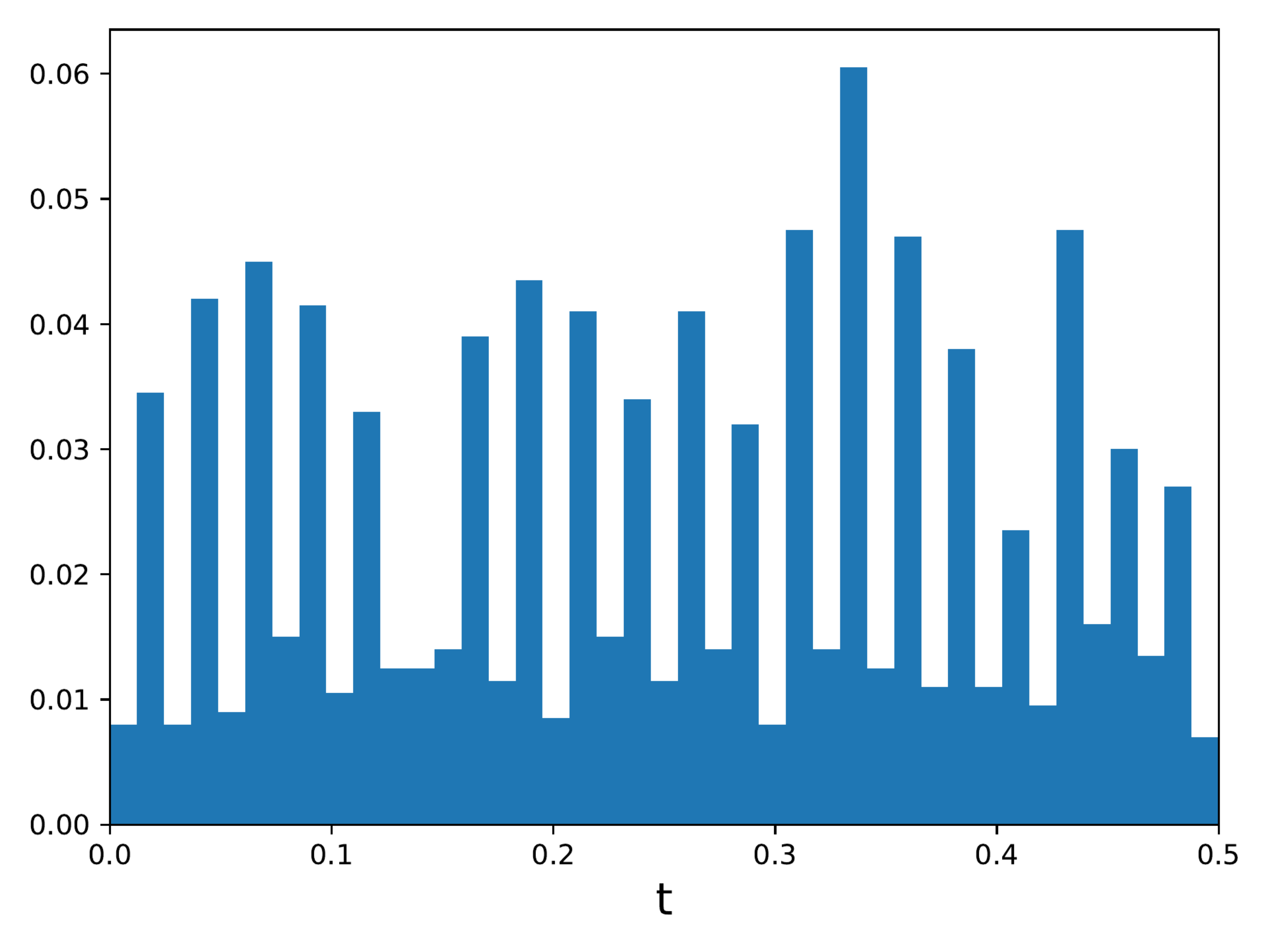}
	\includegraphics[width=0.3\linewidth]{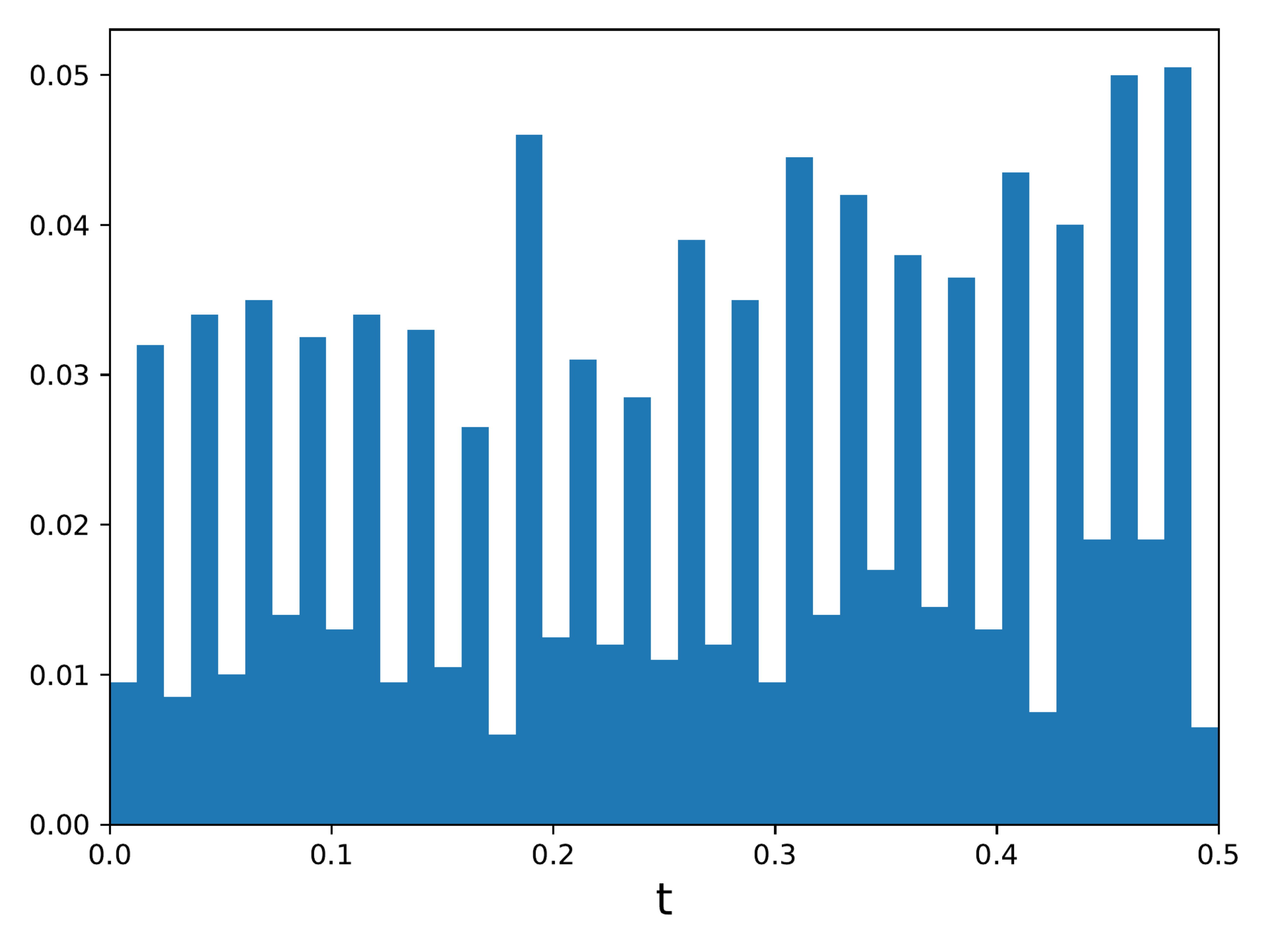}
	\includegraphics[width=0.3\linewidth]{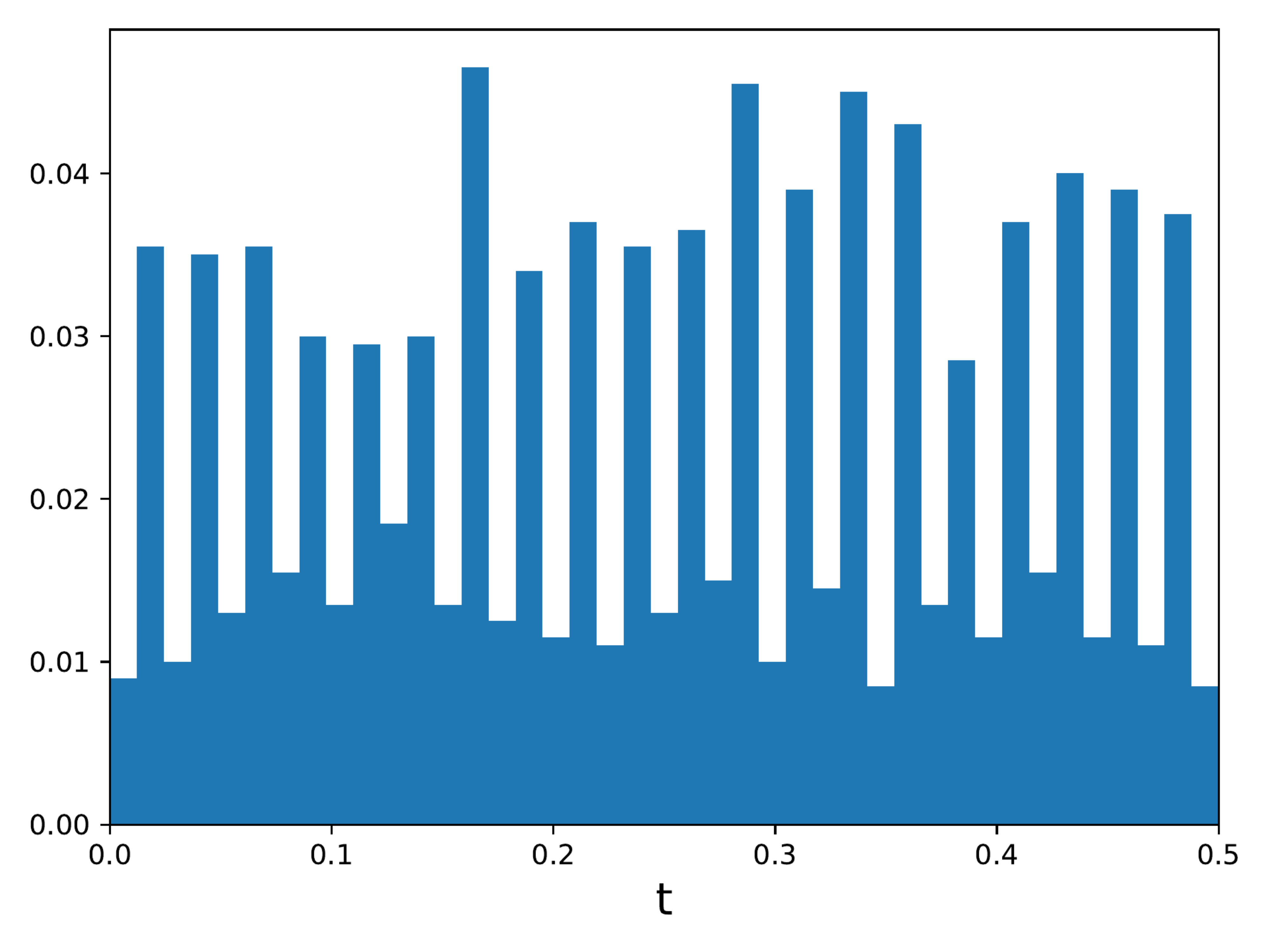}

	\caption{Low regularity test problem. Temporal density of the 1st, 3rd, 5th adaptivity iteration. ($\beta=1000$, Galerkin projection)}
	\label{fig:galerkin_temporal_weights_beta_1000}
\end{figure}

\begin{figure}[H]
	\centering
	\includegraphics[width=0.3\linewidth]{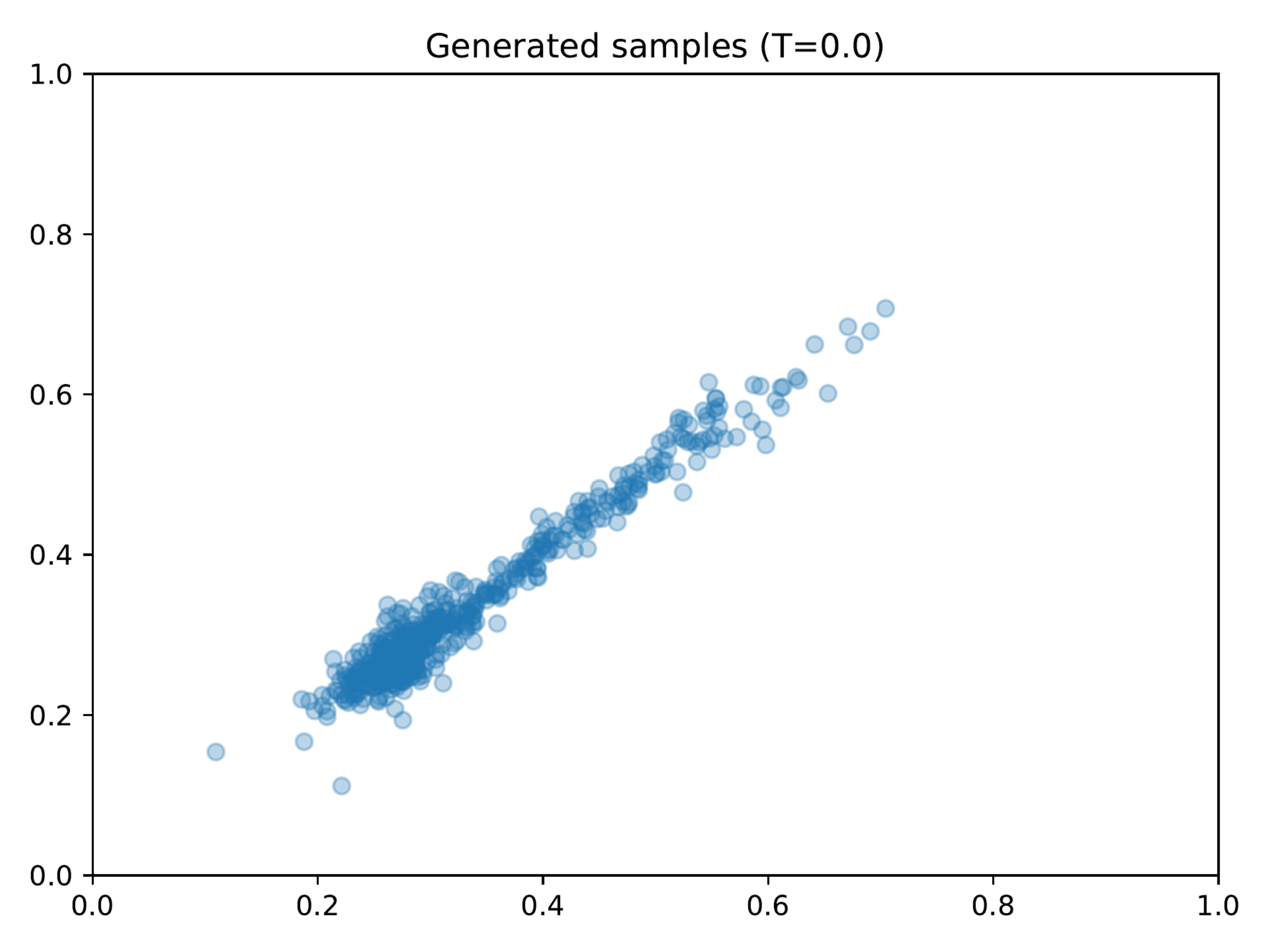}
	\includegraphics[width=0.3\linewidth]{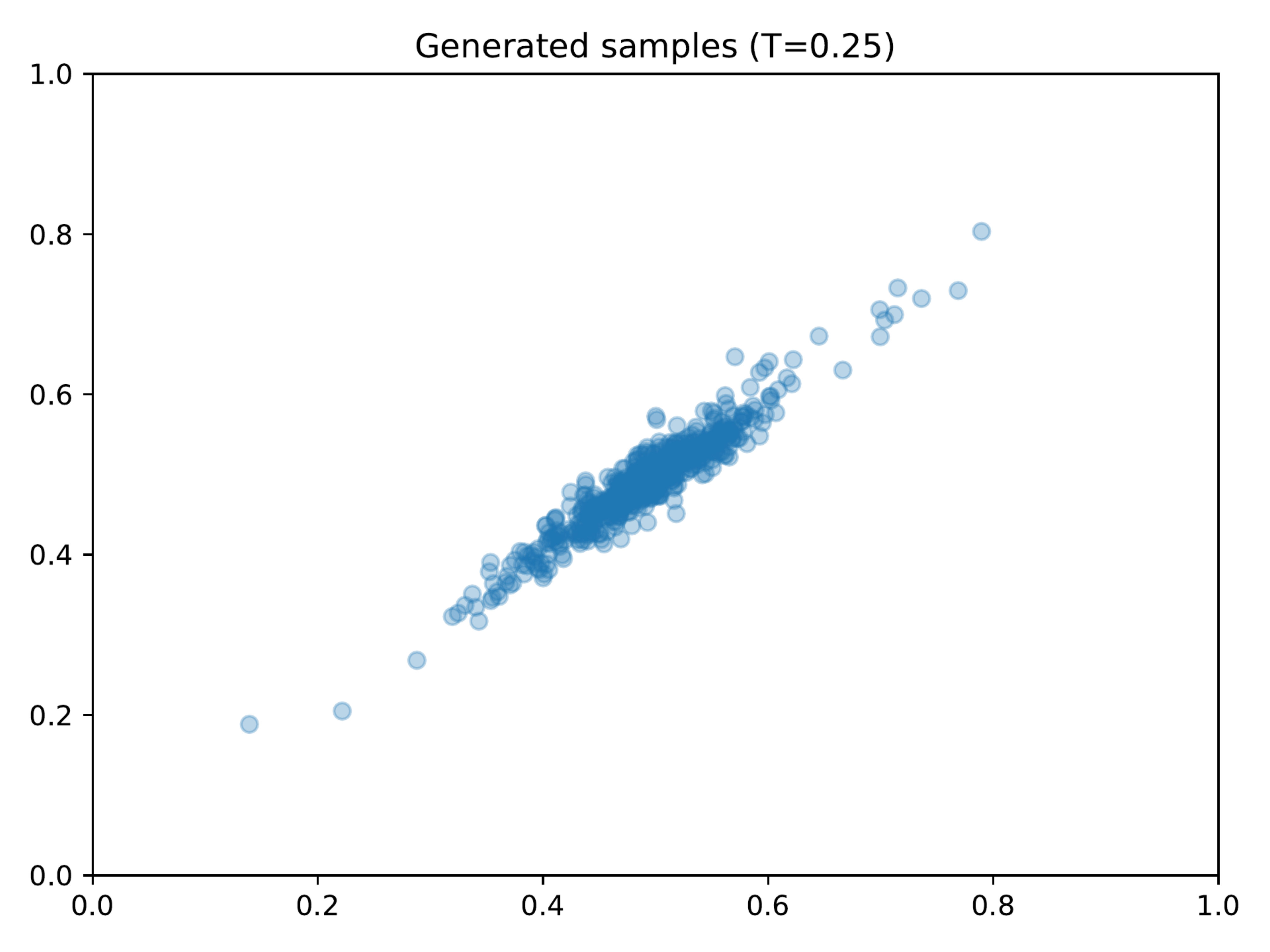}
	\includegraphics[width=0.3\linewidth]{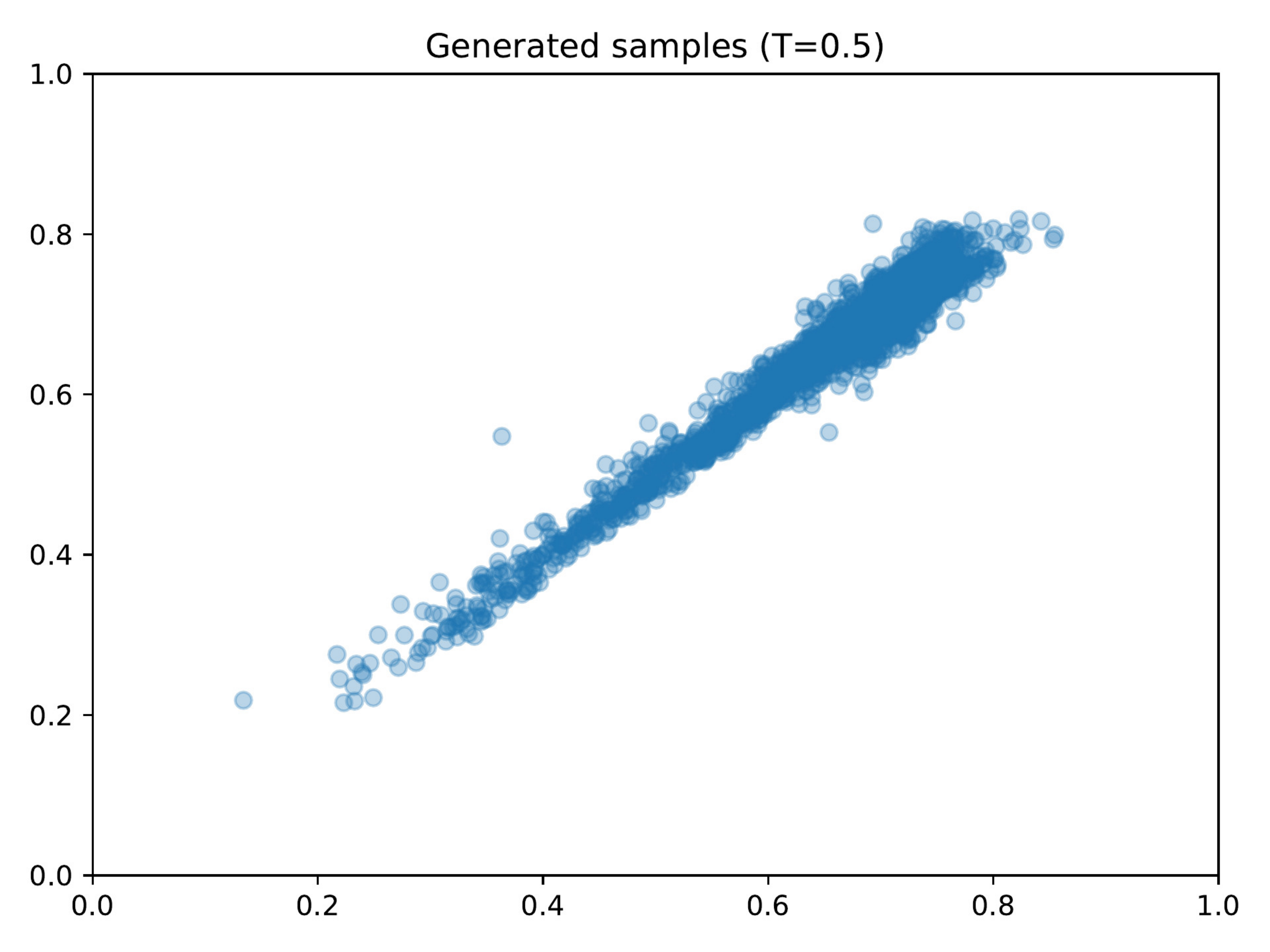}

	\includegraphics[width=0.3\linewidth]{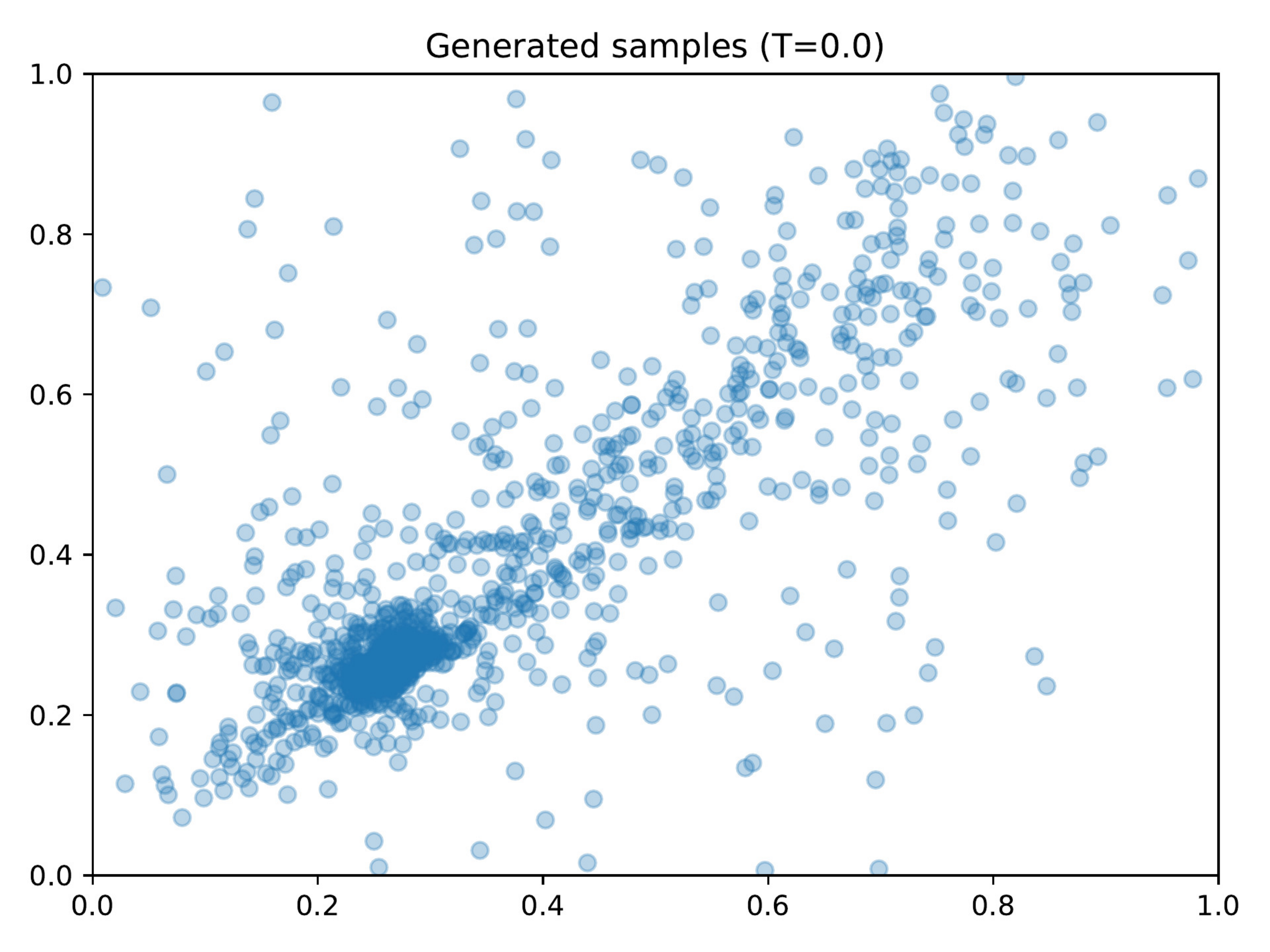}
	\includegraphics[width=0.3\linewidth]{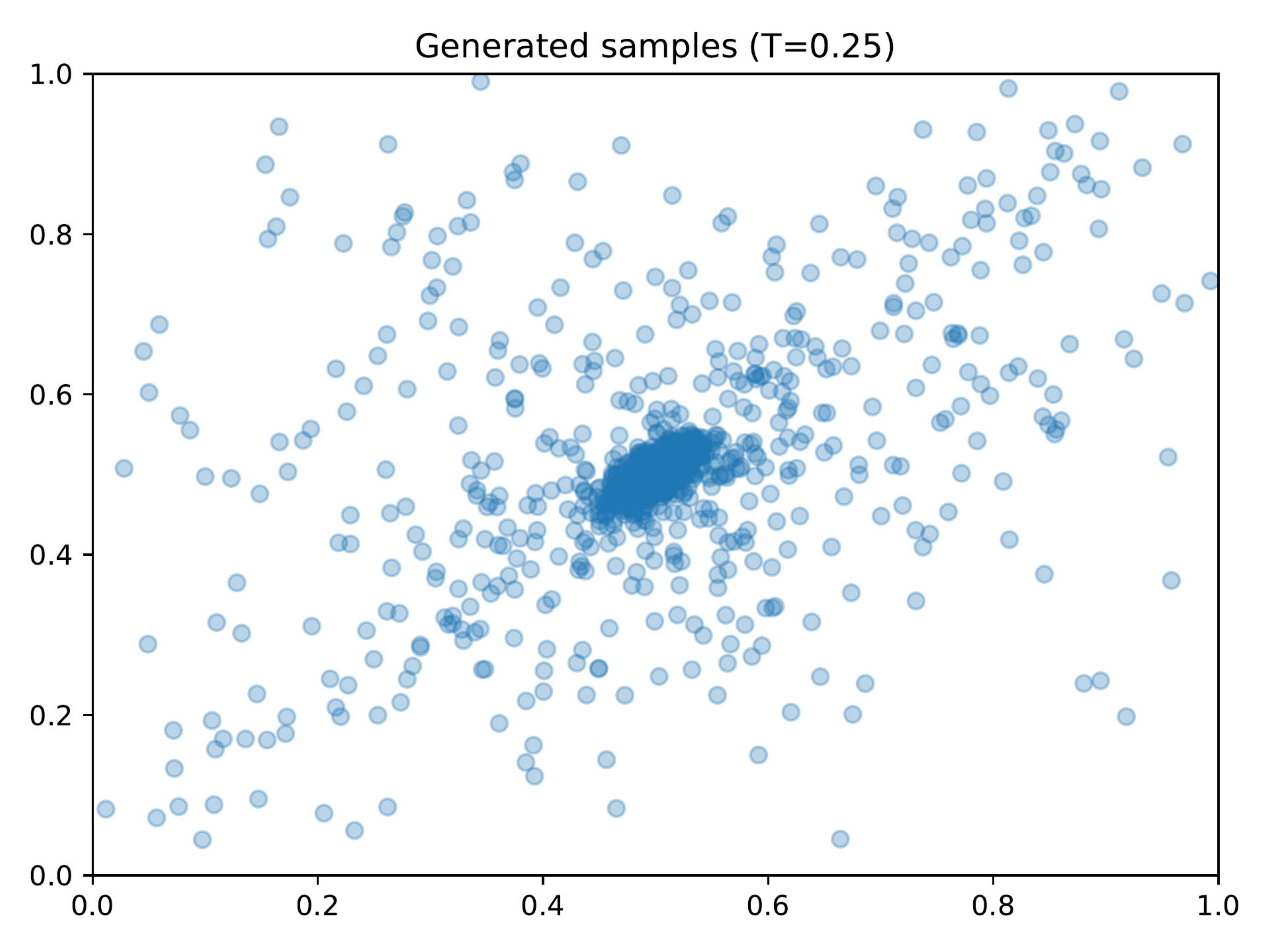}
	\includegraphics[width=0.3\linewidth]{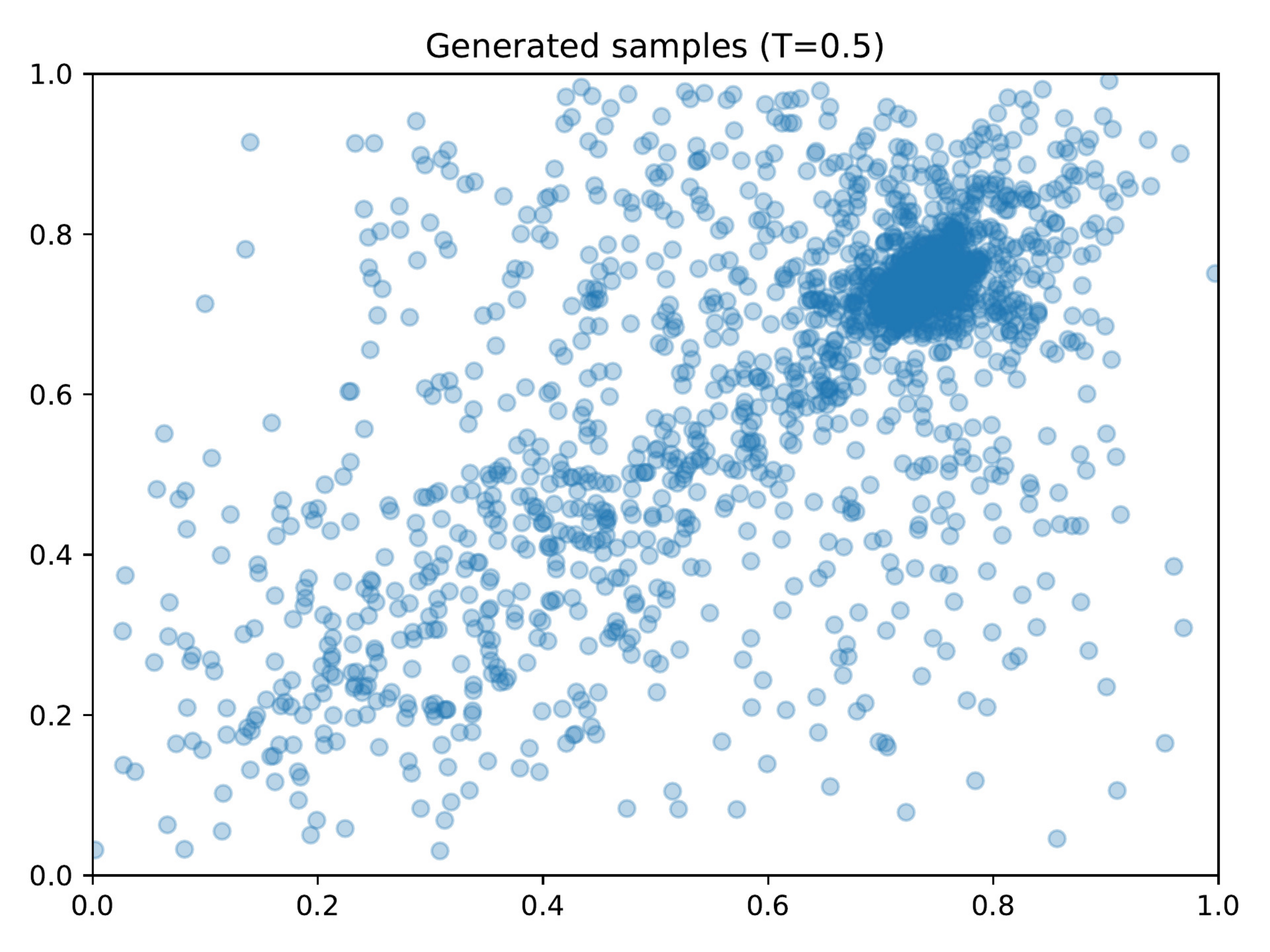}

	\caption{Low regularity test problem. From top to bottom: The generated samples of the 1st and 5th adaptivity iteration. ($\beta=1000,$ Galerkin projection)}
	\label{2D_heat_equation_samples}
\end{figure}

We also investigate the case when $\beta=1000$. Here we take the Galerkin projection as an example.
\begin{figure}[H]
	\centering
	\includegraphics[width=0.4\linewidth]{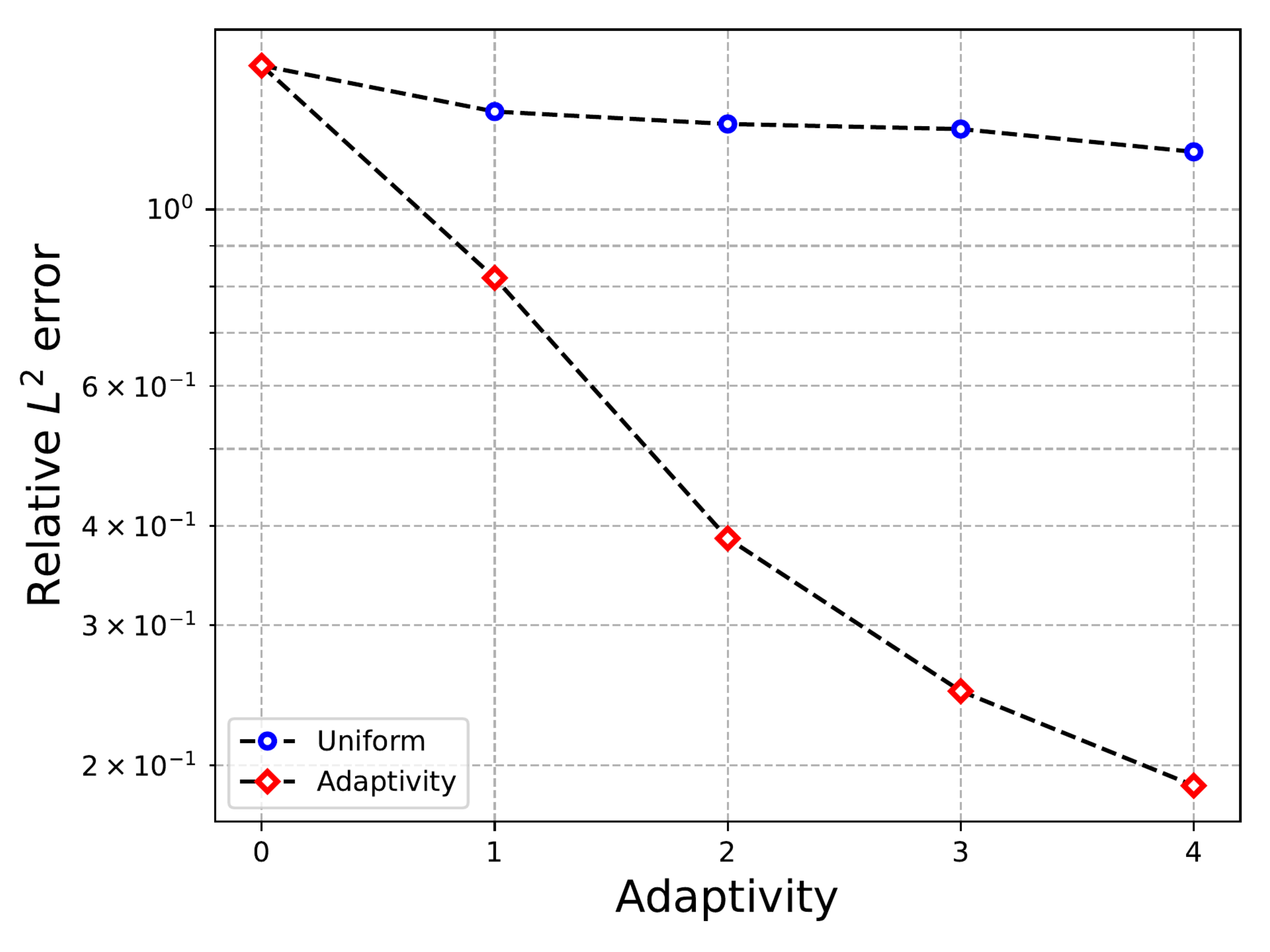}
	\caption{Low regularity test problem. Relative $L_2$ error for different sampling strategies. ($\beta$=1000, Galerkin projection).}
	\label{fig:heatadaptivitytwodim1000}
\end{figure}
In Figure \ref{fig:heatadaptivitytwodim1000}, we compare the performance of uniform sampling and adaptive sampling using the relative $L_2$ error, which demonstrates the necessity of adaptive sampling for this case.
Figure \ref{fig:galerkin_temporal_weights_beta_1000} shows the evolution of temporal weights $\lambda_i$ for $k=1,3,5$. The observed temporal weights exhibit an alternating high-low pattern since the strong form residuals tend to be uniform, taking the inner product with two different types of the test piecewise quadratic basis functions will result in different magnitudes. We note that both high-amplitude $\lambda_i$ and low-amplitude $\lambda_i$ tend to become uniform. One also finds that the error profiles become more flat as the adaptivity iterations increase.
Figure \ref{2D_heat_equation_samples} shows the
generated samples of bounded KRnet with respect to adaptivity iterations $k=1, 5$.
It is seen that the largest density of generated samples at the first update is
around the peak of the reference solution, which is consistent with the residual-induced
distribution. Moreover, as $k$ increases, we expect the tail of the residual-induced
distribution becomes heavier since the adaptivity tries to make the residual-induced
distribution more uniform. It is illustrated by the last update in the Figure
\ref{2D_heat_equation_samples}, which is consistent with previous findings
reported in \cite{tang2023pinns}.

\section{Conclusion}\label{section:conclusion}
In this paper we have developed a hybrid numerical method for solving evolution partial differential equations
(PDEs) by merging the time finite element method with deep neural networks. The key idea of the proposed method is to represent the
solution as a tensor product comprising a series of
predetermined local finite element basis functions in
time and a sequence of unspecified neural networks in
space. Subsequently, we apply the Galerkin or collocation
formulation to the original evolution equation, eliminating
the temporal variable and resulting in a differential system
exclusively involving unknown neural networks with respect to
the spatial variable. Furthermore, to address the evolution
problems characterized by high dimensionality and low regularity, we have developed an adaptive
sampling strategy where the training set and the model are updated alternately to improve both efficiency and accuracy.
Numerical experiments have demonstrated the effectiveness of our proposed approaches.
Several issues deserve further investigation. First, a rigorous error analysis is still missing especially for the collocation projection. Second, while herein we only consider linear and quadratic finite element methods over a uniform mesh, the strategy can also be generalized to $hp$-finite element method in the temporal domain. Finally, the proposed methods need to be further investigated for convection-dominated problems.
\section*{Acknowledgements}
  X. Feng was partially supported by Guangdong and Hong Kong Universities ``1+1+1'' Joint Research Collaboration Scheme project (No. 2025A0505000014).
T. Tang was partially
supported by Science Challenge Project (Grant No. TZ2018001) and NSFC 11731006 and
K20911001
\appendix
\section{Proof of Theorem \ref{theorem:error_estimate_heat_equation}}\label{app:proof_of_error_estimate}
\begin{proof}
	We use the notation $A \lesssim B$ to denote $A \leq CB$ for some constant $C>0$ only depending on $\Omega, T, \sigma, u, \nu$. As discussed in the error decomposition \eqref{eqn:error_decomposition}, we will estimate the first three terms separately.

	\paragraph{Finite element approximation error:} Partition $[0,T]$ with equidistant points $0=t_0<t_1<\cdots<t_N=T$, denote $\tau = \tau_n \coloneqq t_n-t_{n-1} = T/N, n = 1,\ldots,N$. One can define the following time finite element space
	\begin{equation*}
		V_{\tau} \coloneqq \{v\in C([0,T];\mathbb{R})\mid v|_{[t_{n-1},t_n]} \mbox{\; is linear for each } n \}.
	\end{equation*}
	Note that test functions belong to $V_\tau$ and depends only on time. For all $v\in V_{\tau}$, the semi-discrete variational formulation reads
	\begin{equation*}
		\int_0^T \partial_t u(\bm{x},t)v(t) \mathrm{d}t - \int_0^T \Delta u(\bm{x},t) v(t)\mathrm{d}t = \int_0^Tf(\bm{x},t)v(t)\mathrm{d}t,
	\end{equation*}
	with the initial and boundary conditions. The finite element method aims to seek
	\begin{equation*}
		u_{N}(\bm{x}, t) = \sum_{j=0}^N u_j(\bm{x})\phi_j(t)
	\end{equation*}
	satisfying the PDE, where $\{\phi_j\}_{j=0}^{N}$ is the standard basis of $V_{\tau}$ and $u_j(\bm{x})\in H_0^1(\Omega)\cap H^2(\Omega)$. Hence for each basis function $\phi_i(t)$, we have
	\begin{equation*}
		\int_0^T \partial_t u_{N}(\bm{x},t)\phi_i(t) \mathrm{d}t - \int_0^T \Delta u_{N}(\bm{x},t)\phi_i(t)\mathrm{d}t = \int_0^Tf(\bm{x},t)\phi_i(t)\mathrm{d}t, \quad i=0,\cdots,N.
	\end{equation*}
	We mention that the above system has unique solution $u_{N}(\bm{x},t)$. Our goal is to estimate the error between $u(\bm{x}, t)$ and its corresponding finite element approximation solution $u_{N}(\bm{x},t)$. To this end, we define the piecewise linear interpolation operator $\prod_\tau:C([0,T],L^2(\Omega))\to V_\tau \otimes L^2(\Omega)$ by
	\begin{equation*}
		(\textstyle \prod_\tau u)(\bm{x},t_n) = u(\bm{x},t_n), \quad n=0,\cdots,N,
	\end{equation*}
	and linear in time on each subinterval $[t_{n-1},t_n]$.  Then we have the following standard interpolation error estimates
	\begin{equation*}
		\begin{aligned}
			 & \| u - \textstyle \prod _\tau u\|_{L^2(0,T;L^2(\Omega))} \lesssim \tau^2 \|\partial_{tt}u\|_{L^2(0,T;L^2(\Omega))} \lesssim \tau^2,                   \\
			 & \| \partial_t u - \partial_t \textstyle \prod _\tau u\|_{L^2(0,T;L^2(\Omega))} \lesssim \tau \|\partial_{tt}u\|_{L^2(0,T;L^2(\Omega))} \lesssim \tau. \\
		\end{aligned}
	\end{equation*}
	We decompose the finite element approximation error into two parts as follows:
	\begin{equation*}
		e \coloneqq u -u_{N} = (u - \textstyle \prod_\tau u) + (\textstyle \prod_\tau u - u_{N}) \coloneqq e_1 + e_2.
	\end{equation*}
	Since $u_{N}$ satisfies the discrete variational formulation and $u$ satisfies the continuous one, for all $v\in V_{\tau}$, we have
	\begin{equation*}
		\int_0^T \partial_t (u-u_{N})v \mathrm{d}t -\int_0^T \Delta (u-u_{N})v\mathrm{d}t = 0.
	\end{equation*}
	Fixing $\bm{x}\in\Omega$ and plugging $v=e_2=\prod_\tau u - u_{N}\in V_\tau$, this yields
	\begin{equation*}
		\int_0^T \partial_t e_2 \cdot e_2 \mathrm{d}t -\int_0^T \Delta e_2 \cdot e_2\mathrm{d}t =-\int_0^T \partial_t e_1 \cdot e_2\mathrm{d}t + \int_0^T \Delta e_1\cdot e_2\mathrm{d}t.
	\end{equation*}
	Observe that
	\begin{equation*}
		\int_0^T \partial_t e_2 \cdot e_2 \mathrm{d}t =  \frac{1}{2} e_2^2(\bm{x},T) -\frac{1}{2}e_2^2(\bm{x},0),
	\end{equation*}
	and by integration by parts in space
	\begin{equation*}
		-\int_0^T\int_{\Omega} \Delta e_2 \cdot e_2\mathrm{d}\bm{x}\mathrm{d}t = \int_0^T \|\nabla e_2\|^2_{L^2(\Omega)}\mathrm{d}t.
	\end{equation*}
	Since the initial condition is embedded in the solution space, we have $e_2(\bm{x},0) = 0$. Hence, by integrating over $\Omega$, we have
	\begin{equation*}
		\frac{1}{2}\|e_2(\bm{x},T)\|_{L^2(\Omega)}^2 + \int_0^T \|\nabla e_2\|_{L^2(\Omega)}^2\mathrm{d}t = -\int_0^T \int_\Omega \partial_t e_1\cdot e_2\mathrm{d}\bm{x}\mathrm{d}t + \int_0^T\int_\Omega \Delta e_1\cdot e_2 \mathrm{d}\bm{x}\mathrm{d}t.
	\end{equation*}
	Using Cauchy-Schwarz and Young's inequalities with arbitrary $\epsilon>0$:
	\begin{equation*}
		\left|\int_0^T\int_\Omega \partial_t e_1\cdot e_2\mathrm{d}\bm{x}\mathrm{d}t\right| \leq  \epsilon \|e_2\|^2_{L^2(0,T;L^2(\Omega))} + \frac{1}{4\epsilon} \|\partial_t e_1\|^2_{L^2(0,T;L^2(\Omega))}.
	\end{equation*}
	Similarly,
	\begin{equation*}
		\left| \int_0^T\int_\Omega \Delta e_1\cdot e_2\mathrm{d}\bm{x}\mathrm{d}t \right| \leq \epsilon \|e_2\|^2_{L^2(0,T;L^2(\Omega))} + \frac{1}{4\epsilon} \|\Delta e_1\|^2_{L^2(0,T;L^2(\Omega))}.
	\end{equation*}
	By Poincar\'{e} inequality, $\|e_2\|_{L^2(\Omega)}\lesssim \|\nabla e_2\|_{L^2(\Omega)}$, so for sufficiently small $\epsilon$ we absorb $\|e_2\|_{L^2(0,T;L^2(\Omega))}$ terms into the left side:
	\begin{equation*}
		\|e_2(\bm{x},T)\|_{L^2(\Omega)}^2 + \|\nabla e_2\|_{L^2(0,T;L^2(\Omega))}^2 \lesssim \left( \|\partial_t e_1\|_{L^2(0,T;L^2(\Omega))}^2 + \|\Delta e_1\|_{L^2(0,T;L^2(\Omega))}^2\right)
	\end{equation*}
	Since $e_1 = u- \prod_\tau u$ interpolation error estimates yields
	\begin{equation*}
		\begin{aligned}
			 & \|e_1\|_{L^2(0,T;L^2(\Omega))} \lesssim \tau^2 \|\partial_{tt} u\|_{L^2(0,T;L^2(\Omega))} \lesssim \tau^2,               \\
			 & \|\partial_t e_1\|_{L^2(0,T;L^2(\Omega))} \lesssim \tau \|\partial_{tt} u\|_{L^2(0,T;L^2(\Omega))} \lesssim \tau,        \\
			 & \| \Delta e_1\|_{L^2(0,T;L^2(\Omega))} \lesssim \tau^2 \|\partial_{t} \Delta u\|_{L^2(0,T;L^2(\Omega))} \lesssim \tau^2.
		\end{aligned}
	\end{equation*}
	Combining the above estimates, we have
	\begin{equation*}
		\|e_2\|_{L^2(0,T;L^2(\Omega))}^2 \le \|e_2(\bm{x},T)\|_{L^2(\Omega)}^2 + \|\nabla e_2\|_{L^2(0,T;L^2(\Omega))}^2 \lesssim \tau^2,
	\end{equation*}
	and
	\begin{equation*}
		\|e\|_{L^2(0,T;L^2(\Omega))} \leq \|e_1\|_{L^2(0,T;L^2(\Omega))} + \|e_2\|_{L^2(0,T;L^2(\Omega))} \lesssim \tau.
	\end{equation*}
	Finally, we have
	\begin{equation}
		\|u(\bm{x}, t) - u_{N}(\bm{x}, t)\|_{L^2(0,T;L^2(\Omega))} \lesssim \tau \lesssim N^{-1}.
		\label{eqn:finite_element_approximation_error}
	\end{equation}

	\paragraph{Neural network approximation error:} Note that $u_N(\bm{x}, t)$ denotes the temporal finite element approximation of the solution and it satisfies
	\begin{equation*}
		\left\{
		\begin{aligned}
			 & \left(\partial_t u(\bm{x}, t), \phi_i(t)\right) - \left(\Delta u(\bm{x}, t), \phi_i(t)\right) = \left(f(\bm{x}, t), \phi_i(t)\right), & \quad \bm{x} \in \Omega, \quad i = 0,\ldots,N, \\
			 & u(\bm{x}, t) = 0,                                                                                                                     & (\bm{x}, t)\in \partial \Omega \times (0,T],   \\
			 & u(\bm{x},0) = u_0(\bm{x}),                                                                                                            & \bm{x}\in \Omega.
		\end{aligned}
		\right.
	\end{equation*}
	Since we embed the initial and boundary conditions into the hypothesis space in the form of \eqref{eqn:solution_embedding}, the continuous loss function can be simplified as
	\begin{equation*}
		\begin{aligned}
			\mathcal{L}(\theta) & = \sum_{i=0}^N\Vert \left((\partial_t - \Delta)u_{N, M}(\bm{x}, t;\theta) - f(\bm{x}, t), \phi_i(t)\right) \Vert_{L^2(\Omega)}^2      \\
			                    & = \sum_{i=0}^N\Vert \left((\partial_t - \Delta)(u_{N, M}(\bm{x}, t;\theta) - u_N(\bm{x}, t)), \phi_i(t)\right) \Vert_{L^2(\Omega)}^2.
		\end{aligned}
	\end{equation*}
	By minimizing the above continuous loss function, we can get the best approximation \( u_{N, M}(\bm{x}, t;\theta^*) \) of \( u_N(\bm{x}, t) \) and
	\begin{equation*}
		\theta^* = \mathop{\arg\min}\limits_{\theta}\mathcal{L}(\theta).
	\end{equation*}

	First of all, we need to prove that neural network approximation error is bounded by the continuous loss function, that is,
	\begin{equation*}
		\Vert u_{N, M}(\bm{x}, t;\theta) - u_N(\bm{x}, t) \Vert^2_{L^2(0, T; L^2(\Omega))} \lesssim \mathcal{L}(\theta).
	\end{equation*}
	Define the error
	\begin{equation*}
		e(\bm{x}, t; \theta) \coloneqq u_{N, M}(\bm{x}, t;\theta) - u_N(\bm{x}, t) = \sum_{j=1}^N \left(\nu(\bm{x})\omega_j(\bm{x}; \theta)-u_j(\bm{x})\right)\phi_j(t) = \sum_{j=0}^N e_j(\bm{x}, \theta)\phi_j(t),
	\end{equation*}
	where $e_0(\bm{x}, \theta) = 0, e_j(\bm{x}, \theta) = \nu(\bm{x})\omega_j(\bm{x}; \theta)-u_j(\bm{x}), j = 1, \ldots, N$ and the continuous residual
	\begin{equation*}
		r_i(\bm{x}; \theta) \coloneqq \left((\partial_t - \Delta)(u_{N, M}(\bm{x},t;\theta) - u_N(\bm{x},t)), \phi_i(t)\right) = \left((\partial_t - \Delta)e(\bm{x},t), \phi_i(t)\right), \quad i = 0, \ldots, N.
	\end{equation*}
	By inserting the expansion $e(\bm{x}, t; \theta) = \sum_{j=0}^N e_j(\bm{x}; \theta)\phi_j(t)$ and using linearity, we can obtain
	\begin{equation*}
		r_i(\bm{x}; \theta) = \sum_{j=0}^Na_{ji}e_j(\bm{x}; \theta) - \sum_{j=0}^Nb_{ji}\Delta e_j(\bm{x}; \theta), \quad i = 0, \ldots, N,
	\end{equation*}
	where $a_{ji} = (\phi'_j,\phi_i), b_{ji} = (\phi_j,\phi_i), j = 0, \ldots, N$.
	This is a linear system
	\begin{equation*}
		\mathcal{A}\bm{e} = \bm{r},
	\end{equation*}
	where $\bm{e} \coloneqq (e_0, \ldots, e_N)^T$ and $\bm{r} \coloneqq (r_0, \ldots, r_N)^T$. Then we define the Hilbert space
	\begin{equation*}
		H \coloneqq (L^2(\Omega))^{N+1},
	\end{equation*}
	and the corresponding inner product and norm in $H$ by
	\begin{equation*}
		\left<\bm{u}, \bm{v}\right>_H = \sum_{j=0}^N(u_i, v_i)_{L^2(\Omega)}, \quad \forall\bm{u}, \bm{v} \in H.
	\end{equation*}
	and the corresponding norm is
	\begin{equation*}
		\Vert \bm{u} \Vert_{H}^2 = \sum_{j=0}^N\Vert u_j\Vert_{L^2(\Omega)}^2.
	\end{equation*}
	Consider the bilinear form:
	\begin{equation}
		\begin{aligned}
			\left<\mathcal{A}\bm{e}, \bm{e}\right>_H & = \sum_{j, i}a_{ji}(e_j, e_i)_{L^2(\Omega)} - \sum_{j, i}b_{ji}(\Delta e_j, e_i)_{L^2(\Omega)}         \\
			                                         & = \sum_{j, i}a_{ji}(e_j, e_i)_{L^2(\Omega)} + \sum_{j, i}b_{ji}(\nabla e_j, \nabla e_i)_{L^2(\Omega)},
		\end{aligned}
		\label{eqn:bilinear_form}
	\end{equation}
	where the last equality uses integration by parts. Since the matrix $(a_{ji})_{(N+1)\times(N+1)}$ has the following form:
	\begin{equation*}
		a_{ji} = \left\{
		\begin{array}{ll}
			-\frac{1}{2}, & j=i=0 \ \mathrm{or} \ j-i=1, \\
			\frac{1}{2},  & j=i=N \ \mathrm{or} \ i-j=1, \\
			0,            & \mathrm{otherwise},
		\end{array}
		\right.
	\end{equation*}
	the first term of \eqref{eqn:bilinear_form} can be written as
	\begin{equation*}
		\sum_{j, i}a_{ji}(e_j, e_i)_{L^2(\Omega)} = \frac{1}{2}\Vert e_N \Vert_{L^2(\Omega)}^2 - \frac{1}{2}\Vert e_0 \Vert_{L^2(\Omega)}^2 = \frac{1}{2}\Vert e_N \Vert_{L^2(\Omega)}^2 \ge 0.
	\end{equation*}
	By positive definiteness of $(b_{ij})_{(N+1)\times(N+1)}$, the second term of \eqref{eqn:bilinear_form} satisfies
	\begin{equation*}
		\sum_{j, i}b_{ji}(\nabla e_j, \nabla e_i)_{L^2(\Omega)} \gtrsim \sum_{j=0}^N\Vert \nabla e_j\Vert_{L^2(\Omega)}^2 \gtrsim \sum_{j=0}^N\Vert e_j\Vert_{L^2(\Omega)}^2,
	\end{equation*}
	where the last inequality uses Poincar\'{e} inequality on each $e_j$. Therefore, the bilinear form is shown to be coercive:
	\begin{equation*}
		\left<\mathcal{A}\bm{e}, \bm{e}\right>_H \gtrsim \sum_{j=0}^N\Vert e_j\Vert_{L^2(\Omega)}^2 = \Vert \bm{e} \Vert_H^2.
	\end{equation*}
	From the linear system $\mathcal{A}\bm{e} = \bm{r}$, the coercivity and the Cauchy-Schwarz inequality, we have
	\begin{equation*}
		\Vert \bm{e} \Vert_H^2 \lesssim \left<\mathcal{A}\bm{e}, \bm{e}\right>_H \lesssim \Vert \bm{r} \Vert_H\Vert \bm{e} \Vert_H,
	\end{equation*}
	and a stability bound can be obtained:
	\begin{equation*}
		\Vert \bm{e} \Vert_H \lesssim \Vert \bm{r} \Vert_H.
	\end{equation*}
	Consider the full space-time $L^2$ norm of $e$
	\begin{equation*}
		\Vert e \Vert_{L^2(0, T; L^2(\Omega))}^2 = \int_0^T\Vert e(\bm{x}, t; \theta)\Vert_{L^2(\Omega)}^2 \mathrm{d}t = \sum_{i, j} a_{ij}(e_i, e_j)_{L^2(\Omega)} \lesssim \sum_{i=0}^N\Vert e_i\Vert_{L^2(\Omega)}^2 = \Vert \bm{e} \Vert_H^2,
	\end{equation*}
	where the inequality uses the matrix bound of $(a_{ij})_{(N+1)\times(N+1)}$. Therefore, we have
	\begin{equation*}
		\Vert e \Vert_{L^2(0, T; L^2(\Omega))}^2 \lesssim \Vert \bm{e} \Vert_H^2 \lesssim \Vert \bm{r} \Vert_H^2 = \mathcal{L}(\theta).
	\end{equation*}

	Now, we can prove that there exists an upper bound for the loss function corresponding to the best approximation \( u_{N, M}(\bm{x}, t;\theta^*) \). Since $\theta^*$ is the global minimizer, we have
	\begin{equation*}
		\mathcal{L}(\theta^*) \le \mathcal{L}(\theta), \quad \forall \theta,
	\end{equation*}
	and
	\begin{equation*}
		\begin{aligned}
			\mathcal{L}(\theta)
			 & = \sum_{i=0}^N \left\Vert \left((\partial_t - \Delta)e(\bm{x}, t;\theta), \phi_i(t)\right) \right\Vert_{L^2(\Omega)}^2               \\
			 & = \sum_{i=0}^N \left\Vert \left(\sum_{j=1}^N(a_{ji}e_j - b_{ji}\Delta e_j)\right) \right\Vert_{L^2(\Omega)}^2                        \\
			 & = \sum_{i=0}^N \int_{\Omega} \left(\sum_{j=1}^N\left(a_{ji}e_j - b_{ji}\Delta e_j\right)\right)^2 \mathrm{d}\bm{x}                   \\
			 & \le \sum_{i=0}^N 2N \int_{\Omega} \sum_{j=0}^N\left(a_{ji}^2e_j^2 + b_{ji}^2(\Delta e_j)^2 \right)\mathrm{d}\bm{x}                   \\
			 & = \sum_{i=0}^N 2N \sum_{j=0}^N \left(a_{ji}^2\Vert e_j \Vert_{L^2(\Omega)}^2 + b_{ji}^2\Vert \Delta e_j \Vert_{L^2(\Omega)}^2\right) \\
			 & \lesssim N^2 d \max_j\Vert e_j \Vert_{H^2(\Omega)}^2, \quad \forall \theta,
		\end{aligned}
	\end{equation*}
	where the first inequality uses the n-variable mean inequality and the second inequality uses the sparse property of $(a_{ij})_{(N+1)\times(N+1)}$ and $(b_{ij})_{(N+1)\times(N+1)}$.

	Since $u/\nu < \infty$ on $\partial\Omega$, it is easy to verify that $\frac{u_j}{\nu} < \infty$ on $\partial\Omega$ for each $j=1, \ldots, N$ and $\frac{u_j}{\nu}\in \mathcal{B}^3$ when $s=2$ by Lemma \ref{lemma:Barron_space}. Next by Lemma \ref{lemma:NN_approximation}, there exists a $\theta^{c}$ such that
	\begin{equation*}
		\Vert e_j \Vert_{H^2(\Omega)} = \Vert \nu \omega_j-u_j \Vert_{H^2(\Omega)} \lesssim (\Vert \nu\Vert_{\infty}+\Vert D\nu \Vert_{\infty}+\Vert D^2\nu\Vert_{\infty})\Vert \omega_j-\frac{u_j}{\nu} \Vert_{H^2(\Omega)} \lesssim \Vert \frac{u_j}{\nu}\Vert_{\mathcal{B}^{3}}M^{-\frac{1}{2}},
	\end{equation*}
	where the last inequality uses the approximation property of $\mathcal{F}_{\sigma}(M)$.

	Finally, we have
	\begin{equation}
		\Vert u_{N, M}(\bm{x}, t;\theta^*) - u_N(\bm{x}, t) \Vert_{L^2(0, T; L^2(\Omega))} \lesssim \sqrt{\mathcal{L}(\theta^*)} \le \sqrt{L(\theta^{c})} \lesssim N d^{\frac{1}{2}} M^{-\frac{1}{2}}.
		\label{eqn:NN_approximation_error_bound}
	\end{equation}

	\paragraph{Statistical error:} In practice, the continuous loss function is usually discretized numerically before the optimization with respect to $\theta$ is addressed and the empirical loss function is
	\begin{equation*}
		\widehat{\mathcal{L}}(\theta) = \frac{\vert \Omega\vert}{P}\sum_{p=1}^{P}\sum_{i=0}^N r_i^2(\bm{x}_p; \theta),
	\end{equation*}
	where $\{\bm{x}_p\}_{p=1}^P$ are random samples in $\Omega$ and
	\begin{equation*}
		r_i(\bm{x}; \theta) = \sum_{j=0}^Na_{ji}\nu(\bm{x})\omega_j(\bm{x}; \theta) - \sum_{j=0}^Nb_{ji}\Delta \left(\nu(\bm{x})\omega_j(\bm{x}; \theta)\right) + \left(f(\bm{x}, t), \phi_i(t)\right), \quad i = 0, \ldots, N.
	\end{equation*}
	For simplicity, we can let $\vert \Omega\vert = 1$. Then \( u_{N,M}(\bm{x}, t;\theta_P^*) \) denotes the best approximation obtained by minimizing the empirical loss over \( P \) training samples and
	\begin{equation*}
		\theta_P^* = \mathop{\arg\min}\limits_{\theta}\widehat{\mathcal{L}}(\theta).
	\end{equation*}

	Consider the hypothesis space of the residual functions
	\begin{equation*}
		\mathcal{H} \coloneqq \{x\mapsto\sum_{i=0}^Nr_i^2(\bm{x}; \theta):\theta\in\Theta\}
	\end{equation*}
	and the empirical Rademacher complexity of $\widehat{\mathcal{R}}_S(\mathcal{H})$ with respect to the training samples $S = \{\bm{x}_p\}_{p=1}^P$ is
	\begin{equation}
		\widehat{\mathcal{R}}_S(\mathcal{H}) \coloneqq \mathbb{E}_{\bm{\zeta}}\left[\sup_{h\in\mathcal{H}}\frac{1}{p}\sum_{p=1}^P \zeta_p \sum_{i=0}^Nr_i^2(\bm{x}_p; \theta)\right],
		\label{eqn:Rademacher_complexity_H}
	\end{equation}
	where $\bm{\zeta} = (\zeta_1, \ldots, \zeta_P)$ are i.i.d. Rademacher random variables. For each $r_i(\bm{x}; \theta)$, we consider two hypothesis spaces
	\begin{equation*}
		\begin{aligned}
			 & \mathcal{H}_i \coloneqq \{x\mapsto r_i(\bm{x}; \theta):\theta\in\Theta\},                                                                                                                                                            \\
			 & \widehat{\mathcal{H}}_i \coloneqq \{x\mapsto \hat{r}_i(\bm{x}; \theta) \coloneqq \sum_{j=0}^Na_{ji}\nu(\bm{x})\omega_j(\bm{x}; \theta) - \sum_{j=0}^Nb_{ji}\Delta \left(\nu(\bm{x})\omega_j(\bm{x}; \theta)\right):\theta\in\Theta\}
		\end{aligned}
	\end{equation*}
	and the corresponding empirical Rademacher complexities $\widehat{\mathcal{R}}_S(\mathcal{H}_i)$ and $\widehat{\mathcal{R}}_S(\widehat{\mathcal{H}}_i)$. From translation invariance property of Rademacher complexity, we have
	\begin{equation}
		\widehat{\mathcal{R}}_S(\mathcal{H}_i) = \widehat{\mathcal{R}}_S(\widehat{\mathcal{H}}_i), \quad i = 0, \ldots, N.
		\label{eqn:Rademacher_complexity_H_i_translation_invariance}
	\end{equation}
	Expand the spatial Laplacian:
	\begin{equation*}
		\Delta(\nu \omega_j) = (\Delta \nu)\omega_j + 2 \nabla v\cdot\nabla \omega_j + \nu\Delta \omega_j,
	\end{equation*}
	and we can write
	\begin{equation*}
		\hat{r}_i(\bm{x}; \theta) = \sum_{j=0}^N\left(\alpha_{ji}\omega_j(\bm{x}; \theta) + \sum_{k=1}^d\beta_{jik} \frac{\partial \omega_j(\bm{x}; \theta)}{\partial x_k} + \sum_{k=1}^d\gamma_{jik} \frac{\partial^2 \omega_j(\bm{x}; \theta)}{\partial x_k^2} \right),
	\end{equation*}
	where $\alpha_{ji}, \beta_{jik}, \gamma_{jik}$ depend on $a_{ji}, b_{ji}$ and the fixed function $\nu$. We can always find a constant $M$ such that $\vert \alpha_{ji} \vert, \vert \beta_{jik} \vert, \vert \gamma_{jik} \vert \le M$. By the subadditivity property of Rademacher complexity,
	\begin{equation}
		\widehat{\mathcal{R}}_S(\widehat{\mathcal{H}}_i) \le \sum_{j=0}^N\left(\vert\alpha_{ji}\vert \widehat{\mathcal{R}}_S(\mathcal{W}) + \sum_{k=1}^d\vert\beta_{jik}\vert \widehat{\mathcal{R}}_S(\mathcal{W}^{1}) + \sum_{k=1}^d\vert\gamma_{jik}\vert \widehat{\mathcal{R}}_S(\mathcal{W}^{2})\right),
		\label{eqn:Rademacher_complexity_hat_H_i_inequality}
	\end{equation}
	where
	\begin{equation*}
		\begin{aligned}
			 & \mathcal{W} \coloneqq \{x\mapsto \omega(\bm{x}; \theta) = \sum_{m=1}^Mc_m\sigma(\bm{w}_m^Tx+b_m)\},           \\
			 & \mathcal{W}^{1} \coloneqq \{x\mapsto \omega(\bm{x}; \theta^1) = \sum_{m=1}^Mc_m^1\sigma'(\bm{w}_m^Tx+b_m)\},  \\
			 & \mathcal{W}^{2} \coloneqq \{x\mapsto \omega(\bm{x}; \theta^2) = \sum_{m=1}^Mc_m^2\sigma''(\bm{w}_m^Tx+b_m)\},
		\end{aligned}
	\end{equation*}
	Here, we need to select a sufficiently large $Q$ such that $\theta^c$ introduced in the neural network approximation error estimation is included within the neural network parameter space and the path norms of $\mathcal{W}, \mathcal{W}^1, \mathcal{W}^2$ are all bounded by $Q$. Then by Lemma \ref{lemma:generalization_bound}, we have
	\begin{equation}
		\begin{aligned}
			 & \widehat{\mathcal{R}}_S(\mathcal{W}) \le 2\gamma(\sigma)Q\sqrt{\frac{2\ln(2d+2)}{P}},     \\
			 & \widehat{\mathcal{R}}_S(\mathcal{W}^1) \le 2\gamma(\sigma')Q\sqrt{\frac{2\ln(2d+2)}{P}},  \\
			 & \widehat{\mathcal{R}}_S(\mathcal{W}^2) \le 2\gamma(\sigma'')Q\sqrt{\frac{2\ln(2d+2)}{P}}.
		\end{aligned}
		\label{eqn:Rademacher_complexity_W_W1_W2_inequalities}
	\end{equation}
	Next, we need to bound the gap between continuous and empirical losses. Since the path norm $\Vert \theta \Vert_p \leq Q$ for all $\theta$ in the parameter space, there exists a constant $B > 0$ such that for all $\theta$ and $x$,
	\begin{equation*}
		\left(\sum_{i=0}^N r_i^2(\bm{x}; \theta)\right)^{\frac{1}{2}} \le B.
	\end{equation*}
	Then by Lemma \ref{lemma:Rademacher_complexity}, we can get the gap between continuous and empirical losses based on Rademacher complexity, with probability at least $1-\delta$,
	\begin{equation}
		\sup_{\theta}\vert \mathcal{L}(\theta)-\widehat{\mathcal{L}}(\theta)\vert \le 2\widehat{\mathcal{R}}_S(\mathcal{H}) + 3B\sqrt{\frac{\ln(2/\delta)}{2P}}.
		\label{eqn:generalization_bound_Rademacher_complexity}
	\end{equation}
	Note that the sum over $i$ in \eqref{eqn:Rademacher_complexity_H} is finite, we split:
	\begin{equation*}
		\frac{1}{P}\sum_{p=1}^P\zeta_p\sum_{i=0}^Nr_i^2(\bm{x}_p; \theta) = \sum_{i=0}^N\frac{1}{P}\sum_{p=1}^P\zeta_pr_i^2(\bm{x}_p; \theta)
	\end{equation*}
	and we have
	\begin{equation*}
		\mathbb{E}_{\bm{\zeta}}\left[\sup_{\theta}\frac{1}{P}\sum_{p=1}^P\zeta_p\sum_{i=0}^Nr_i^2(\bm{x}_p; \theta)\right] \le \sum_{i=0}^N\mathbb{E}_{\bm{\zeta}}\left[\sup_{\theta}\frac{1}{P}\sum_{p=1}^P\zeta_p r_i^2(\bm{x}_p; \theta)\right].
	\end{equation*}
	Since $\vert r_i\vert < B$, for fixed $i$, the map $r_i\mapsto r_i^2$ is $2B$-Lipschitz continuous on $[-B,B]$. Thus by Ledoux-Talagrand contraction,
	\begin{equation*}
		\mathbb{E}_{\bm{\zeta}}\left[\sup_{\theta}\frac{1}{P}\sum_{p=1}^P\zeta_p r_i^2(\bm{x}_p; \theta)\right] \le 2B \mathbb{E}_{\bm{\zeta}}\left[\sup_{\theta}\frac{1}{P}\sum_{p=1}^P\zeta_p r_i(\bm{x}_p; \theta)\right].
	\end{equation*}
	Therefore, we have
	\begin{equation}
		\widehat{\mathcal{R}}_S(\mathcal{H}) \le 2B \sum_{i=0}^N\mathbb{E}_{\bm{\zeta}}\left[\sup_{\theta}\frac{1}{P}\sum_{p=1}^P\zeta_p r_i(\bm{x}_p; \theta)\right] = 2B \sum_{i=0}^N\widehat{\mathcal{R}}_S(\mathcal{H}_i).
		\label{eqn:Rademacher_complexity_H_inequality}
	\end{equation}
	Combining \eqref{eqn:Rademacher_complexity_H_i_translation_invariance}, \eqref{eqn:Rademacher_complexity_hat_H_i_inequality}, \eqref{eqn:Rademacher_complexity_W_W1_W2_inequalities}, \eqref{eqn:generalization_bound_Rademacher_complexity} and \eqref{eqn:Rademacher_complexity_H_inequality}, with probability at least $1-\delta$, we have
	\begin{equation}
		\sup_{\theta}\vert \mathcal{L}(\theta)-\widehat{\mathcal{L}}(\theta)\vert \le 4B \sum_{i=0}^N\widehat{\mathcal{R}}_S(\widehat{\mathcal{H}}_i) + 3B\sqrt{\frac{\ln(2/\delta)}{2P}} \lesssim N^2d \sqrt{\frac{\ln d}{P}} + \sqrt{\frac{\ln(2/\delta)}{P}}.
		\label{eqn:generalization_bound_final}
	\end{equation}
	For the difference between the continuous loss function values taken on $\theta^*$ and $\theta_P^*$, we have
	\begin{equation*}
		\mathcal{L}(\theta^*_P) - \mathcal{L}(\theta^*) \le 2\sup_{\theta}\vert \mathcal{L}(\theta)-\widehat{\mathcal{L}}(\theta)\vert
	\end{equation*}
	since
	\begin{equation*}
		\mathcal{L}(\theta_P^*) \le \widehat{\mathcal{L}}(\theta_P^*) + \sup\vert \mathcal{L}-\widehat{\mathcal{L}}\vert \le \widehat{\mathcal{L}}(\theta^*) + \sup\vert \mathcal{L}-\widehat{\mathcal{L}}\vert \le \mathcal{L}(\theta^*) + 2\sup\vert \mathcal{L}-\widehat{\mathcal{L}}\vert.
	\end{equation*}
	Similar to the proof of neural network approximation error estimation, we have
	\begin{equation*}
		\begin{aligned}
			         & \Vert u_{N, M}(\bm{x}, t;\theta_P^*) - u_{N, M}(\bm{x}, t;\theta^*)\Vert_{L^2(0,T;L^2(\Omega))}^2                   \\
			\lesssim & \sum_{i=0}^N\Vert r_i(\bm{x}; \theta_P^*)-r_i(\bm{x}; \theta^*)\Vert_{L^2(\Omega)}^2                                \\
			\lesssim & \sum_{i=0}^N(\Vert r_i(\bm{x}; \theta_P^*)\Vert_{L^2(\Omega)}^2 + \Vert r_i(\bm{x}; \theta^*)\Vert_{L^2(\Omega)}^2) \\
			=        & \mathcal{L}(\theta_P^*)+\mathcal{L}(\theta^*)                                                                       \\
			=        & 2\mathcal{L}(\theta^*) + \left(\mathcal{L}(\theta_P^*)-\mathcal{L}(\theta^*)\right)                                 \\
			\lesssim & N^2 d M^{-1} + N^2d \sqrt{\frac{\ln d}{P}} + \sqrt{\frac{\ln(2/\delta)}{P}},
		\end{aligned}
	\end{equation*}
	where the last inequality uses \eqref{eqn:NN_approximation_error_bound} and \eqref{eqn:generalization_bound_final}.

	Finally, with probability at least $1-\delta$, we have
	\begin{equation}
		\Vert u_{N, M}(\bm{x}, t;\theta_P^*) - u_{N, M}(\bm{x}, t;\theta^*)\Vert_{L^2(0,T;L^2(\Omega))} \lesssim \sqrt{N^2 d M^{-1} + N^2d \sqrt{\frac{\ln d}{P}} + \sqrt{\frac{\ln(2/\delta)}{P}}}.
		\label{eqn:statistical_error_bound}
	\end{equation}
	\paragraph{Total error:} Combining \eqref{eqn:finite_element_approximation_error}, \eqref{eqn:NN_approximation_error_bound} and \eqref{eqn:statistical_error_bound} and using the triangle inequality, with probability at least $1-\delta$, we have
	\begin{equation*}
		\begin{aligned}
			         & \Vert u_{N, M}(\bm{x}, t;\theta_P^*) - u(\bm{x}, t)\Vert_{L^2(0,T;L^2(\Omega))}                                                                                                   \\
			\le      & \Vert u_{N, M}(\bm{x}, t;\theta_P^*) - u_{N, M}(\bm{x}, t;\theta^*)\Vert_{L^2(0,T;L^2(\Omega))} + \Vert u_{N, M}(\bm{x}, t;\theta^*) - u_N(\bm{x}, t)\Vert_{L^2(0,T;L^2(\Omega))} \\
			         & + \Vert u_N(\bm{x}, t) - u(\bm{x}, t)\Vert_{L^2(0,T;L^2(\Omega))}                                                                                                                 \\
			\lesssim & N^{-1} + N d^{\frac{1}{2}} M^{-\frac{1}{2}} + \sqrt{N^2 d M^{-1} + N^2d \sqrt{\frac{\ln d}{P}} + \sqrt{\frac{\ln(2/\delta)}{P}}}.
		\end{aligned}
	\end{equation*}

\end{proof}

\section{Exact periodic boundary conditions} \label{app_periodic_boundary_condition}
Following the work from \cite{wang2024respecting,penwarden2023unified}, we can exactly enforce $C^{\infty}$ periodic boundary conditions by applying a Fourier feature encoding to the spatial input of the network. The spatial encoding is
\begin{equation*}
	v(x) = \left\{ 1,\cos(\omega x),\sin(\omega x) ,\cdots,\cos(M\omega x),\sin(M\omega x)\right\}
\end{equation*}
where $\omega =\frac{2\pi}{L}, L = x_{\mathrm{max}} - x_{\mathrm{min}}$, and $M$ is a non-negative integer representing the sinusoidal frequency of the input. In this work, we take $M=10.$

\section{Supplementary Figures}\label{app:adaptivity_figure}
\begin{figure}[H]
	\centering
	\includegraphics[width=0.4\linewidth]{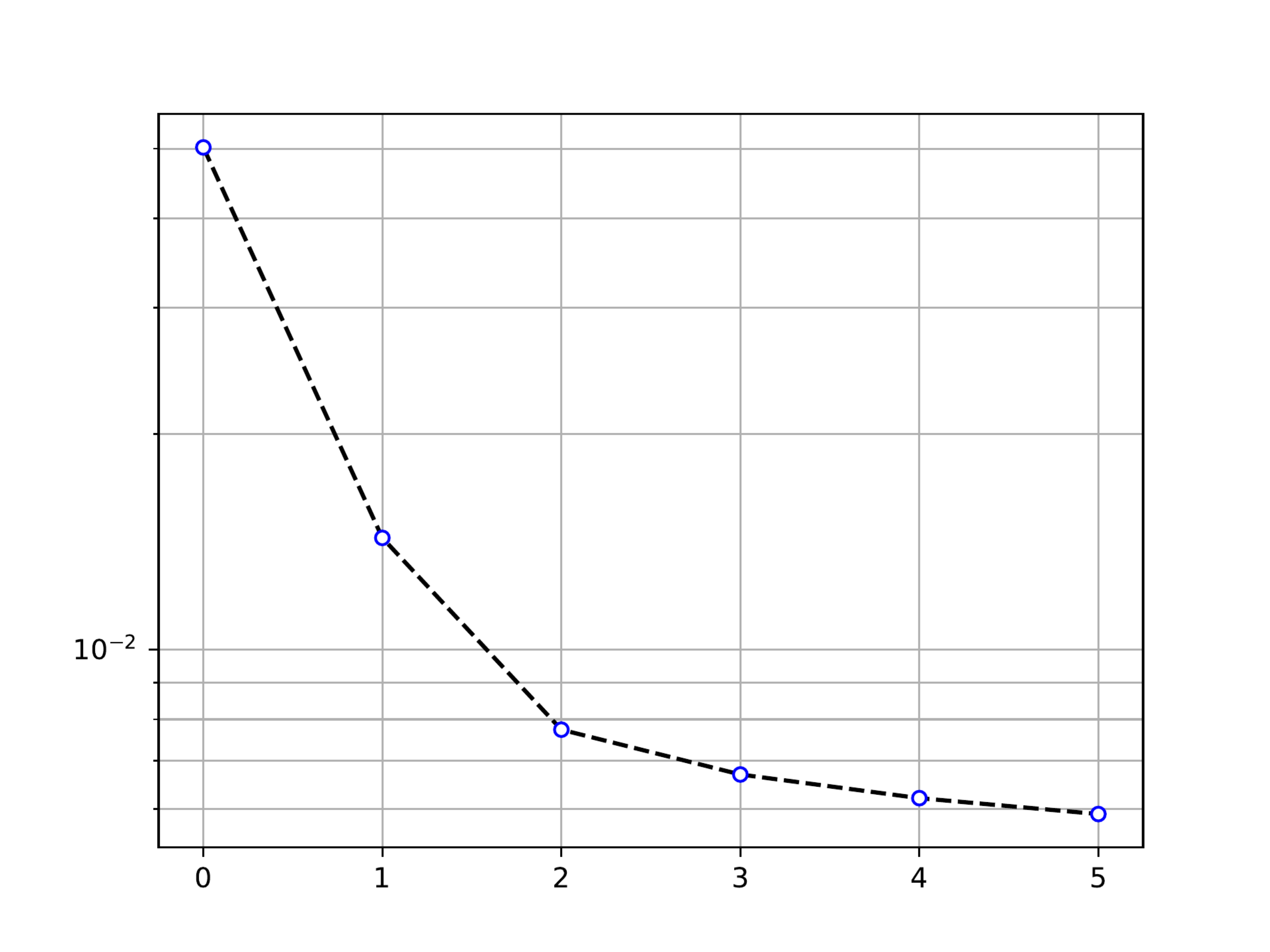}
	\caption{Low regularity test problem. Relative $L_2$ error of different adaptivity iterations. ($\beta=200$, Galerkin projection)}
	\label{fig:galerkin_error_beta_200}
\end{figure}

\begin{figure}[H]
	\centering
	\includegraphics[width=0.3\linewidth]{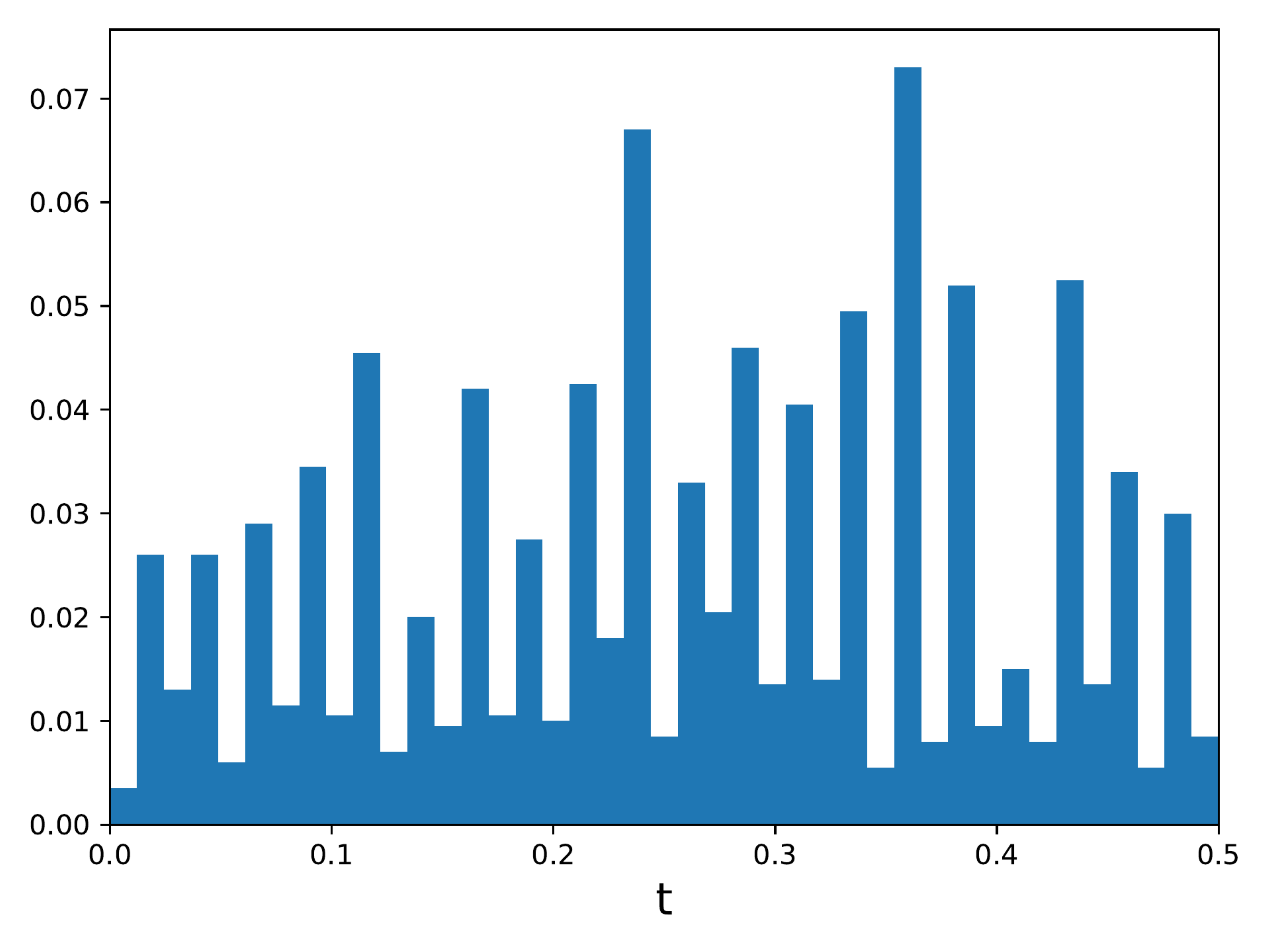}
	\includegraphics[width=0.3\linewidth]{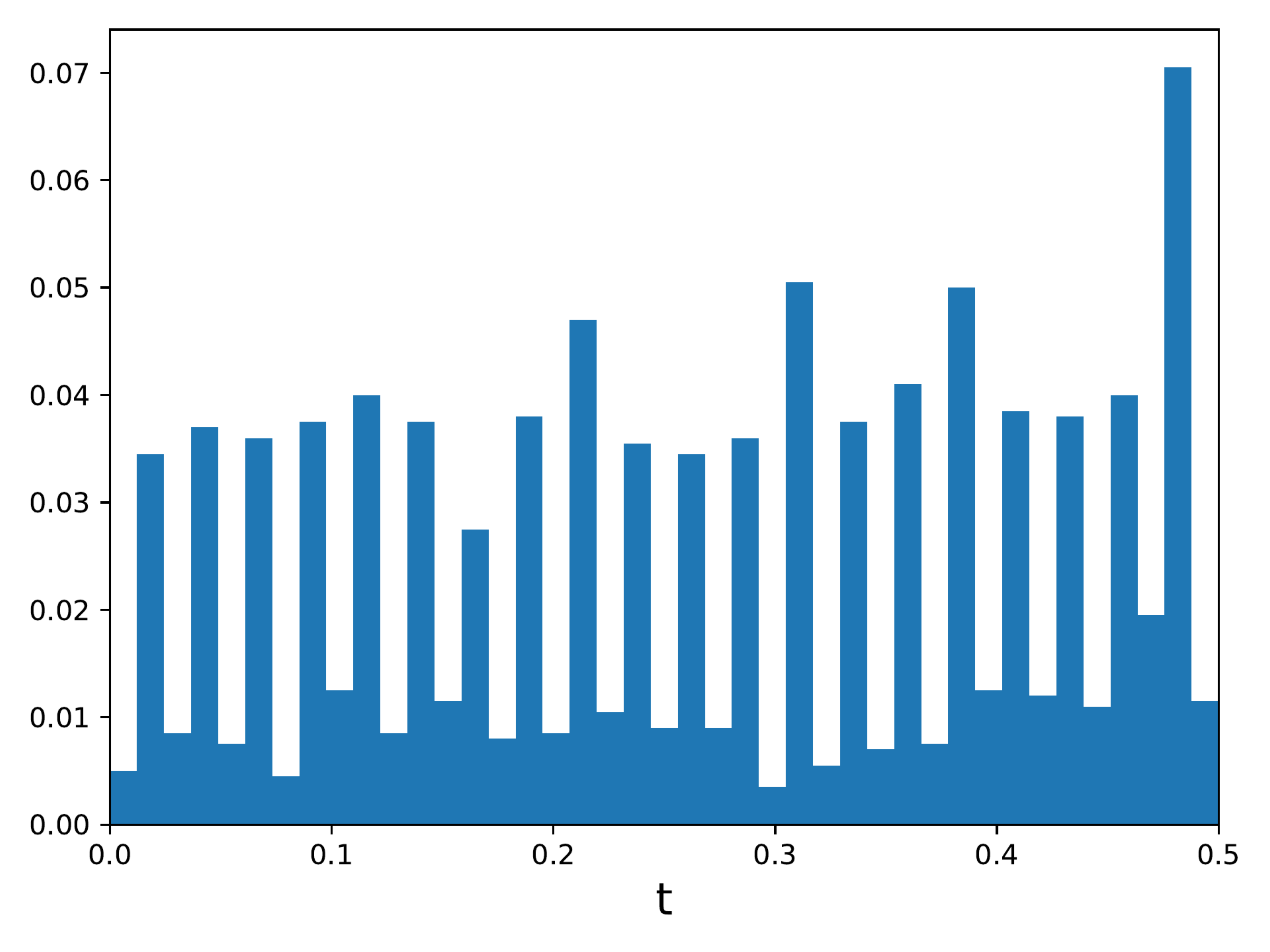}
	\includegraphics[width=0.3\linewidth]{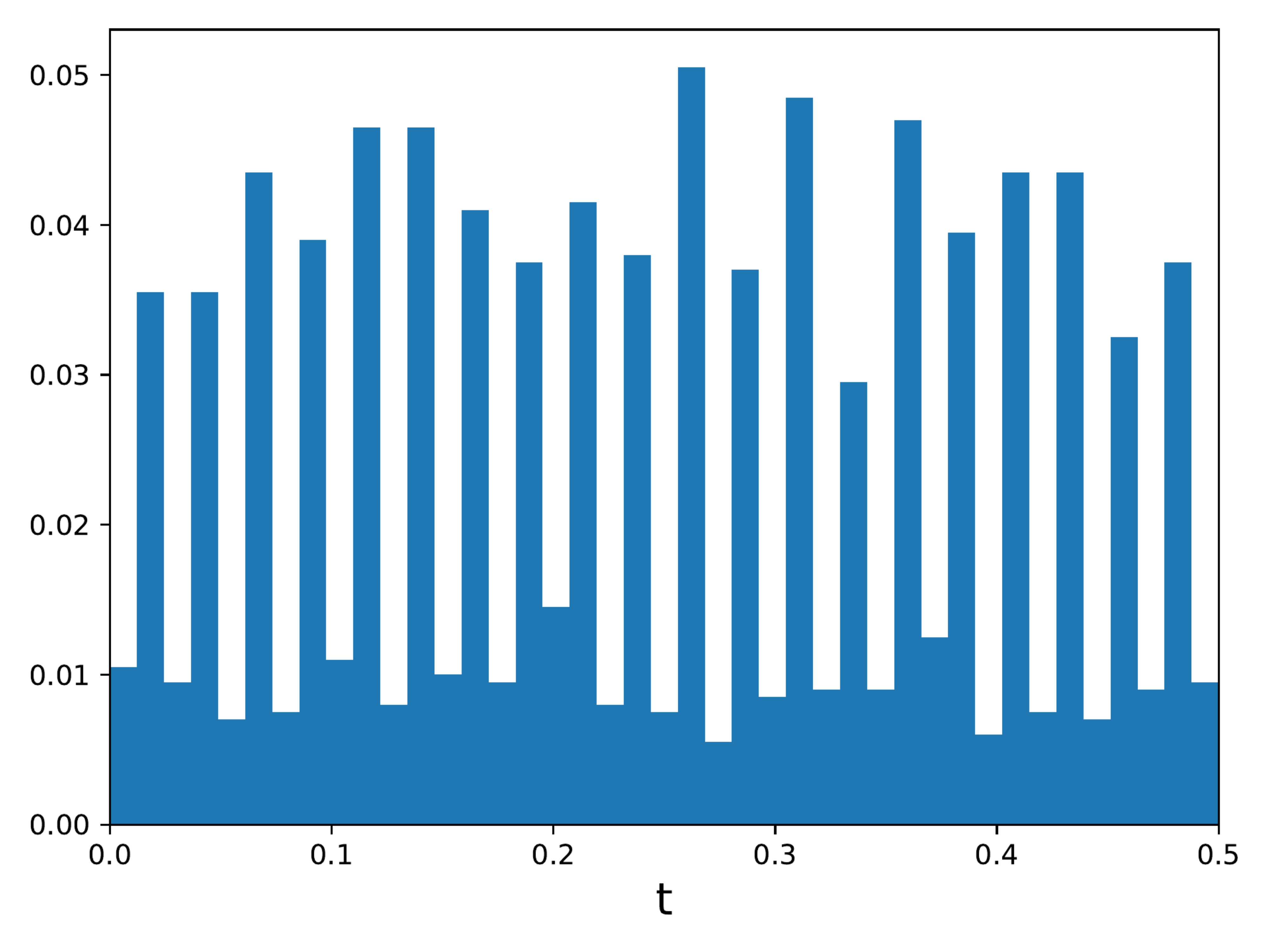}
	\caption{Low regularity test problem. Temporal weights $\lambda_i$ of the 1st, 3rd, 5th adaptivity iteration. ($\beta=200$, Galerkin projection)}
	\label{fig:galerkin_temporal_weights_beta_200}
\end{figure}

\begin{figure}[H]
	\centering
	\includegraphics[width=0.3\linewidth]{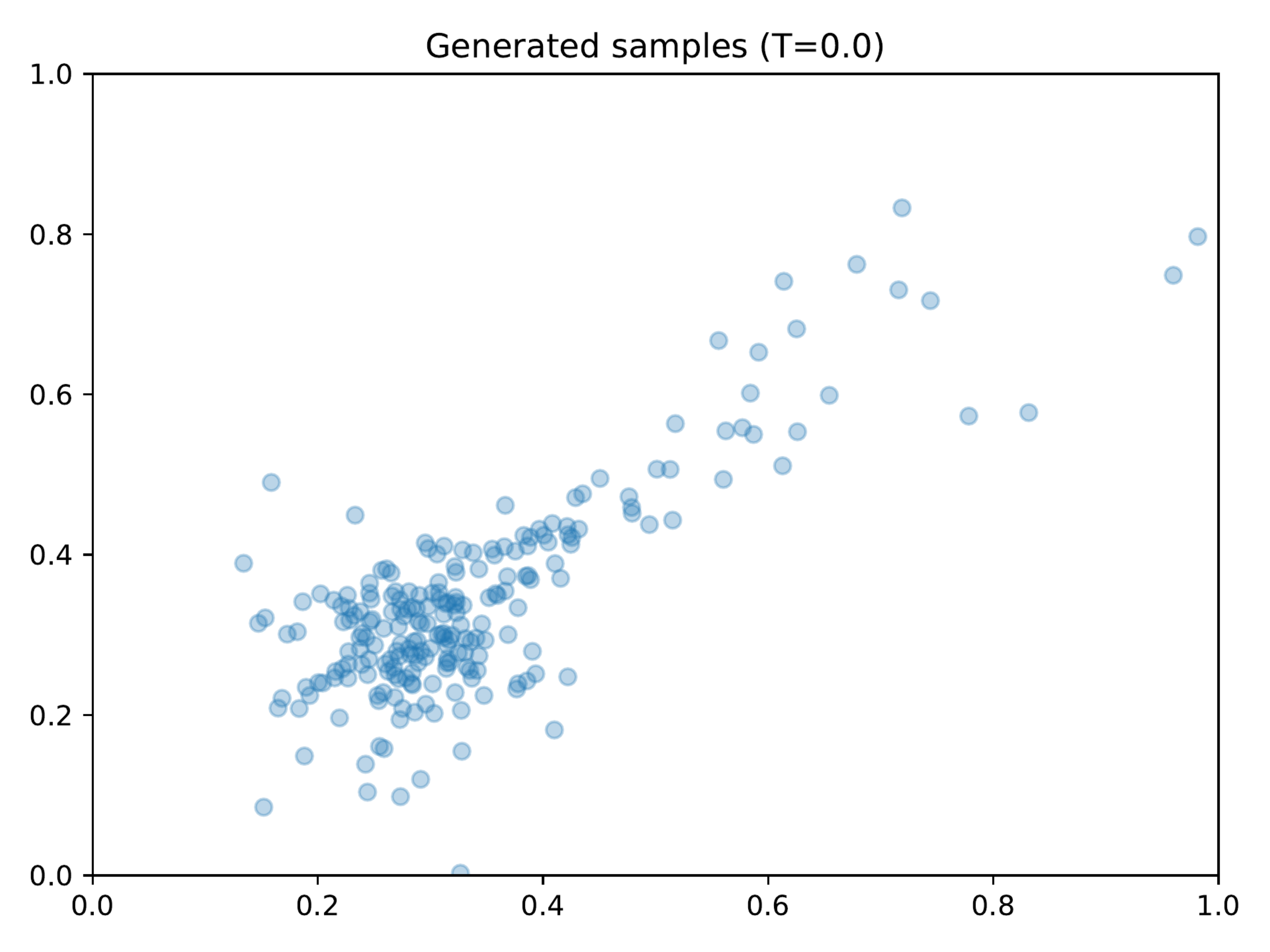}
	\includegraphics[width=0.3\linewidth]{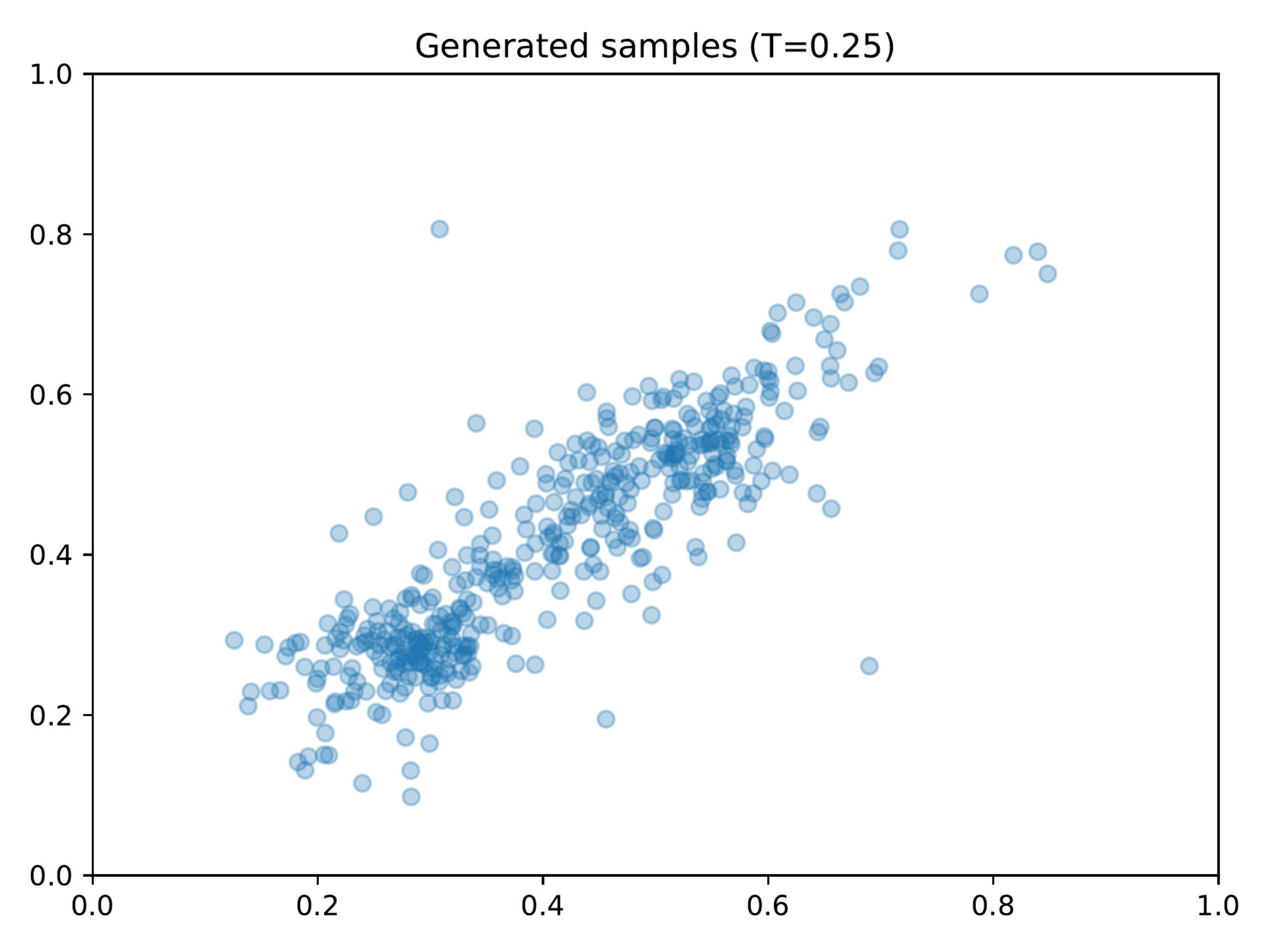}
	\includegraphics[width=0.3\linewidth]{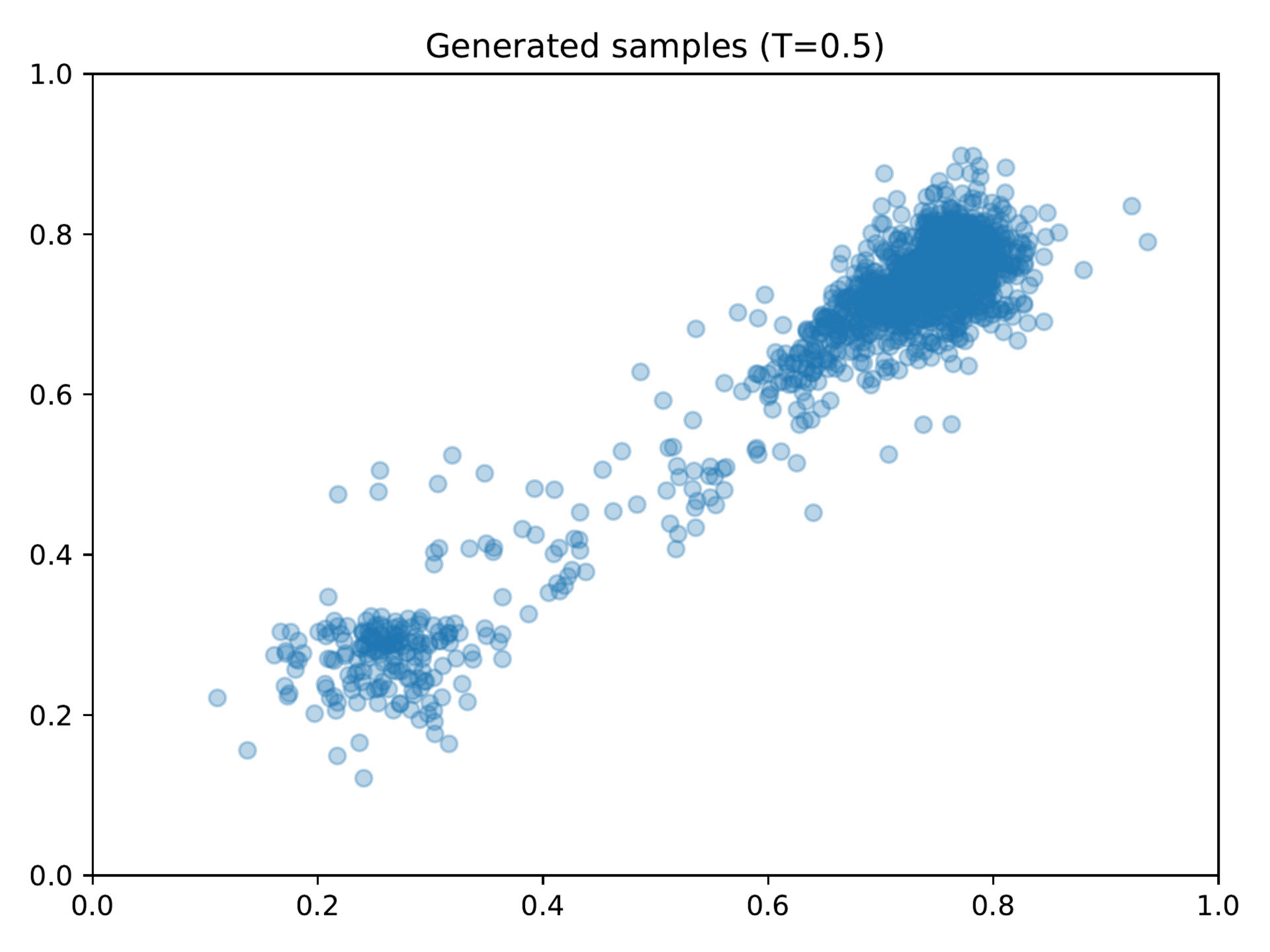}

	\includegraphics[width=0.3\linewidth]{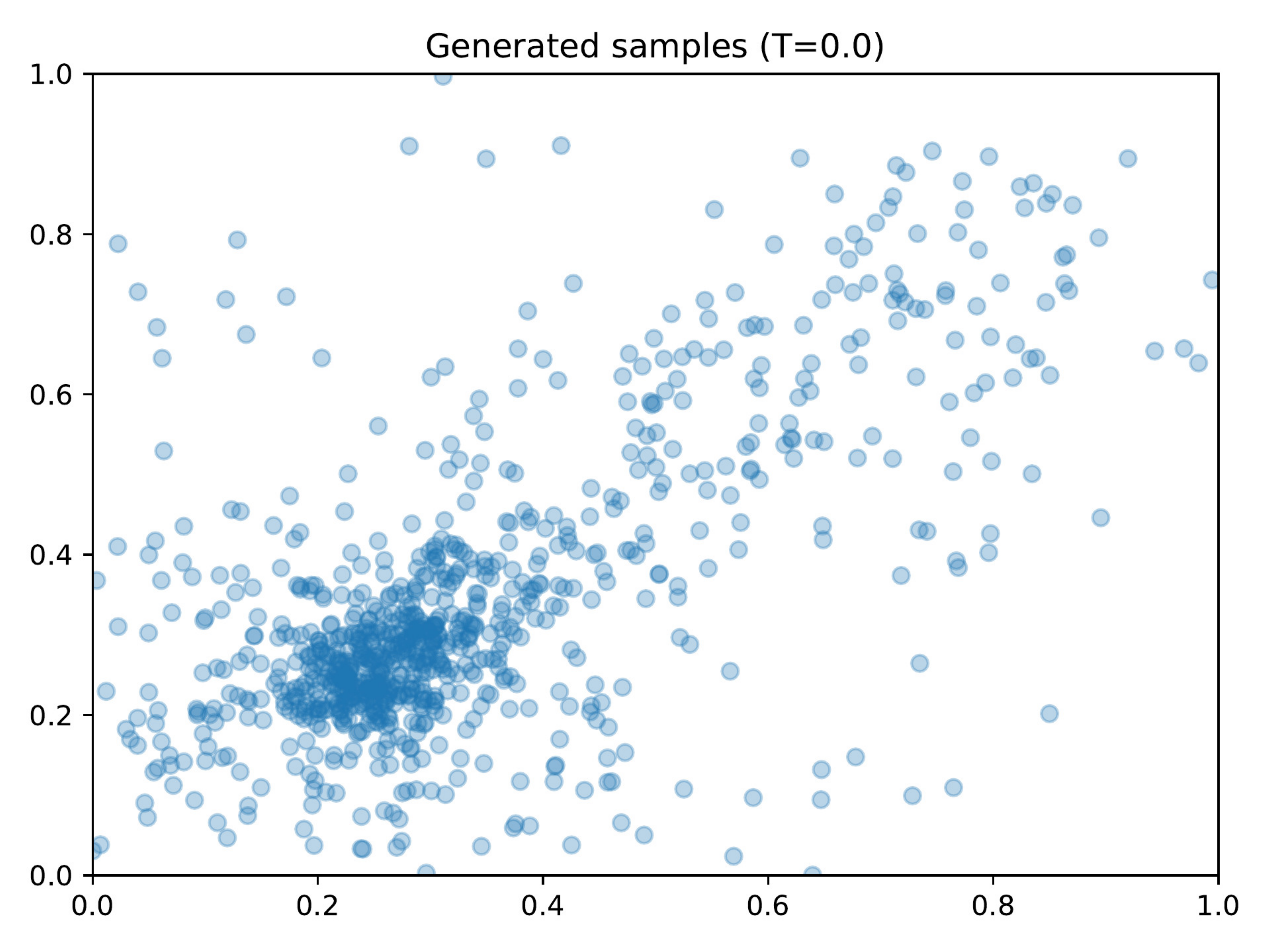}
	\includegraphics[width=0.3\linewidth]{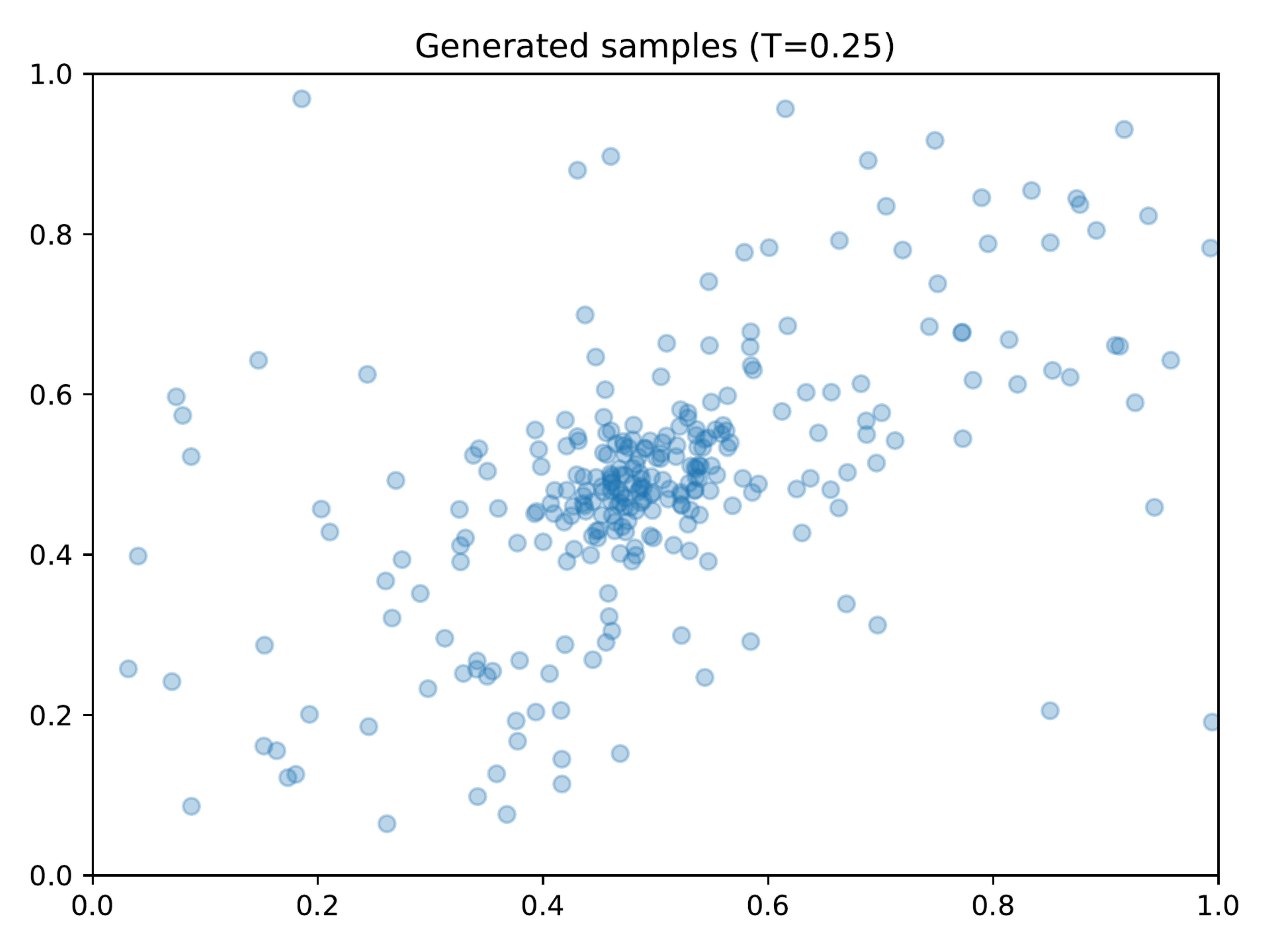}
	\includegraphics[width=0.3\linewidth]{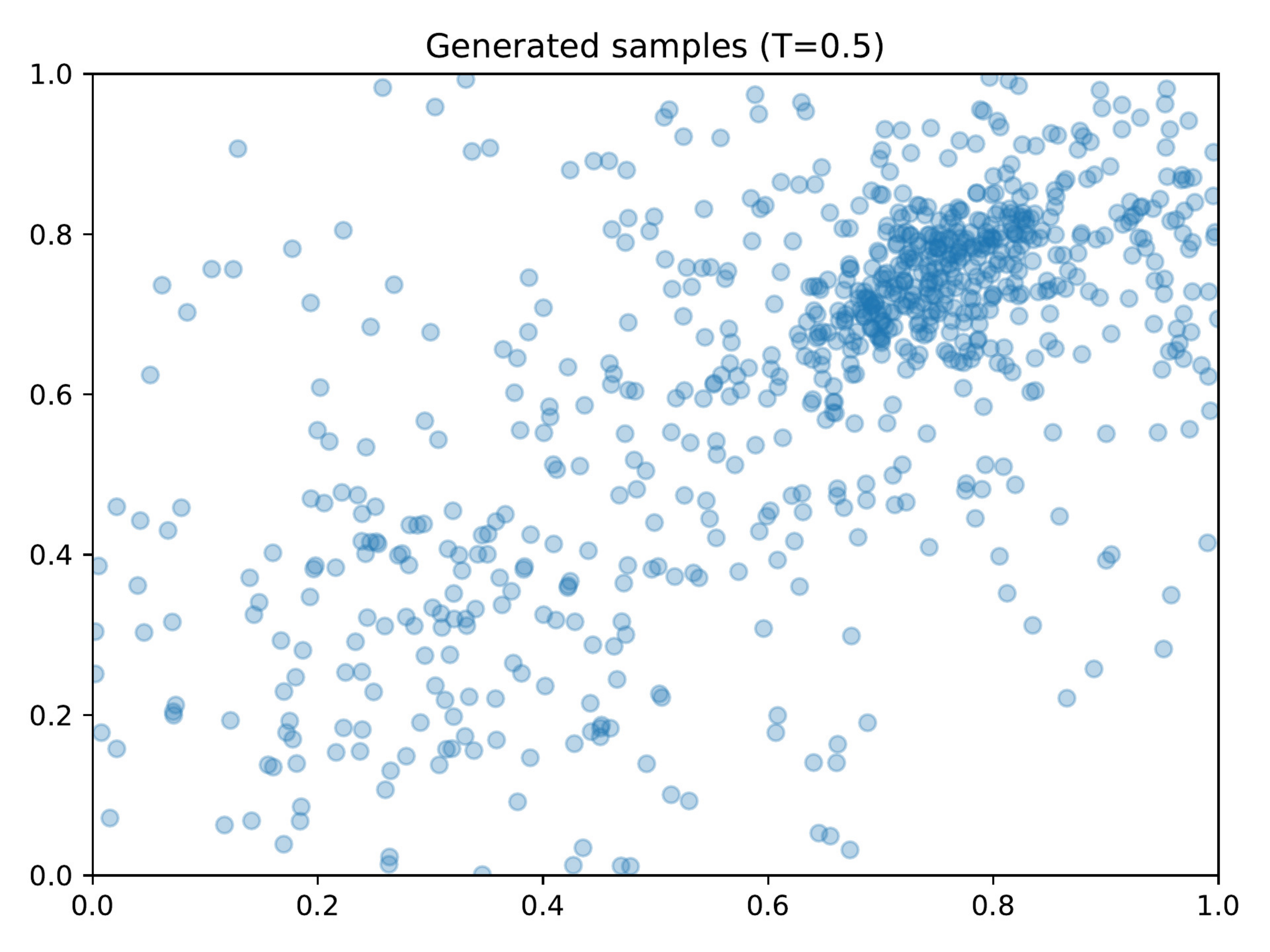}

	\caption{Low regularity test problem. From top to bottom: The generated samples of the first, and the last adaptivity iteration. ($\beta=200$, Galerkin projection)}
	\label{fig:galerkin_generated_samples_beta_200}
\end{figure}

\bibliographystyle{unsrt}
\bibliography{ref.bib}

\end{document}